\documentclass[11pt]{article}
\usepackage{amsmath,amssymb,eucal}
\usepackage[dvips]{graphicx}

\newtheorem{theorem}{\bf Theorem}[section]
\newtheorem{corollary}[theorem]{\bf Corollary}
\newtheorem{lemma}[theorem]{\bf Lemma}

\newtheorem{remark}[theorem]{\bf Remark}
\newtheorem{proposition}[theorem]{\bf Proposition}
\newtheorem{example}[theorem]{\bf Example}
\newtheorem{definition}[theorem]{\bf Definition}
\newtheorem{notation}[theorem]{\bf Notation}
\newtheorem{conjecture}[theorem]{\bf Conjecture}

\newtheorem{problem}[theorem]{\bf Problem}

\numberwithin{equation}{section}

\def\qed{$\Box$}

\begin{document}

\title{
Lens space surgeries along certain $2$-component links
related with Park's rational blow down,
and Reidemeister--Turaev torsion}
\author{Teruhisa KADOKAMI and Yuichi YAMADA}
%\date{\today}
\date{April 19, 2012}
\footnotetext[0]{%
}
\footnotetext[1]{%
The first author was supported by
a grant (No.10801021/a010402) of NSFC. 
The second author was supported by KAKENHI (Grant-in-Aid for
Scientific Research) No.21540072.
}
\maketitle
%%%%%%%%%%%%%%%%
\begin{abstract}{
We study lens space surgeries along two different families of 
2-component links, denoted by $A_{m,n}$ and $B_{p,q}$,
related with the rational homology $4$-ball used in 
J.\ Park's (generalized) rational blow down.
We determine which coefficient $r$ of the knotted component of
the link yields a lens space by Dehn surgery. 
The link $A_{m,n}$ yields a lens space
only by the known surgery with $r=mn$
and unexpectedly with $r=7$ for $(m,n)=(2,3)$.
On the other hand, $B_{p,q}$ yields a lens space by infinitely many $r$.
Our main tool for the proof is the Reidemeister--Turaev torsions, 
i.e.\ Reidemeister torsions with combinatorial Euler structures.
Our results can be extended to the links whose Alexander polynomials
are same with those of $A_{m,n}$ and $B_{p,q}$.
%We also point out that $A_{m,n}$ is closely related with
%Berge's ^^ ^^ knots in genus one fiber surface" on lens space surgery.
}\end{abstract}
%%%%%%%%%%%%%%%

\tableofcontents

%\setcounter{section}{4} 
%%%%%%%%%%%%%%%%%%%%%%%%%%%%%%%%%%%%%%%%%%
%%%                   Section 1
%%%%%%%%%%%%%%%%%%%%%%%%%%%%%%%%%%%%%%%%%%
\section{Introduction}~\label{sec:intro}
\begin{figure}[h]
\begin{center}
\includegraphics[scale=0.4]{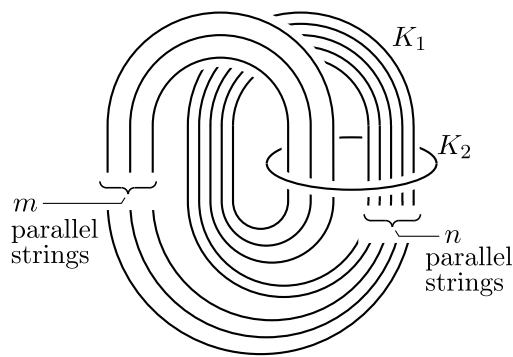}
\caption{$A_{m,n}$ \ (ex. $(m,n)=(3,5)$)}
\label{fig:Amn}
\end{center}
\end{figure}

\begin{figure}[t]
\begin{center}
\includegraphics[scale=0.4]{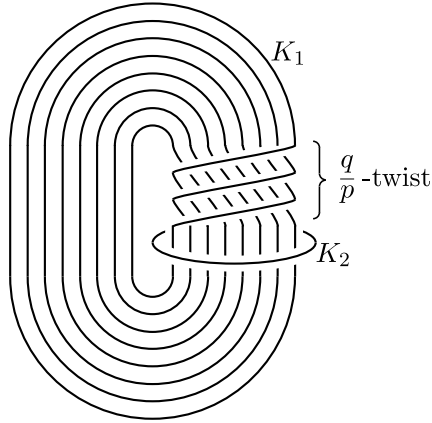}
\caption{$B_{p,q}$ \ (ex. $(p,q)=(8,3)$)}
\label{fig:Bpq}
\end{center}
\end{figure}

For a coprime pair of non-zero positive integers $(m,n)$,
let $A_{m,n}$ be a 2-component link in $S^3$ in Figure~\ref{fig:Amn},
where $K_1$ is the $(m,n)$-torus knot $T_{m,n}$ and $K_2$ is an unknot.
The linking number of $K_1$ and $K_2$ is $m+n$.
Next, for a coprime pair of non-zero integers $(p,q)$,
let $B_{p,q}$ be a 2-component link in $S^3$ in Figure~\ref{fig:Bpq},
where $K_1$ is the closure of the $(p,q)$-torus braid 
(the standard $p$-braid of the $(p,q)$-torus knot $T_{p,q}$) 
and $K_2$ is the braid axis.
The linking number of $K_1$ and $K_2$ is $p$.
%
%
%%%%%%%%%%%%%%%%%
%%%%%%%%%%%%%%%%%
%
%
%
%
%%%%%%%%%%%%%%%%%
%%%%%%%%%%%%%%%%%
%
%
For a $\mu$-component link $L = K_1 \cup K_2 \cup \cdots \cup K_\mu$,
by $(L; r_1, r_2, \ldots , r_{\mu})$, we denote
the result of $(r_1, r_2, \ldots, r_{\mu})$-surgery along $L$, where
$r_i\in \mathbb{Q}\cup \{\infty, \emptyset\}$ ($i = 1, 2, \ldots, \mu$).
A Dehn surgery along a link is called {\it a lens space surgery}
if  the result is a lens space.
We study lens space surgeries along the links 
$A_{m,n}$ and $B_{p,q}$, fixing the surgery coefficient of $K_2$ as $0$,
except in Section~\ref{sec:otherA}.
Our convention on lens spaces is
^^ ^^ $L(a, b)$ is the result of $-a/b$-surgery along the trivial knot".

The result of $(r, 0)$-surgery along a $2$-component link
$L=K_1\cup K_2$ in $S^3$ (as $\partial B^4$) such that
$K_2$ is an unknot and $r\in \mathbb{Z}$ bounds
a $4$-manifold by attaching a $1$-handle along $K_2$, and
a $2$-handle along $K_1$ with a framing $r$,
to a $0$-handle $B^4$ along $S^3$.
We denote the $4$-manifold by
$W^4(L; r, \dot{0})$ following S.~Akbulut  \cite{Ak} (see also \cite{GS}, \cite{Kir}).
Then $\pi_1(W^4(L; r, \dot{0}))\cong \mathbb{Z}/|l|\mathbb{Z}$
and $H_1(\partial W^4(L; r, \dot{0}); \mathbb{Z})\cong 
\mathbb{Z}/|l|^2\mathbb{Z}$, where $l$ is the linking number of
$K_1$ and $K_2$, and
$W^4(L; r, \dot{0})$ is a rational homology $4$-ball
if and only if $l\ne 0$.
We also note that $K_1$ can be regarded as a knot in
$S^1\times S^2$.
%
%Such a Dehn surgery $(K_1 \cup K_2; r, 0)$ is regarded as a lens space surgery 
%along a knot of $K_1$ in $S^1 \times S^2$.
%An integral lens space surgery along a knot in $S^1 \times S^2$
%is related with attaching a $4$-dimensional $2$-handle to $S^1 \times D^3$,
%thus defines a rational homology $4$-ball bounded by the resulting lens space
%of the lens space surgery. 
%

We explain a background of our targets $A_{m,n}$ and $B_{p,q}$.
J.~Park \cite{Pa} discussed {\it generalized rational blow down},
which is an operation on a $4$-manifold cutting a certain submanifold $C_{p,q}$ 
and pasting a $4$-manifold $W_{p,q}$ along $\partial C_{p,q} \cong \partial W_{p,q}$.
The $4$-manifold $W_{p,q}$ is a rational homology $4$-ball 
that is characterized by
$\pi_1 (W_{p,q}) \cong \Bbb{Z}/p\Bbb{Z}$ and 
$\partial W_{p,q}\cong L(p^2, pq-1)$ (cf. \cite{FS, CH}).
As far as the authors' knowledge, uniqueness of $W_{p,q}$ is not known well.
A lens space surgery along $B_{p,q}$ whose result is $L(p^2, pq-1)$
is often used to describe $W_{p,q}$ (cf. \cite{Li}).
On the other hand, the second author \cite{Yam3}
found a lens space surgery along $A_{m,n}$
whose result is $L(p^2, pq-1)$
and defined an algorithm to determine 
$(m,n)$ from $(p,q)$ using the Euclidean algorithm.
We also remark that the link $A_{m,n}$ appears as a bi-product of
Stipsicz--Szab\'{o}--Wahl's construction \cite[Remark~7.1]{SSW}, 
see also Endo--Mark--Horn-Morris \cite{EMM}.
We compare the links $A_{m,n}$ and $B_{p,q}$
by studying lens space surgeries along them.
Our problem is:
%
%
%%%%%%%%%%%%%%%%%%%%%%%%%%%%%%%%%%%%
\begin{problem}~\label{prob1} 
When is $(A_{m,n}; r, 0)$ ($(B_{p,q}; r, 0)$, respectively) a lens space ?
\end{problem}
%%%%%%%%%%%%%%%%%%%%%%%%%%%%%%%%%%%%
%
%
There are some trivial and overlapping cases: 
$B_{1,q}$ is the Hopf link and $A_{1,n}=B_{n+1,1}$.
The links have some symmetries: $A_{m,n}=A_{n,m}$ and
$B_{p,-q}$ is the mirror image of $B_{p,q}$.
Throughout the paper, we assume the following: 
\par \medskip \noindent
%%%%%%%%%%%%%%%%%%%%%%%%%%%%%%%%%%%%
{\bf Assumption}\ 
For $A_{m,n}$, $\gcd(m, n)=1$ and $2\le m<n$. 
For $B_{p,q}$, $\gcd(p, q)=1$, $p\ge 2$ and $q\ge 1$.
%%%%%%%%%%%%%%%%%%%%%%%%%%%%%%%%%%%%
\par \medskip

The following theorem asserts that
the link $A_{m,n}$ has at least one lens space surgery.
%which is one of the motivations of our research.
%
%
%%%%%%%%%%%%%%%%%%%%%%%%%%%%%%%%%%%%
\begin{theorem}~\label{thm:(mn)}
{\rm (\cite[Theorem 1.1]{Yam3})} \ 
For a pair $(m,n)$ satisfying the assumption,
there exists a pair $(p,q)$ satisfying the assumption
such that $(A_{m,n}; mn, 0)\cong L(p^2, pq-1)$.
\end{theorem}
%%%%%%%%%%%%%%%%%%%%%%%%%%%%%%%%%%%%
%
%
%
%
%%%%%%%%%%%%%%%%%%%%%%%%%%%%%%%%%%%%
\begin{notation}~\label{nota:modulo}
For integers $x$ and $N$, 
we denote the multiple inverse of $x$ modulo $N$ by 
$\overline{x}\ (\mathrm{mod}\ \! N)$,
i.e.\ $x\overline{x}\equiv 1\ (\mathrm{mod}\ \! N)$.
Note that, for a divisor $d$($\geq 2$) of $N$, 
both $x\ (\mathrm{mod}\ \! d)$ and 
$\overline{x}\ (\mathrm{mod}\ \! d)$
are uniquely determined by $x\ (\mathrm{mod}\ \! N)$.
We also use $\overline{x}$ as a representing integer
of $\overline{x}\ (\mathrm{mod}\ \! N)$.
\end{notation}
%%%%%%%%%%%%%%%%%%%%%%%%%%%%%%%%%%%%
%
%

Our first main theorem is the answer to Problem~\ref{prob1} 
for the link $A_{m,n}$.
%
%
%%%%%%%%%%%%%%%%%%%%%%%%%%%%%%%%%%%%
\begin{theorem}~\label{thm:MT1}
We assume $r \in \Bbb{Q}$.
\begin{enumerate}
\item[(1)]
The result of $(r, 0)$-surgery along $A_{m,n}$ is a lens space
if and only if
\begin{enumerate}
\item[(i)] $r=mn$, or 
\item[(ii)] $r=7$ for $(m, n)=(2, 3)$.
\end{enumerate}

\item[(2)]
The resulting lens spaces are as follows:
\begin{enumerate}
\item[(i)] 
$(A_{m,n}; mn, 0)\cong L((m+n)^2, m \overline{n})$,
where $n\overline{n}\equiv 1\ (\mathrm{mod}\ \! (m+n)^2)$, and 
\item[(ii)] $(A_{2,3}; 7, 0)\cong L(25, 7)$.
\end{enumerate}
\end{enumerate}
\end{theorem}
%%%%%%%%%%%%%%%%%%%%%%%%%%%%%%%%%%%%
%
%
The ``if part" of Theorem~\ref{thm:MT1} (1) (i) follows from
Theorem~\ref{thm:(mn)}.
Thus our purpose is to show the ``only if part" of (1), and (2). 
\par

\medskip

Our second main theorem is the answer to Problem~\ref{prob1} 
for the link $B_{p,q}$, which is contrast to $A_{m,n}$.
%
%
%%%%%%%%%%%%%%%%%%%%%%%%%%%%%%%%%%%%
\begin{theorem}~\label{thm:MT2}
We set $\alpha/\beta \in \Bbb{Q}$, 
where $\alpha$ and $\beta$ are coprime integers.
\begin{enumerate}
\item[(1)]
The result of $(\alpha/\beta, 0)$-surgery along $B_{p,q}$ is a lens space
if and only if $|\alpha-pq\beta|=1$.

\item[(2)]
For $\alpha/\beta$ with $|\alpha-pq\beta|=1$, the resulting lens space is
$(B_{p,q}; \alpha/\beta, 0)\cong L(p^2\beta, \alpha)$.

\end{enumerate}
\end{theorem}
%%%%%%%%%%%%%%%%%%%%%%%%%%%%%%%%%%%%
%
%
%%%%%%%%%%%%%%%%%%%%%%%%%%%%%%%%%%%%
\begin{remark}~\label{rmk:}
{\rm 
We remark on surgeries along the mirror images of the links.
Naturally, we have $B_{-p,-q}=B_{p,q}, B_{-p,q}=B_{p,-q}$
as unoriented links, and $B_{p,-q}$ is the mirror image of $B_{p,q}$.
Theorem~\ref{thm:MT2} can be extended to the cases $p<0$ or/and $q<0$.
Similarly, Theorem~\ref{thm:(mn)} can be extended to the mirror image
$A_{m,n}!$ of $A_{m,n}$. 
We note that $A_{m,n}!$ is not included in the
family $\{ A_{m,n} \}$.
}
%By Theorem~\ref{thm:(mn)}, the result of Dehn surgery
%$(A_{m,n}!; r,0)$ is a lens space if and only if
%$r =-mn$, or $r=-7$ for the case $(m,n)=(2,3)$.
%But we have not define $A_{m,n}$ for the case $m<0$ or/and $n<0$.
%(cf. \cite{Yam4}). 
%
%But we have not define The reader may be able to extend the theorem
%on $A_{m,n}$ for the case $m<0$ or/and $n<0$,
%Contrastingly, the links $A_{m,n}$ and their mirror images $A_{m,n}!$ 
%are kind of different families of links. 
\end{remark}
%%%%%%%%%%%%%%%%%%%%%%%%%%%%%%%%%%%%
%
%

We list some corollaries without the proofs.
%
%
%%%%%%%%%%%%%%%%%%%%%%%%%%%%%%%%%%%%
\begin{corollary}~\label{cor:L257}
The lens space $(A_{2,3}; 7, 0) \cong L(25, 7)$ cannot be obtained
by any $(r, 0)$-surgery along $B_{p,q}$.
\end{corollary}
%%%%%%%%%%%%%%%%%%%%%%%%%%%%%%%%%%%%
%
%
%%%%%%%%%%%%%%%%%%%%%%%%%%%%%%%%%%%%
\begin{corollary}[Integral lens space surgery along $B_{p,q}$]~\label{cor:Btype}
Suppose $\beta=1$ in Theorem~\ref{thm:MT2}.
Then $(B_{p,q}; \alpha, 0)$ is a lens space if and only if
$\alpha=pq-1$ or $pq+1$.
The resulting lens space is 
$L(p^2, pq-1)$ or $L(p^2, pq+1)$, respectively.
\end{corollary}
%%%%%%%%%%%%%%%%%%%%%%%%%%%%%%%%%%%%
%
%
%For a link $L= K_1 \cup K_2$ such that the component $K_2$ is an unknot
%and the linking number $l$ of $K_1$ and $K_2$ satisfies
%$\vert l\vert >1$, an integral $(r,0)$-surgery along $L$ defines 
%a rational homology $4$-ball $W^4(L; r ,\dot{0})$
%bounded by the $3$-manifold $(L; r,0)$, where the dotted $0$ means a $4$-dimensional $1$-handle.
%This method is introduced by S.~Akbulut \cite{Ak},
%see also \cite{GS} for the method ^^ ^^ $4$-dimensional Kirby diagrams".
%The framed knot $K_1$ is regarded as an attaching part of a $2$-handle. 
%It holds that $\pi_1(W^4(L; r ,\dot{0})) \cong \mathbb{Z}/l\mathbb{Z}$.
%Note that, if the framed knot $K_1$ moves isotopically in 
%$S^1 \times S^2$ (the result of $0$-surgery along $K_2$), 
%then the rational homology $4$-ball is diffeomorphically unchanged from $W^4(L; r ,\dot{0})$,
%since such a move corresponds to a handle slide of the $2$-handle over the $1$-handle.
%
%The $(pq-1,0)$-surgery along $B_{p,q}$
%%lens space surgery $(B_{p,q};pq-1,0)$
%in Corollary~\ref{cor:Btype}
%defines a rational homology $4$-ball $W^4(B_{p,q}; pq-1,\dot{0})$.
%On the other hand, under the correspondence $(m,n)=A(p-q, q)$,
%%$(mn,0)$-surgery along $A_{m,n}$
%the lens space surgery in Theorem~\ref{thm:(mn)}
%defines a rational homology $4$-ball $W^4(A_{m,n}; mn,\dot{0})$.

By Theorem~\ref{thm:(mn)} and Corollary~\ref{cor:Btype},
both $W^4(A_{m,n}; mn, \dot{0})$ and 
$W^4(B_{p,q}; pq-1, \dot{0})$ represent $W_{p,q}$,
under the correspondence between $(m,n)$ and $(p,q)$
in Theorem~\ref{thm:(mn)}.
Our second problem is:
%
%
%%%%%%%%%%%%%%%%%%%%%%%%%%%%%%%%%%%%
\begin{problem}~\label{prob:4ball}
Are $W^4(A_{m,n}; mn,\dot{0})$ and $W^4(B_{p,q}; pq-1,\dot{0})$
diffeomorphic, homeomorphic or homotopic 
relative to the boundaries?
%
%Does there exist a homotopy equivalent map whose 
%restriction over the boundary is a homeomorphism?
%
%Are they diffeomorphic or homeomorphic to each other?
\end{problem}
%%%%%%%%%%%%%%%%%%%%%%%%%%%%%%%%%%%%
%
%

The lens space $(A_{2,3};7,0)$ in Corollary~\ref{cor:L257} satisfies $L(25, 7)\cong -L(25, 7)$. 
On the other hand, on the integral lens spaces surgeries along $B_{p,q}$ in Corollary~\ref{cor:Btype},
it is easy to see that $L(p^2, pq-1) \cong - L(p^2, pq+1)$ and that $L(p^2, pq-1)\not \cong L(p^2, pq+1)$.
Thus we have 
%
%
%%%%%%%%%%%%%%%%%%%%%%%%%%%%%%%%%%%%
\begin{corollary}~\label{cor:Nopair}
Only in this corollary, we regard the link $A_{m,n}$ or $B_{p,q}$ 
as a knot $K_1$ in $S^1 \times S^2$, the result of $0$-surgery along $K_2$.
We assume that $2\le m<n, 2\le p$ and $1\le q$.
Then, any knot $A_{m, n}$ is not isotopic to any knot $B_{p,q}$.
\end{corollary}
%%%%%%%%%%%%%%%%%%%%%%%%%%%%%%%%%%%%
%
%
Corollary~\ref{cor:Nopair} asserts that 
the attaching parts $K_1$ of the $2$-handles 
of $W^4(A_{m,n}; mn, \dot{0})$ and $W^4(B_{p,q}; pq-1,\dot{0})$
are not isotopic in $S^1 \times S^2$. 
Thus, the handle decompositions of the rational homology $4$-balls 
do not move to each other by only handle slides of the $2$-handles over the $1$-handles.
Problem~\ref{prob:4ball} may be still open.
\par

\medskip

Our results can be regarded as lens space surgeries along 2-component links.
For usual lens space surgery along knots, see 
\cite{Ba}, 
\cite{Ber},
\cite{CGLS},
\cite{Go}
and so on. 
We point out that Theorem~\ref{thm:MT1} (on $A_{m,n}$) can also be obtained 
by the results
^^ ^^ $(A_{m,n}; \emptyset, 0)$ is a hyperbolic manifold"
in \cite{DMM1}, and the Cyclic Surgery Theorem \cite{CGLS}.
If we use them, then the proof of the theorem can be
shortened (see Subsection~\ref{ssec:CGLS}).
The reason why we do not use them is 
to clarify effectivity of Alexander polynomials and 
Reidemeister--Turaev torsions from technical point of view.
As consequences, they preserved information of lens space 
surgeries completely in the present case, and our results
are generalized to wider situations, see Theorem~\ref{thm:extMT1} and 
Theorem~\ref{thm:extMT2}.
The links $A_{m,n}$ are related to
subfamilies of knots, called TypeVII and Type VIII
in Berge's list \cite{Ber},
see 
\cite{Kad4}, \cite{Yam1} and \cite{Yam2}.
On the other hand, the results $(B_{p,q}; r,0)$ with any $r$ is a Seifert manifold
(or a graph manifold), thus Theorem~\ref{thm:MT2} looks
like L.~Moser's result \cite{Mos} on lens space surgeries along torus knots.

\medskip

We also study a generalization of Theorem~\ref{thm:MT1}.
%Contrasted to the case of $B_{p,q}$ in Subsection~\ref{ssec:otherB},
To determine all Dehn surgeries along $A_{m,n}$
is a hard problem.
As far as the authors' knowledge,
the complete answer is not given.
From our present results and some known results,
we feel like to raise the following conjecture:
%
%
%%%%%%%%%%%%%%%%%%%%%%%%%%%%%%%%%%%%
\begin{conjecture}\label{conj:Amn}
Let $M=(A_{m,n}; \alpha_1/\beta_1, \alpha_2/\beta_2)$
be the result of 
$(\alpha_1/\beta_1, \alpha_2/\beta_2)$-surgery
along $A_{m,n}$, where 
$\alpha_i$ and $\beta_i$ $(i=1, 2)$ are coprime integers 
with $\beta_i\ge 1$.
Then $M$ is a lens space if and only if
(1) $\alpha_1/\beta_1=mn$ and $\beta_2=1$, or
(2) $\alpha_1/\beta_1=7$ for $(m, n)=(2, 3)$.
\end{conjecture}
%%%%%%%%%%%%%%%%%%%%%%%%%%%%%%%%%%%%
%
%
K.~Ichihara (\cite{IS}) informed to the authors that,
if we fix $\alpha_1/\beta_1=mn$, then $M$
is a lens space for any integer ($\beta_2=1$).
He says that it can be shown by a method in \cite{Yam3}
(i.e.\ a geometric method).
In Section~\ref{sec:otherA}, we compute the Reidemeister torsions of $M$
under the case (1) $\alpha_1/\beta_1=mn$ and (2) $\alpha_1/\beta_1=7$ for $(m, n)=(2, 3)$, respectively.
%By Kirby calculus, ^^ ^^ if part" of Conjecture~\ref{conj:Amn} is settled,
A partial affirmative answer for ^^ ^^ only if part" of the case (1) is given.
%In Section~\ref{sec:otherA}, we compute the Reidemeister torsions of $M$
%under the case (1) $\alpha_1/\beta_1=mn$ and (2) $\alpha_1/\beta_1=7$ for $(m, n)=(2, 3)$, respectively.
%The result support Conjecture \ref{cj:Amn} above.

\medskip

In Section~\ref{sec:torsion}, we explain about the Reidemeister torsion
and its basic properties such as surgery formulae, $d$-norm
and combinatorial Euler structure (i.e.\ Reidemeister--Turaev torsion).
In Section~\ref{sec:Alex}, we compute the Alexander polynomial of $A_{m,n}$.
In Section~\ref{sec:key}, 
we compute the Reidemeister torsions of $(A_{m,n}; r, 0)$
by using the results in Section~\ref{sec:torsion} and Section \ref{sec:Alex}.
In Section~\ref{sec:proofA}, 
we prove the ``only if part" of Theorem \ref{thm:MT1} (1)
by using the Reidemeister--Turaev torsions, and
Theorem~\ref{thm:MT1} (2) 
by using the values of the Reidemeister torsions,
Theorem~\ref{thm:(mn)}, and Kirby moves.
In Section~\ref{sec:Bpq}, we prove Theorem \ref{thm:MT2} 
by using Seifert structures of the link complement
and the Reidemeister torsions.
In Section~\ref{sec:otherA}, we study some 
lens space surgeries along $A_{m,n}$ other than $(r, 0)$-surgery, related to
Conjecture~\ref{conj:Amn}.
%In Section~\ref{sec:lensS}, we point out that $A_{m,n}$ is closely related with
%Berge knots in genus one fiber surface (cf. \cite{Ber}).
In Section~\ref{sec:rem}, we give an alternative proof of
Theorem \ref{thm:MT1} (1) under some assumption, 
and we generalize Theorem \ref{thm:MT1}
and Theorem \ref{thm:MT2}.

%%%%%%%%%%%%%%%%%%%%%%%%%%%%%%%%%%%%%%%%%%
%%%                   Section 2
%%%%%%%%%%%%%%%%%%%%%%%%%%%%%%%%%%%%%%%%%%
\section{Reidemeister torsion}~\label{sec:torsion}
Our method to prove the main theorems 
is to deduce necessary conditions from 
the Reidemeister torsions of both
the surgered manifolds and lens spaces.
In this section, we state surgery formulae
of the Reidemeister torsions (Subsection~\ref{ssec:surgery}), 
and define derived invariants: 
one is the {\it $d$-norm} (Subsection~\ref{ssec:norm}), and
the other is the {\it Reidemeister--Turaev torsion}
which is a lift of the Reidemeister torsions
by fixing a combinatorial Euler structure
(Subsection~\ref{ssec:Euler}).
%To prove the main theorem,
%we will compare the Reidemeister torsion
%of the result $M$ of the Dehn surgery along the link 
%and that of a lens space $L(p,q)$.
%%
%Some necessary conditions are obtained 
%from the value $\tau^{\psi_d}(M)$ of the Reidemeister torsion
%in the $d$-th cyclotomic field $\Bbb{Q}(\zeta_d)$
%by $d$-norm, where $d$($\geq 2$) is a divisor of $p$.
%%
%From the sequence of the equalities on $\tau^{\psi_d}(M)$ in $\Bbb{Q}(\zeta_d )$
%for all divisors $d$ of $p$
%(with a fixed combinatorial Euler structure of $M$),
%we take an identity on symmetric Laurent polynomials, as a lift of the equalities.
%Finally, we regard the identity as an equation of the surgery coefficient
%for $M$ to be a lens space.
%
%In this section,
%we review 
%definitions and some formulae of the Reidemeister torsion 
%in Subsection~\ref{ssec:surgery},
%$d$-norms in the $d$-th cyclotomic field
%in Subsection~\ref{ssec:norm}.
%In Subsection~\ref{ssec:Euler}, 
%we study on a certain uniqueness of 
%a symmetric polynomial
%as a lift of the sequence of the equalities in $\Bbb{Q}(\zeta_d)$,
%respectively.

%%%%%%%%%%%%%%%%%%%%%%%%%%%%%%%%%%%%%%%%%%
%%%                   SubSection 2-1
%%%%%%%%%%%%%%%%%%%%%%%%%%%%%%%%%%%%%%%%%%
\subsection{Surgery formulae}~\label{ssec:surgery}
For a precise definition of the Reidemeister torsion,
the reader refer to V.~Turaev \cite{Tur1, Tur2}. 
Throughout this paper, we use the following notations.
%
%
%%%%%%%%%%%%%%%%%%%%%%%%%%%%%%%%%%%%
\begin{notation}~\label{nota:mfd} 
{\rm (for manifolds and homologies)}\\
Let  $L=K_1\cup \cdots \cup K_{\mu}$
be an oriented $\mu$-component link in a homology $3$-sphere.
\begin{center}
\begin{tabular}{ll}
$E_L$ & the complement of $L$. \\
$m_i, l_i$ & a meridian and a longitude of the $i$-th component $K_i$. \\
$[m_i], [l_i]$ & their homology classes. \\
${\it \Delta}_L(t_1, \ldots, t_{\mu})$ & the Alexander polynomial of $L$, 
where $t_i$ is represented by $m_i$. \\
$(L ; r_1, \ldots, r_{\mu})$ &
the result of $(r_1, \ldots, r_{\mu})$-surgery along $L$, \\
& \quad where $r_i\in \mathbb{Q}\cup \{\infty, \emptyset \}$ 
is the surgery coefficient of $K_i$. \\
$V_i$ &
the solid torus attached along $K_i$ in the Dehn surgery.\\
$m_i'$, $[m_i']$ &
a meridian of $V_i$, and its homology class.\\
$l_i'$, $[l_i']$ &
an oriented core curve of $V_i$, and its homology class.
\end{tabular}
\end{center}
\end{notation}
%%%%%%%%%%%%%%%%%%%%%%%%%%%%%%%%%%%%
%
%

Let $X$ be a finite CW complex
and $\pi : {\tilde X}\to X$ its maximal abelian covering.
Then ${\tilde X}$ has a CW structure induced by that of $X$ and $\pi$,
and the cell chain complex $\mathbf{C}_{\ast}$ of ${\tilde X}$
has a $\mathbb{Z}[H]$-module structure, where
$H=H_1(X ; \mathbb{Z})$ is the first homology of $X$.
For an integral domain $R$ and a ring homomorphism
$\psi : \mathbb{Z}[H]\to R$,  
\lq\lq the chain complex of ${\tilde X}$ related with $\psi$\rq\rq,
denoted by $\mathbf{C}_{\ast}^{\psi}$,
is $\mathbf{C}_{\ast}\otimes_{\mathbb{Z}[H]}Q(R)$,
where $Q(R)$ is the quotient field of $R$.
{\it The Reidemeister torsion of $X$ related with $\psi$}, 
denoted by $\tau^{\psi}(X)$,
is calculated from $\mathbf{C}_{\ast}^{\psi}$, and 
is an element of $Q(R)$ determined up to multiplication of $\pm \psi(h)\ (h\in H)$.
If $R=\mathbb{Z}[H]$ and $\psi$ is the identity map,
then we denote $\tau^{\psi}(X)$ by $\tau(X)$.
We note that $\tau^{\psi}(X)$ is not zero if and only if 
$\mathbf{C}_{\ast}^{\psi}$ is acyclic.
%
%
%%%%%%%%%%%%%%%%%%%%%%%%%%%%%%%%%%%%
\begin{notation}~\label{nota:algebra}
{\rm (for algebra)}\\
For a pair of elements $A, B$ in $Q(R)$, if there exists an element $h\in H$
such that $A=\pm \psi(h)B$, then we denote the equality by $A\doteq B$.
We will often take a field $F$
and a ring homomorphism $\psi : \mathbb{Z}[H_1(M)]\to F$.
We mainly use the $d$-th cyclotomic fields $\mathbb{Q}(\zeta_d)$ as $F$,
where $\zeta_d$ is a primitive $d$-th root of unity.
\end{notation}
%%%%%%%%%%%%%%%%%%%%%%%%%%%%%%%%%%%%
%
%

For the first lemma, we need a little general setting: 
Let $E$ be a compact $3$-manifold whose boundary $\partial E$ 
consists of tori.
We study the $3$-manifold $M=E\cup V_1\cup \cdots \cup V_n$
obtained by attaching solid tori $V_i$ to $E$ by attaching maps
$f_i : \partial V_i\to \partial E$
($\mathrm{Im}\ \! (f_i)\cap \mathrm{Im}\ \! (f_j)=\emptyset$ for $i\ne j$).
We let $\iota : E\hookrightarrow M$ denote the natural inclusion.
%
%
%%%%%%%%%%%%%%%%%%%%%%%%%%%%%%%%%%%%
\begin{lemma}[Surgery formula I]~\label{lem:surgery1}
% Turaev 70s?
If $\psi([l_i'])\ne 1$ for every $i=1, \ldots, n$, then we have
\[
\tau^{\psi}(M)\doteq \tau^{\psi'}(E)
\prod_{i=1}^n(\psi([l_i'])-1)^{-1},
\]
where $\psi'=\psi \circ \iota_{\ast}$ 
($\iota_{\ast}$ is a ring homomorphism induced by $\iota$).
\end{lemma}
%%%%%%%%%%%%%%%%%%%%%%%%%%%%%%%%%%%%
%
%
For the case of the link complement of a homology $3$-sphere, 
the Reidemeister torsion is closely related with the Alexander polynomial.
%
%
%%%%%%%%%%%%%%%%%%%%%%%%%%%%%%%%%%%%
\begin{lemma}~\label{lem:Alexander}{\rm (Milnor \cite{Mil})} \ 
For a $\mu$-component link 
$L=K_1\cup \cdots \cup K_{\mu}$ in a homology $3$-sphere, 
we have
\[
\tau(E_L)\doteq
\begin{cases}
\Delta_L(t_1)(t_1-1)^{-1} & (\mu =1),\\
\Delta_L(t_1, \ldots, t_{\mu}) &  (\mu\ge 2).
\end{cases}
%\left\{
%\begin{array}{cl}
%\Delta_L(t_1)(t_1-1)^{-1} & \quad (\mu =1),\\
%\Delta_L(t_1, \ldots, t_{\mu}) & \quad (\mu\ge 2).
%\end{array}
%\right.
\]
\end{lemma}
%%%%%%%%%%%%%%%%%%%%%%%%%%%%%%%%%%%%
%
%
By Lemma 2.1 and Lemma 2.2, 
we have the following:
%
%
%%%%%%%%%%%%%%%%%%%%%%%%%%%%%%%%%%%%
\begin{lemma}[Surgery formula II]~\label{lem:surgery2}
{\rm (T.~Sakai \cite{Sa}, V.~G.~Turaev \cite{Tur1})}
\begin{enumerate}
\item[(1)]
Let $K$ be a knot in a homology $3$-sphere.
We set $M=(K ; p/q)\ (|p|\ge 2)$, 
where $p$ and $q$ are coprime integers.
Let $T$ be a generator of $H_1(M)$ represented by a meridian of $K$, 
and $\psi_d : \mathbb{Z}[H_1(M)]\to \mathbb{Q}(\zeta_d)$
a ring homomorphism defined by $\psi_d(T)=\zeta_d$,
where $d$ ($\ge 2$) is a divisor  of $p$.
Then we have
%In the case $M=(K ; p/q)\ (|p|\ge 2)$, 
%we have $H=H_1(M)\cong \langle T\ \vert \ T^p=1\rangle 
%\cong \mathbb{Z}/|p|\mathbb{Z}$,
%where $T$ is represented by the meridian $[m]$.
%For a divisor $d$ ($\ge 2$) of $p$, we define a ring homomorphism 
%$\psi_d : \mathbb{Z}[H]\to \mathbb{Q}(\zeta_d)$ by $\psi_d(T)=\zeta_d$.
%Then we have
\[
\tau^{\psi_d}(M)\doteq {\mit \Delta}_K(\zeta_d)
(\zeta_d-1)^{-1}(\zeta_d^{{\bar q}}-1)^{-1}
\]
where $q\overline{q}\equiv 1\ (\mathrm{mod}\ \! p)$.

\item[(2)]
Let $L$ be a $\mu$-component link in a homology $3$-sphere.
We set $M=(L ; p_1/q_1, \ldots, p_{\mu}/q_{\mu})\ (\mu \ge 2)$,
where $p_i$ and $q_i$ are coprime integers for every $i=1, \ldots, \mu$.
Let $F$ be a field, 
$\psi : \mathbb{Z}[H_1(M)]\to F$ a ring homomorphism
with $\psi([m_i]^{r_i}[l_i]^{s_i})\ne 1$ for every $i=1, \ldots, \mu$,
where $r_i$ and $s_i$ are integers satisfying $p_i s_i -q_i r_i = -1$.
Then we have
%In the case $M=(L ; p_1/q_1, \ldots, p_{\mu}/q_{\mu})\ (\mu \ge 2)$.
%Let $F$ be a field and $\psi : \mathbb{Z}[H_1(M)]\to F$ a ring homomorphism.
%If $\psi([m_i]^{r_i}[l_i]^{s_i})\ne 1$ for every $i=1, \ldots, \mu$, then we have
\[
\tau^{\psi}(M)\doteq {\it \Delta}_L(\psi([m_1]), \ldots, \psi([m_{\mu}]))
\prod_{i=1}^{\mu}(\psi([m_i]^{r_i}[l_i]^{s_i})-1)^{-1}.
\]
\end{enumerate}
\end{lemma}
%%%%%%%%%%%%%%%%%%%%%%%%%%%%%%%%%%%%
%
%
%
%
%%%%%%%%%%%%%%%%%%%%%%%%%%%%%%%%%%%%
\begin{example}~\label{ex:lens}
The lens space $L(p,q)$ is obtained as $-p/q$-surgery along the unknot.
By Lemma~\ref{lem:surgery2} (1), for a divisor $d\ge 2$ of $p$, we have
\[
\tau^{\psi_d}(L(p, q))\doteq (\zeta_d-1)^{-1}(\zeta_d^{{\bar q}}-1)^{-1},
\]
where $q\overline{q}\equiv 1\ (\mathrm{mod}\ \! p)$.
\end{example}
%%%%%%%%%%%%%%%%%%%%%%%%%%%%%%%%%%%%
%
%

%%%%%%%%%%%%%%%%%%%%%%%%%%%%%%%%%%%%
%%%%%%%%%%%%%%%%%%%%%%%%%%%%%%%%%%%%
We recall the Torres formula for the Alexander polynomials.
%
%
%%%%%%%%%%%%%%%%%%%%%%%%%%%%%%%%%%%%
\begin{lemma}[Torres formula]~\label{lem:Torres}{\rm (\cite{Tor})} \
Let $L=K_1\cup \cdots \cup K_{\mu} \cup K_{\mu+1}\ 
(\mu \ge 1)$ be an oriented $(\mu+1)$-component link,
$L'=K_1\cup \cdots \cup K_{\mu}$ a $\mu$-component sublink,
and $\ell_i=\mathrm{lk}\ \! (K_i, K_{\mu+1})\ (i=1, \ldots, \mu)$.
Then we have\[
\Delta_L(t_1, \ldots, t_{\mu}, 1)\doteq
\begin{cases}
{\displaystyle \frac{t_1^{\ell}-1}{t_1-1}{\mit \Delta}_{L'}(t_1)} & (\mu=1), \\
(t_1^{\ell_1}\cdots t_{\mu}^{\ell_{\mu}}-1)
{\mit \Delta}_{L'}(t_1, \ldots, t_{\mu}) & (\mu \ge 2).
\end{cases}
%\left\{
%\begin{array}{cl}
%{\displaystyle \frac{t_1^{\ell}-1}{t_1-1}{\mit \Delta}_{L'}(t_1)} & (\mu=1),
%\medskip\\
%(t_1^{\ell_1}\cdots t_{\mu}^{\ell_{\mu}}-1)
%{\mit \Delta}_{L'}(t_1, \ldots, t_{\mu}) & (\mu \ge 2).
%\end{array}
%\right.
\]
\end{lemma}

%%%%%%%%%%%%%%%%%%%%%%%%%%%%%%%%%%%%%%%%%%
%%%                   SubSection 2-2
%%%%%%%%%%%%%%%%%%%%%%%%%%%%%%%%%%%%%%%%%%
\subsection{$d$-norm}~\label{ssec:norm}
About algebraic fields,
the reader refer to L.~C.~Washington \cite{Was} for example.
\par

\medskip

For an element $x$ in the $d$-th cyclotomic field $\mathbb{Q}(\zeta_d)$,
the {\it $d$-norm} of $x$ is defined as 
\[
N_d(x)=\prod_{\sigma \in \mathrm{Gal}\ \! 
(\mathbb{Q}(\zeta_d)/\mathbb{Q})}
\sigma(x),
\]
where $\mathrm{Gal}\ \! (\mathbb{Q}(\zeta_d)/\mathbb{Q})$
is the Galois group related with a Galois extension
$\mathbb{Q}(\zeta_d)$ over $\mathbb{Q}$.
The following is well-known.
%
%
%%%%%%%%%%%%%%%%%%%%%%%%%%%%%%%%%%%%
\begin{proposition}~\label{pr:norm}
\begin{enumerate}
\item[(1)]
If $x\in \mathbb{Q}(\zeta_d)$, then $N_d(x)\in \mathbb{Q}$.
The map $N_d : \mathbb{Q}(\zeta_d)\setminus \{0\}\to 
\mathbb{Q}\setminus \{0\}$ is a group homomorphism.

\item[(2)]
If $x\in \mathbb{Z}[\zeta_d]$, then $N_d(x)\in \mathbb{Z}$.
\end{enumerate}
\end{proposition}
%%%%%%%%%%%%%%%%%%%%%%%%%%%%%%%%%%%%
%
%

By easy calculations, we have the following.
%
%
%%%%%%%%%%%%%%%%%%%%%%%%%%%%%%%%%%%%
\begin{lemma}~\label{lem:cyclotomic}
\begin{enumerate}
\item[(1)]
${\displaystyle
N_d(\pm \zeta_d)=
\begin{cases}
\pm 1 & (d=2),\\
1 & (d\ge 3).
\end{cases}
}$

\item[(2)]
${\displaystyle
N_d(1-\zeta_d)=
\begin{cases}
\ell & (\mbox{$d$ is a power of a prime $\ell \ge 2$}),\\
1 & (\mbox{otherwise}).
\end{cases}
}$
\end{enumerate}
\end{lemma}
%%%%%%%%%%%%%%%%%%%%%%%%%%%%%%%%%%%%
%
%

About applications of $d$-norms,
for example, see \cite{Kad1, Kad2, Kad3, KMS, KY1, KY2}.
\par 

\medskip

W.~Franz \cite{Fz} showed the following, and
classified lens spaces by using it.
We state a modified version (cf. \cite{KY1}).
%
%
%%%%%%%%%%%%%%%%%%%%%%%%%%%%%%%%%%%%
\begin{lemma}~\label{lem:Franz}{\rm (Franz \cite{Fz})}
Let $p\ge 2$ be an integer, and $(\mathbb{Z}/p\mathbb{Z})^{\times}$
the multiplicative group of a ring $\mathbb{Z}/p\mathbb{Z}$.
For $a_i, b_i\in (\mathbb{Z}/p\mathbb{Z})^{\times}$\ $(i=1, \ldots, n)$,
suppose
\[
\prod_{i=1}^n(\zeta_p^{a_i}-1)
\doteq \prod_{i=1}^n(\zeta_p^{b_i}-1),
\]
where $\zeta_p$ is a primitive $p$-th root of unity.
Then there exists a permutation $\sigma$ of $\{1, \ldots, n\}$
such that $a_i=\pm b_{\sigma(i)}$ for all $i=1, \ldots, n$.
In other words, 
$\{ \pm a_i \, (\mathrm{mod}\ \! p)\} = \{ \pm b_i \, \, (\mathrm{mod}\ \! p)\}$
as multiple sets.
\end{lemma}
%%%%%%%%%%%%%%%%%%%%%%%%%%%%%%%%%%%%
%
%
We will use this lemma in Lemma~\ref{lem:special} and in Section~\ref{sec:rem}.

%%%%%%%%%%%%%%%%%%%%%%%%%%%%%%%%%%%%%%%%%%
%%%                   SubSection 2-3
%%%%%%%%%%%%%%%%%%%%%%%%%%%%%%%%%%%%%%%%%%
\subsection{Combinatorial Euler structure 
(Reidemeister--Turaev torsion)}~\label{ssec:Euler}
Let $M$ be a homology lens space with 
$H=H_1(M)\cong \mathbb{Z}/p\mathbb{Z}\ (p\ge 2)$.
Then the Reidemeister torsion $\tau^{\psi_d}(M)$ of $M$ 
related with $\psi_d$ is determined up to multiplication of 
$\pm \zeta_d^m\ (m\in \mathbb{Z})$,
where $d\ge 2$ is a divisor of $p$ and 
$\psi_d$ is the same ring homomorphism as in Lemma \ref{lem:surgery2} (1).
Once we fix a basis of a cell chain complex for the maximal abelian covering
of $M$ as a $\mathbb{Z}[H]=\mathbb{Z}[t, t^{-1}]/(t^p-1)$-module, 
the value $\tau^{\psi_d}(M)$ is uniquely determined as
an element of $\mathbb{Q}(\zeta_d)$ for every $d$.
The choice of the basis up to ^^ ^^ base change equivalence" is called 
a {\it combinatorial Euler structure} of $M$ (cf. Turaev \cite{Tur2}).
The Reidemeister torsion of a manifold 
with a fixed combinatorial Euler structure
is said the {\it Reidemeister--Turaev torsion}.

We consider the sequence of 
the values $\tau^{\psi_d}(M)$ in $\mathbb{Q}(\zeta_d)$
of the Reidemeister--Turaev torsion
for every divisor $d\ge 2$ of $p$, and regard them
as a value sequence $\{\tau^{\psi_d}(M)\}_{d|p, d\ge 2}$
defined as below.
%
%
%%%%%%%%%%%%%%%%%%%%%%%%%%%%%%%%%%%%
\begin{definition}~\label{defini:seq}
A sequence of values $\boldsymbol{x}=\{x_d\}_{d|p, d\ge 2}$ 
is a {\it value sequence (of degree $p$)} 
if $x_d\in \mathbb{Q}(\zeta_d)$ for every $d$.
Two value sequences
$\boldsymbol{x}=\{x_d\}_{d|p, d\ge 2}$ and $\boldsymbol{y}=\{y_d\}_{d|p, d\ge 2}$ 
are {\it equal} ($\boldsymbol{x}=\boldsymbol{y}$) if $x_d=y_d$ for every $d$.
We are mainly concerned with the value sequence of type 
$\boldsymbol{x}=\{F(\zeta_d)\}_{d|p, d\ge 2}$
for a rational function $F(t)\in \mathbb{Q}(t)$.
In such a case, we say that $\boldsymbol{x}$ is {\it induced by $F(t)$}
and that $F(t)$ is a {\it lift} of $\boldsymbol{x}$.
A {\it control} of 
$\boldsymbol{x}=\{x_d\}_{d|p, d\ge 2}$ by a trivial unit
$u=\eta t^m\in \mathbb{Q}[t, t^{-1}]/(t^p-1)$
is defined by
\[
u\boldsymbol{x}=\{\eta \zeta_d^m x_d\}_{d|p, d\ge 2},
\]
where $\eta=1$ or $-1$ (constant) and $m\in \mathbb{Z}$.
Two value sequences
$\boldsymbol{x}=\{x_d\}_{d|p, d\ge 2}$ and $\boldsymbol{y}=\{y_d\}_{d|p, d\ge 2}$ 
are {\it control equivalent} if 
there is a trivial unit $u\in \mathbb{Q}[t, t^{-1}]/(t^p-1)$
such that $\boldsymbol{y}=u\boldsymbol{x}$.
A value sequence $\boldsymbol{x}=\{x_d\}_{d|p, d\ge 2}$ is 
a {\it real value sequence} if
$x_d$ is a real number for every $d$.
\end{definition}
%%%%%%%%%%%%%%%%%%%%%%%%%%%%%%%%%%%%
%
%
Let $M$ be a homology lens space with 
$H_1(M)\cong \mathbb{Z}/p\mathbb{Z}\ (p\ge 2)$.
Then a sequence $\{\tau^{\psi_d}(M)\}_{d|p, d\ge 2}$
of the Reidemeister torsions of $M$
with a combinatorial Euler structure
is a value sequence of degree $p$.
We say the value sequence
a {\it torsion sequence} of $M$.
%
%
%%%%%%%%%%%%%%%%%%%%%%%%%%%%%%%%%%%%
\begin{lemma}~\label{lem:lensEx}
\begin{enumerate}
\item[(1)]
Let $M$ and $M'$ be homeomorphic homology lens spaces 
with $H_1(M)\cong H_1(M')\cong \mathbb{Z}/p\mathbb{Z}\ (p\ge 2)$.
Then torsion sequences 
$\{\tau^{\psi_d}(M)\}_{d|p, d\ge 2}$
and $\{\tau^{\psi'_d}(M')\}_{d|p, d\ge 2}$
related with the corresponding 
representations $\psi_d$ and $\psi'_d$
(i.e., $\psi_d = \psi'_d \circ h_{\ast}$, where $h_{\ast}$ is the induced
homomorphism of the homeomorphism)
are control equivalent.

\item[(2)]
Let $M$ be a homology lens space
with $H_1(M)\cong \mathbb{Z}/p\mathbb{Z}\ (p\ge 2)$.
Then we can control a torsion sequence of $M$
into a real value sequence.
\end{enumerate}
\end{lemma}
%%%%%%%%%%%%%%%%%%%%%%%%%%%%%%%%%%%%
%
%
\medskip \noindent
{\bf Proof} \ 
(1)\ It is easy to see.

\medskip \noindent
(2)\ 
Here we let $\zeta$ denote any $d$-th primitive root ($\zeta_d$) of unity.
Since $M$ is obtained by $p/q$-surgery along a knot $K$
in a homology $3$-sphere for some $q$ (cf.\ \cite{BL}).
By Lemma 2.5 (1),
we have
\[
\tau^{\psi_d}(M)\doteq
{\mit \Delta}_K(\zeta)(\zeta-1)^{-1}(\zeta^{\bar{q}}-1)^{-1},
\]
where $q\bar{q}\equiv 1\ (\mathrm{mod}\ \! p)$.
By the duality of the Alexander polynomial
(cf.\ \cite{Mil, Tur1, Tur2}), we may assume
\[
{\mit \Delta}_K(t)={\mit \Delta}_K(t^{-1}).
\]
This is also a control of the combinatorial Euler structure
of the exterior of $K$, which induces 
a control of a torsion sequence of $M$.
We take an odd integer lift of $\bar{q}$.
Then 
\[
\zeta^{\frac{1+\bar{q}}{2}}
{\mit \Delta}_K(\zeta)(\zeta-1)^{-1}(\zeta^{\bar{q}}-1)^{-1}
\]
is a real number for every $d$.
\qed 
\par
%%
%%
%%%%%%%%%%%%%%%%%%%%%%%%%%%%%%%%%%%%%
%\begin{lem}~\label{lem:multiplication}
%{\sl
%Let
%$\mathbf{x}$, $\mathbf{y}$ and $\mathbf{z}$ 
%be three value sequences of degree $p$.
%Suppose that $\mathbf{x}$ and $\mathbf{y}$
%are (control) equivalent.
%Then $\mathbf{x}\cdot \mathbf{z}$ and 
%$\mathbf{y}\cdot \mathbf{z}$
%are (control) equivalent.
%}\end{lem}
%%%%%%%%%%%%%%%%%%%%%%%%%%%%%%%%%%%%%
%%
%%

%
%
%%%%%%%%%%%%%%%%%%%%%%%%%%%%%%%%%%%%
\begin{lemma}~\label{lem:control}
If two \underline{real} value sequences $\boldsymbol{x}=\{x_d\}_{d|p, d\ge 2}$ and 
$\boldsymbol{y}=\{y_d\}_{d|p, d\ge 2}$ of degree $p$
are control equivalent satisfying
$\boldsymbol{y}=u\boldsymbol{x}$
for a trivial unit $u = \eta t^m \in \mathbb{Z}[t, t^{-1}]/(t^p-1)$,
where $\eta= \pm 1$ and $m\in \mathbb{Z}$,
then the possibility of $u$ is restricted as follows:
\begin{enumerate}
\item[(i)]
If $p$ is odd, then $u=1$ or $-1$.

\item[(ii)]
If $p$ is even, then 
$u=1$, $-1$, $t^{p/2}$ or $-t^{p/2}$.

\end{enumerate}
\end{lemma}
%%%%%%%%%%%%%%%%%%%%%%%%%%%%%%%%%%%%
%
%
\medskip \noindent
{\bf Proof} \ 
Since the ratio $\zeta_p^m = \pm y_p/x_p$ is a real number,
we have (i) $m\equiv 0\ (\mathrm{mod}\ \! p)$ if $p$ is odd,
and (ii) $m\equiv 0\ \mbox{or}\ p/2\ (\mathrm{mod}\ \! p)$
if $p$ is even.
\qed 
\par

\bigskip
%
%
%%%%%%%%%%%%%%%%%%%%%%%%%%%%%%%%%%%%
\begin{definition}~\label{defini:symLPoly}
{\rm (Symmetric Laurent polynomial)} \
A Laurent polynomial $F(t)\in \mathbb{Q}[t, t^{-1}]$ is
{\it symmetric} if it is of the form
%%%%%
\[
F(t)=a_0+\sum_{i=1}^{\infty} a_i(t^i+t^{-i}),
\]
%%%%%
where $a_i$ is a rational number for all $i=1, 2, \ldots$
and $a_i=0$ for every sufficiently large $i$.
Note that, if $F(t)$ is a symmetric Laurent polynomial,
the induced value sequence $\{F(\zeta_d)\}_{d|p, d\ge 2}$
is a real value sequence.
We are concerned with symmetric Laurent polynomials
that are lifts (in $\mathbb{Q}[t, t^{-1}]$) of a polynomial 
in the quotient ring $\mathbb{Q}[t, t^{-1}]/(t^p-1)$.
We say that $F(t)$ (as above) is {\it reduced} if 
$a_i =0$ for all $i > [p/2]$. 
We often {\it reduce} the symmetric polynomials 
by using $t^i + t^{-i} = t^{p+i} + t^{-(p+i)}$ modulo $(t^p -1)$.
We let $\textrm{red}(F(t))$ denote the reduction of $F(t)$
(i.e., $\textrm{red}(F(t))$ is reduced and
$\textrm{red}(F(t)) = F(t)$ in $\Bbb{Q}[t, t^{-1}]/(t^p-1)$).

For a Laurent polynomial $F(t)\in \mathbb{Q}[t, t^{-1}]$,
the {\it span} of $F(t)$ is the difference of
the maximal degree and the minimal degree of $F(t)$,
and we denote it by $\mathrm{span}\ \! (F(t))$.
Note that 
the span of a symmetric Laurent polynomial is always even, 
and that
the span of a reduced symmetric Laurent polynomial
is less than or equal to $2[p/2]$. 
\end{definition}
%%%%%%%%%%%%%%%%%%%%%%%%%%%%%%%%%%%%
%
%

%
%
%%%%%%%%%%%%%%%%%%%%%%%%%%%%%%%%%%%%
\begin{lemma}~\label{lem:Euler2}
Let $N \ge 2$ be an integer.
Let $F(t), G(t)$ be symmetric Laurent polynomials and 
$\boldsymbol{x} = \{F(\zeta_d)\}_{d|N, d\ge 2}, \
\boldsymbol{y} = \{G(\zeta_d)\}_{d|N, d\ge 2}$
the induced real value sequences, respectively.
%We assume that  $\mathrm{span}\ \! (G(t))\le 2([N/2]-1)$.
If $\boldsymbol{x}$ and $\boldsymbol{y}$ are control equivalent, i.e., $u \boldsymbol{x} = \boldsymbol{y}$
for a trivial unit $u$ 
(here, $u=1$ or $-1$ if $N$ is odd, 
 $u=1, -1, t^{N/2}$ or $- t^{N/2}$ if $N$ is even, by Lemma~\ref{lem:control}),
then we have a congruence
\[
u F(t)\equiv G(t)\ \mod \ t^{N-1}+t^{N-2}+\cdots +t+1.
\]
Furthermore, assuming $\mathrm{span}\ \! (G(t))\le 2([N/2]-1)$,
\begin{enumerate}
\item[(i)] In the case that $u = 1$ or $-1$ and $\mathrm{span}(F(t)) \le N-2$,
we have an identity $u F(t) = G(t)$ in $\mathbb{Q}[t, t^{-1}]$.
\item[(ii)]
Otherwise (in the case that $N$ is even and $u = \eta t^{N/2}$ with $\eta = 1$ or $-1$),
if $\mathrm{span}\ \! ( \textrm{red}(t^{N/2}F(t))) \le N-2$,
then we have $\textrm{red}(t^{N/2} F(t))= \eta G(t)$ in $\mathbb{Q}[t, t^{-1}]$.
\end{enumerate}
\end{lemma}
%%%%%%%%%%%%%%%%%%%%%%%%%%%%%%%%%%%%
%
%
\medskip \noindent
{\bf Proof} \ 
By the Chinese Remainder Theorem, we have a ring isomorphism:
%%%%%
\begin{equation*}~\label{eq:Chinese1}
\mathbb{Q}[t, t^{-1}]/\left( t^{N-1}+t^{N-2}+\cdots +t+1\right)
\cong
\bigoplus_{d|N, d\ge 2}\mathbb{Q}(\zeta_d),
\end{equation*}
%%%%%
where $f(t)$ in the left-hand side maps to the value sequences 
$\{ f(\zeta_d)\}_{d|N, d\ge 2}$ in the right-hand side.
The isomorphism implies the required congruence.

Since $2([N/2]-1) <  N-1 = \textrm{span}(t^{N-1}+t^{N-2}+\cdots +t+1)$,
we have the identities.
\qed 
\par 

%
%
%%%%%%%%%%%%%%%%%%%%%%%%%%%%%%%%%%%%
\begin{lemma}~\label{lem:Euler}
Let $N \ge 2$ be an integer.
Let $F(t), G(t)$ be symmetric Laurent polynomials and 
$\boldsymbol{x} = \{F(\zeta_d)\}_{d|N, d\ge 2}, \
\boldsymbol{y} = \{G(\zeta_d)\}_{d|N, d\ge 2}$
the induced real value sequences, respectively.
If $\boldsymbol{x}$ and $\boldsymbol{y}$ are control equivalent, i.e., $u \boldsymbol{x} = \boldsymbol{y}$
for a trivial unit $u$ 
(here, $u=1$ or $-1$ if $N$ is odd, 
 $u=1, -1, t^{N/2}$ or $- t^{N/2}$ if $N$ is even, by Lemma~\ref{lem:control}),
and $F(1)=G(1)=0$,
then we have a congruence
\[
u F(t)\equiv G(t)\ \mod \ t^N -1.
\]
Furthermore, assuming $\mathrm{span}\ \! (G(t))\le 2[N/2]$,
\begin{enumerate}
\item[(i)] In the case that $u = 1$ or $-1$ and $\mathrm{span}(F(t)) \le N-1$,
we have an identity $u F(t) = G(t)$ in $\mathbb{Q}[t, t^{-1}]$.
\item[(ii)]
Otherwise (in the case that $N$ is even and $u = \eta t^{N/2}$ with $\eta = 1$ or $-1$),
we have $\textrm{red}(t^{N/2} F(t))= \eta G(t)$ in $\mathbb{Q}[t, t^{-1}]$.
\end{enumerate}
\end{lemma}
%%%%%%%%%%%%%%%%%%%%%%%%%%%%%%%%%%%%
%
%
\medskip \noindent
{\bf Proof} \ 
We use the same argument with the proof of Lemma~\ref{lem:Euler2}, 
but here we use the Chinese Remainder Theorem for the following ring isomorphism
%%%%%
\begin{equation*}~\label{eq:Chinese2}
\mathbb{Q}[t, t^{-1}]/\left( t^N-1\right)
\cong
\bigoplus_{d|N, d\ge 1}\mathbb{Q}(\zeta_d).
\end{equation*}
%%%%%
\qed 
\par 

\medskip

Note that $F(t)$ and $t^{N/2} F(t)$ 
induce the control equivalent real value sequences by $u = t^{N/2}$, 
but $\textrm{red}(t^{N/2}F(t)) \not= F(t)$ in general. 
Thus we have to care the case (ii) in Lemma~\ref{lem:Euler2} and \ref{lem:Euler}.
Here, we study relation between the coefficients of $F(t)$
and those of $\textrm{red}(t^{N/2}F(t))$.
%
%
%%%%%%%%%%%%%%%%%%%%%%%%%%%%%%%%%%%%
\begin{lemma}~\label{lem:symeven}
Let $N$ be an even integer. 
\[
\textrm{If } \ F(t)= a_0+\sum_{i=1}^{N/2} a_i(t^i+t^{-i}), 
\ \textrm{then } \ 
\textrm{red}(t^{N/2}F(t))= b_0+\sum_{i=1}^{N/2} b_i(t^i+t^{-i})
\]
with
\begin{equation*}
b_0 = 2 a_{N/2}, \ b_{N/2} = a_0/2 \textrm{ and }
b_j = a_{N/2 - j }
\quad  (j =1,2, \ldots, N/2-1).
\end{equation*}
\end{lemma}
%%%%%%%%%%%%%%%%%%%%%%%%%%%%%%%%%%%%
%
%
\medskip \noindent
{\bf Proof} \ 
It is because
\[
t^{N/2}(t^j + t^{-j}) \ = \ 
t^{N/2 +j } +  t^{N/2 -j} \ \equiv \ 
t^{(N/2 -j)} +  t^{-(N/2 -j)} \ \mod t^N -1.
\]
\qed 
\par 

\medskip

This will be used in the proof of Lemma~\ref{lem:special},
see also Remark~\ref{rem:Nes(b)}.

%%%%%%%%%%%%%%%%%%%%%%%%%%%%%%%%%%%%%%%%%%
%%%                   Section 3
%%%%%%%%%%%%%%%%%%%%%%%%%%%%%%%%%%%%%%%%%%
\section{Alexander polynomial of $A_{m,n}$}~\label{sec:Alex}
We compute the Alexander polynomial of the link $A_{m,n}$. 
%
%
%%%%%%%%%%%%%%%%%%%%%%%%%%%%%%%%%%%%
\begin{definition}~\label{defini:Imn}
For a coprime positive pair $(m, n)$, 
we define a set ${\frak I}(m, n)$ by 
%%%%%
\begin{eqnarray*}
{\frak I}(m, n) 
& := & 
(m\mathbb{Z} \cup n\mathbb{Z}) \cap 
\{ k \in \mathbb{Z}\ \vert \ 0 \le k \le mn \} \\
& = &
\{0, m, 2m, \ldots, nm \} \cup 
\{0, n, 2n, \ldots, mn \}.
\end{eqnarray*}
%%%%%
Note that the cardinality of ${\frak I}(m, n)$ is $m+n$.
We sort the all elements in ${\frak I}(m, n)$ as 
\[
0 = k_0 < k_1 < k_2 < \cdots < k_{m+n-1} = mn \qquad
(k_i \in {\frak I}(m, n)).
\]
Here, $k_1$ is the smaller one in $m$ and $n$.
\end{definition}
%%%%%%%%%%%%%%%%%%%%%%%%%%%%%%%%%%%%
%
%
The goal of this section is: 
%
%
%%%%%%%%%%%%%%%%%%%%%%%%%%%%%%%%%%%%
\begin{theorem}~\label{thm:AmnAlex}
The Alexander polynomial of $A_{m,n}$ is 
\[
{\it \Delta}_{A_{m,n}} (t,x) 
\ \doteq \
\sum_{i=0}^{m+n-1} t^{k_i} x^i,
\]
where $t$ (and $x$, respectively) is represented by
a meridian of $K_1$ (that of $K_2$).
\end{theorem}
%%%%%%%%%%%%%%%%%%%%%%%%%%%%%%%%%%%%
%
%
%%%%%%%%%%%%%%%%%%%%%%%%%%%%%%%%%%%%
\begin{example}~\label{ex:m3n5_1}
In the case $(m, n)=(3,5)$, we have 
${\frak I}(3,5) = \{0, 3, 5, 6, 9, 10, 12, 15 \}$ and 
\[
{\it \Delta}_{A_{3,5}} (t,x) \doteq
t^{15}x^7
+ t^{12}x^6
+ t^{10}x^5
+ t^9x^4
+ t^6x^3
+ t^5x^2
+ t^3x
+1.
\]
\end{example}
%%%%%%%%%%%%%%%%%%%%%%%%%%%%%%%%%%%%
%
%

%%%%%%%%%%%%%%%%%%%%%%%%%%%%%%%%%%%%%%%%%%
%%%                   SubSection 3-1
%%%%%%%%%%%%%%%%%%%%%%%%%%%%%%%%%%%%%%%%%%
\subsection{Alexander matrix of $A_{m,n}$}~\label{ssec:matrix}
We start the proof of Theorem~\ref{thm:AmnAlex}
%(the formula of the Alexander polynomial of the link $A_{m,n}$) 
with the following lemma: 
%
%
%
%%%%%%%%%%%%%%%%%%%%%%%%%%%%%%%%%%%%
\begin{lemma}~\label{lem:AmnAlexMat}
The Alexander matrix of $A_{m,n}$ is 
\[
I_{m+n-1} - x M(m, n)
\] 
and the Alexander polynomial of $A_{m,n}$ is obtained by
\[
{\it \Delta}_{A_{m,n}} (t,x) 
\ \doteq \
\det \left (
I_{m+n-1} - x M(m, n)
\right ),
\] 
where $M(m, n)$ is the $(m+n-1)\times (m+n-1)$-matrix of the form
\[
M(m, n) := 
\begin{bmatrix}
O_{n-1,m-1} & -\overrightarrow{T_{n-1}} & I_{n-1} \\[5pt]
o_{m-1} & -t^n & o_{n-1}  \\[5pt]
t^n I_{m-1} &  -t^n\overrightarrow{i_{m-1}} & O_{m-1,n-1} \\[5pt]
\end{bmatrix},
\]
$\overrightarrow{T_{n-1}}$, $\overrightarrow{i_{m-1}}$ are
the following column vectors of size $(n-1) \times 1$
and $(m-1) \times 1$ respectively
\[
\overrightarrow{T_{n-1}}
:= 
{\small
\begin{bmatrix}
t \\ t^2 \\ \vdots \\ t^{n-1}
\end{bmatrix}
},
\qquad
\overrightarrow{i_{m-1}}
:= 
{\small
\begin{bmatrix}
1 \\1 \\ \vdots \\ 1
\end{bmatrix}
},
\]
$O_{s,s'}$ (and $o_s$, respectively)
is the zero matrix of size $s \times s'$ (and of size $1 \times s$)
and $I_s$ is the identity matrix of size $s \times s$.
\end{lemma}
%%%%%%%%%%%%%%%%%%%%%%%%%%%%%%%%%%%%
%
%
%%%%%%%%%%%%%%%%%%%%%%%%%%%%%%%%%%%%
\begin{example}~\label{ex:m3n5_Mat}\\
$M(3, 5) =
{\footnotesize
\begin{bmatrix}
0 & 0 & -t     & 1 & 0 & 0 & 0 \\
0 & 0 & -t^2 & 0 & 1 & 0 & 0 \\
0 & 0 & -t^3 & 0 & 0 & 1 & 0 \\
0 & 0 & -t^4 & 0 & 0 & 0 & 1 \\
0 & 0 & -t^5 & 0 & 0 & 0 & 0 \\
t^5 & 0 & -t^5 & 0 & 0 & 0 & 0 \\
0 & t^5 & -t^5 & 0 & 0 & 0 & 0 \\
\end{bmatrix}
}$, \ 
$ I - xM(3, 5) =
{\footnotesize
\begin{bmatrix}
1 & 0 & tx     & -x & 0 & 0 & 0 \\
0 & 1 & t^2x & 0 & -x & 0 & 0 \\
0 & 0 & 1+t^3x & 0 & 0 & -x & 0 \\
0 & 0 & t^4x & 1 & 0 & 0 & -x \\
0 & 0 & t^5x & 0 & 1 & 0 & 0 \\
-t^5x & 0 & t^5x & 0 & 0 & 1 & 0 \\
0 & -t^5x & t^5x & 0 & 0 & 0 & 1 \\
\end{bmatrix}
}$.
\end{example}
%%%%%%%%%%%%%%%%%%%%%%%%%%%%%%%%%%%%
%
%
\medskip \noindent
{\bf Proof of Lemma~\ref{lem:AmnAlexMat}} \ 
We assumed that $2 \le m<n$ (Section~\ref{sec:intro}). 
%
%
%%%%%%%%%%%%%%%%%
\begin{figure}[h]
\begin{center}
\includegraphics[scale=0.4]{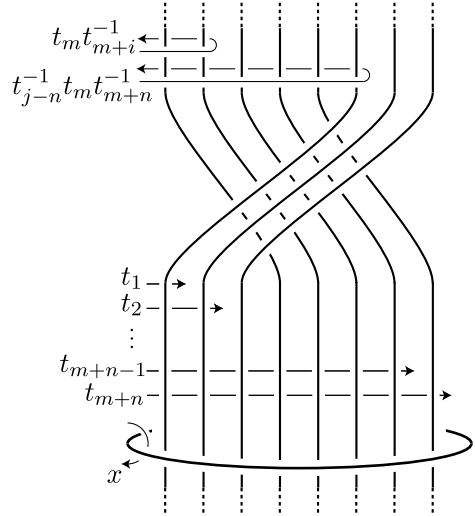}
\caption{Generators $t_i, t_j $ \quad (ex. $(m, n)=(3, 5)$)}
\label{fig:AmnBraid}
\end{center}
\end{figure}
%%%%%%%%%%%%%%%%%
%
%
We take the generators of the fundamental group 
$\pi_1 (S^3 \backslash A_{m,n})$ as in Figure~\ref{fig:AmnBraid}.
Then we have a presentation of the group:
%%%%%
\begin{eqnarray*}
& \left \langle
\begin{matrix}
t_1, t_2, \ldots, t_{m+n} \\
x
\end{matrix}  
\ \left \vert \
\begin{matrix}
L_i := t_i x t_{m} t_{m+i}^{-1} x^{-1} \hfill (i = 1, 2, \ldots , n)  \\[5pt]
R_j := t_j x t_{j-n}^{-1} t_{m} t_{m+n}^{-1} x^{-1} \quad  (j= n+1, n+2, \ldots , m+n) \\
\end{matrix} 
\right.
\right \rangle \\
%%%
& = 
\left \langle
\begin{matrix}
t_1, t_2, \ldots, t_{m+n-1} \\
x
\end{matrix}  
\ \left \vert \
\begin{matrix}
L_i := t_i x t_{m} t_{m+i}^{-1} x^{-1} \hfill (i = 1, 2, \ldots , n-1)  \\[5pt]
C := t_n x t_m x t_m^{-1}x^{-1} t_n^{-1}x^{-1} \hfill \\
R_j := t_j x t_{j-n}^{-1} x^{-1} t_{n}^{-1}\quad  (j = n+1, \ldots , m+n-1) \\
\end{matrix} 
\right.
\right \rangle,
\end{eqnarray*}
%%%%%
where we canceled $t_{m+n} = x^{-1}t_n x t_m$ by the relation $L_n$. 
Then $R_{m+n}$ is changed to the relation $C$, 
and the others $R_j$ are also changed.
The non-zero free differentials of the relations $L_i$ and $R_j$ 
by $t_{\ast}$ are:
%%%%%
\begin{eqnarray}~\label{eq:jrow}
\begin{matrix}
\dfrac{\partial L_i}{\partial t_i} = 1, 
& 
\dfrac{\partial L_i}{\partial t_m} = \delta_{im} + t_i x, \hfill
& 
\dfrac{\partial L_i}{\partial t_{m+i}} = -t_i x t_m t_{m+i}^{-1}, \hfill
& \\
\\
\dfrac{\partial R_j}{\partial t_j} = 1,
& 
\dfrac{\partial R_j}{\partial t_{j-n}} = -t_j x t_{j-n}^{-1}, \hfill
& 
\dfrac{\partial R_j}{\partial t_{n}} = -t_j x t_{j-n}^{-1} x^{-1} t_{n}^{-1}, \\
\end{matrix}
\end{eqnarray}
%%%%%
where $\delta_{im}$ is Kronecker's delta.
The non-zero free differentials of $C$ by $t_{\ast}$ are
%%%%%
\begin{eqnarray}~\label{eq:Crow}
\dfrac{\partial C}{\partial t_{m}} = 
t_n x - t_n x t_m x t_m^{-1}, \qquad
\dfrac{\partial C}{\partial t_{n}}  =
1 - t_n x t_m x t_m^{-1}x^{-1} t_n^{-1}.
\end{eqnarray}
%%%%%
We let $\gamma$ denote the Hurewicz epimorphism 
$\gamma :  \pi_1 (S^3 \backslash A_{m,n})
\rightarrow H_1 (S^3 \backslash A_{m,n}; \mathbb{Z}) \cong \langle t, x \vert \, -\rangle$,
defined as $\gamma (t_i) = t^i, \gamma (x) = x$.
We redefine $L_n := C$ and $L_i := R_i$ for $n+1 \le i \le m+n-1$.
Then the well-known formula (\cite[Lemma~7.3.2 (p.93)]{Kaw}) on 
the Alexander polynomial of links says
\[
\det 
\begin{bmatrix}
\gamma \left ( 
\dfrac{\partial L_i}{\partial t_j}
\right ) 
\end{bmatrix}
\dot= \ (1-x)
{\it \Delta}_{A_{m,n}}(t,x),
\]
where $j$ (columns) runs in $1 \le j \le m+n-1$.

\medskip

A submatrix from the first row to the $(n-1)$-th row of 
$\begin{bmatrix}
\gamma ( \partial L_i/ \partial t_j ) 
\end{bmatrix}$ 
coincides with  $I - x M(m, n)$
by the first half of (\ref{eq:jrow}).
In the $n$-th row of 
$\begin{bmatrix}
\gamma ( \partial L_i/ \partial t_j ) 
\end{bmatrix}$, 
non-zero entries are only
%%%%%
\begin{eqnarray*}
\gamma \left ( 
\dfrac{\partial C}{\partial t_{m}} 
\right ) 
= (1-x) t^n x 
\ \textrm{ at $(n, m)$ \ and } \
\gamma \left ( 
\dfrac{\partial C}{\partial t_{n}} 
\right ) 
= (1-x)
\ \textrm{ at $(n, n)$},
\end{eqnarray*}
%%%%%
by (\ref{eq:Crow}).
Thus the $n$-th row of 
$\begin{bmatrix}
\gamma ( \partial L_i/ \partial t_j ) 
\end{bmatrix}$ 
coincides with  
$(1-x)$ times the $n$-th row of 
$I - x M(m, n)$.
We add $1/(1-x)$ times the $n$-th row of 
$\begin{bmatrix}
\gamma ( \partial L_i/ \partial t_j ) 
\end{bmatrix}$ 
to each $j$-th row of 
$\begin{bmatrix}
\gamma ( \partial L_i/ \partial t_j ) 
\end{bmatrix}$ 
with $j \geq n+1$.
The resulting $j$-th row coincides with  
the $j$-th row of $I - x M(m, n)$ by
the second half of (\ref{eq:jrow}).
Therefore we have
\[
\det 
\begin{bmatrix}
\gamma \left ( 
\dfrac{\partial L_i}{\partial t_j}
\right ) 
\end{bmatrix}
\ \doteq \
(1-x) \cdot \det ( I - x M(m, n)).
\]
\qed 
\par
%
%
%%%%%%%%%%%%%%%%%%%%%%%%%%%%%%%%%%%%
\begin{remark}~\label{rem:Burau}
{\rm 
The matrix $M(m, n)$ can also be obtained by 
the Burau representation (see Birman~\cite[p.121]{Bir}) of 
the braid group.
}
\end{remark}
%%%%%%%%%%%%%%%%%%%%%%%%%%%%%%%%%%%%
%
%

%%%%%%%%%%%%%%%%%%%%%%%%%%%%%%%%%%%%%%%%%%
%%%                   SubSection 3-2
%%%%%%%%%%%%%%%%%%%%%%%%%%%%%%%%%%%%%%%%%%
\subsection{Properties of ${\frak I}(m, n)$}~\label{ssec:Imn}
We need some properties on elements $k_i$ of ${\frak I}(m, n)$, see 
Definition~\ref{defini:Imn}.
They will be also used in Subsection~\ref{ssec:R-tor}, Subsection~\ref{ssec:type}.

\medskip 
%%%%%%%%%%%%%%%%%%%%%%%%%%%%%%%%%%%%
\begin{definition}~\label{defini:uj}
We define $u_j\ (j=0, 1, 2, \ldots, n)$ and
$w_j\ (j=0, 1, 2, \ldots, m)$ by
\[
k_{u_j}=jm\quad \mbox{and}\quad k_{w_j}=jn.
\]
\end{definition}
%%%%%%%%%%%%%%%%%%%%%%%%%%%%%%%%%%%%
It is easy to see:
%
%
%%%%%%%%%%%%%%%%%%%%%%%%%%%%%%%%%%%%
\begin{proposition}~\label{prop:Imn}
\begin{enumerate}
\item[(1)]
$k_i+k_{m+n-1-i}=mn$,
\ ($i=0, 1, \ldots , m+n-1$).

\item[(2)]
$u_0 = w_0 =0$, $u_n = w_m =m+n-1$, and 
\[
u_j =
\left[\frac{jm}{n}\right]+ j
\quad
(j = 1,2, \ldots, n-1),
\qquad 
w_j =
\left[\frac{jn}{m}\right]+ j
\quad
(j = 1,2, \ldots, m-1),
\]
where $[\ \cdot\ ]$ is the gaussian symbol.

\item[(3)]
Both $u_j, w_j$ are increasing sequences.

\item[(4)]
$u_j+u_{n-j} = w_j+w_{m-j}=m+n-1$.

\item[(5)]
\[
\sum_{i=0}^{m+n-1}t^{k_i}x^i
=1+\sum_{i=1}^{n-1}t^{im}x^{u_i}+\sum_{j=1}^{m-1}t^{jn}x^{w_j}
+t^{mn}x^{m+n-1}.
\]
\end{enumerate}
\end{proposition}
%%%%%%%%%%%%%%%%%%%%%%%%%%%%%%%%%%%%
%
%
%%%%%%%%%%%%%%%%%%%%%%%%%%%%%%%%%%%%
\begin{notation}~\label{nota:N} 
For an integer $N$, we denote by $[N]_n$ the unique integer satisfying
\[
[N]_n\equiv N\quad (\mathrm{mod}\ \! n)
\quad \mbox{and}\quad
0\le [N]_n\le n-1.
\]
We let $\sigma$ (and its inverse $\sigma'$, respectively) denote 
the bijection
\[
\sigma,\ \sigma': \{0, 1, 2, \ldots, n-1\}\to \{0, 1, 2, \ldots, n-1\}
\]
defined by $\sigma(i) :=[i\overline{m}]_n$ (and $\sigma'(i) :=[im]_n$),
where $\overline{m}$ is regarded as an integral lift of
$\overline{m}\ (\mathrm{mod}\ \! n)$.
We also define
\[
\rho(i):=\frac{\sigma(i)m-i}{n}\quad (i=1, 2, \ldots, n-1).
\]
By the definition of $\sigma(i)$, $\rho(i)$ is an integer.
\end{notation}
%
%
%%%%%%%%%%%%%%%%%%%%%%%%%%%%%%%%%%%%
\begin{lemma}~\label{lem:Imn2}
\begin{enumerate}
\item[(1)]
$\sigma(0)=0$, $\sigma(m)=1$ and $\sigma(n-m)=n-1$.

\item[(2)]
For $i=1, 2, \ldots, n-1$ and $i\ne n-m$,
we have $u_{\sigma(i)+1}-u_{\sigma(i)}=1$ or $2$.
Furthermore, $u_{\sigma(i)+1}-u_{\sigma(i)}=2$ if and only if
$n-m+1\le i\le n-1$.

\item[(3)]
$\sigma(m+i)=\sigma(i)+1$,
$u_{\sigma(m+i)}=u_{\sigma(i)}+1$, \
$(i=1, 2, \ldots, n-m-1)$.

\item[(4)]
$\sigma(n-m+i)=\sigma(i)-1$,
$u_{\sigma(n-m+i)}=u_{\sigma(i)}-2$, \
$(i=1, 2, \ldots, m-1)$.

\item[(5)]
$\rho : \{1, 2, \ldots, m-1\}\to \{1, 2, \ldots, m-1\}$
is a bijection.

\item[(6)]
$w_{\rho(i)}=u_{\sigma(i)}-1, \
(i=1, 2, \ldots, m-1)$.

\end{enumerate}
\end{lemma}
%%%%%%%%%%%%%%%%%%%%%%%%%%%%%%%%%%%%
%
%
\medskip \noindent
{\bf Proof} \ 
(1) It is easy to see.

\medskip
\noindent
(2)\ Suppose that $i=1, 2, \ldots, n-1$ and $i\ne n-m$.
Since $\sigma$ is a bijection, (1) and Proposition~\ref{prop:Imn} (3), 
we have $1\le \sigma(i)\le n-2$ and
$u_{\sigma(i)+1}-u_{\sigma(i)}=1$ or 2.
We set $j=\sigma(i)$.
It holds that $u_{j+1}-u_j=2$ if and only if 
there exists an integer $j'$ such that $jm<j'n<(j+1)m$,
which implies $j'n-m<jm<j'n$ and
\[
n-m<[jm]_n<n.
\]
Since $[jm]_n=\sigma'(j)=\sigma'\circ \sigma(i)=i$,
we have $n-m+1\le i\le n-1$.

\medskip
\noindent
(3)\ Suppose that $i=1, 2, \ldots, n-m-1$.
We have $1\le \sigma(i)\le n-2$ again.
Since 
$\sigma(m+i)\equiv (m+i)\overline{m}\equiv \sigma(i)+1\ 
(\mathrm{mod}\ \! n)$, we have $\sigma(m+i)=\sigma(i)+1$.
By (2), we have 
$u_{\sigma(m+i)} =u_{\sigma(i)+1} = u_{\sigma(i)}+1$.

\medskip
\noindent
(4)\ Suppose that $i=1, 2, \ldots, m-1$.
We have $2\le \sigma(i)\le n-1$.
Since 
$\sigma(n-m+i)\equiv (n-m+i)\overline{m}\equiv \sigma(i)-1\ 
(\mathrm{mod}\ \! n)$, we have $\sigma(n-m+i)=\sigma(i)-1$.
By (2) and $n-m+1\le n-m+i\le n-1$, 
we have $u_{\sigma(n-m+i)}=u_{\sigma(i)}-2$.

\medskip
\noindent
(5)\ Suppose that $i= 1, 2, \ldots, m-1$.
Since 
$\sigma(i)m=[i\overline{m}]_nm\equiv i\ (\mathrm{mod}\ \! n)$,
and $0\le \sigma(i)\le n-1$, we have
\[
0\le \rho(i)=\frac{\sigma(i)m-i}{n}\le \frac{(n-1)m}{n}<m,
\]
and hence the image of the map $\rho$ is included in $\{0, 1, 2, \ldots, m-1\}$.
By the definition, we have $\rho(i)n\equiv -i\ (\mathrm{mod}\ \! m)$,
thus $\rho$ is bijective. In fact, $\rho(i) = [(m-i)\overline{n}]_m$.

\medskip
\noindent
(6)\ Suppose that $i= 1, 2, \ldots, m-1$.
By (4), we have $u_{\sigma(n-m+i)}=u_{\sigma(i)-1}=u_{\sigma(i)}-2$.
Then, as in the proof of (2), there exists an integer $j$ such that $w_j=u_{\sigma(i)}-1$
such that 
\[
(\sigma(i)-1)m 
= k_{u_{\sigma(i)-1}} 
< k_{w_j} 
= jn
<k_{u_{\sigma(i)}}
=\sigma(i)m.
\]
Since $\sigma(i)m-i\equiv 0\ (\mathrm{mod}\ \! n)$,
we have $jn=\sigma(i)m-i$ and $j=\rho(i)$.
\qed 
\par
%
%
%%%%%%%%%%%%%%%%%%%%%%%%%%%%%%%%%%%%
\begin{definition}~\label{defini:etx}
We define the monomial $e_i$ of $t, x$ by
\begin{equation*}
e_i = t^{\sigma(i)m-i}x^{u_{\sigma(i)}},\quad (i=1, 2, \ldots, n-1).
\end{equation*}
\end{definition}
%%%%%%%%%%%%%%%%%%%%%%%%%%%%%%%%%%%%
%
%
%%%%%%%%%%%%%%%%%%%%%%%%%%%%%%%%%%%%
\begin{lemma}~\label{lem:ex}
\begin{enumerate}
\item[(1)]
$e_m=x$ and $e_{n-m}=t^{mn-n}x^{m+n-2}$.

\item[(2)]
$e_{m+i}=xe_i, \
(i=1, 2, \ldots, n-m-1)$.

\item[(3)]
$e_{n-m+i}=t^{-n}x^{-2}e_i, \ 
(i=1, 2, \ldots, m-1)$.

\item[(4)]
$t^{\sigma (i)m-i}x^{u_{\sigma(i)}-1}
= t^{\rho(i)n}x^{w_{\rho(i)}}, \
(i=1, 2, \ldots, m-1)$.
\end{enumerate}
\end{lemma}
%%%%%%%%%%%%%%%%%%%%%%%%%%%%%%%%%%%%
%
%
\medskip \noindent
{\bf Proof} \ 
They are proved by the Lemma~\ref{lem:Imn2} (1),(3),(4) and (6) respectively.
\qed

%%%%%%%%%%%%%%%%%%%%%%%%%%%%%%%%%%%%%%%%%%
%%%                   SubSection 3-3
%%%%%%%%%%%%%%%%%%%%%%%%%%%%%%%%%%%%%%%%%%
\subsection{Proof of Theorem~\ref{thm:AmnAlex}\ 
(Alexander polynomial of $A_{m,n}$)}~\label{ssec:AmnAlex}
We regard 
\begin{equation}~\label{eq:algeq}
\sum_{i=0}^{m+n-1}t^{k_i}x^i=0
\end{equation}
as an algebraic equation over $\mathbb{Z}[t, t^{-1}]$,
and $x$ as a root of the equation.
%
%
%%%%%%%%%%%%%%%%%%%%%%%%%%%%%%%%%%%%
\begin{lemma}~\label{lem:root}
The algebraic equation (\ref{eq:algeq})
over $\mathbb{Z}[t, t^{-1}]$ has no multiple root.
\end{lemma}
%%%%%%%%%%%%%%%%%%%%%%%%%%%%%%%%%%%%
%
%
\medskip \noindent
{\bf Proof} \ 
If we substitute $t=1$ into the equation (\ref{eq:algeq}), then we have
\[
\sum_{i=0}^{m+n-1}x^i=0.
\]
Since this equation has no multiple root, 
and its degree is equal to that of (\ref{eq:algeq}) ($=m+n-1$),
the equation (\ref{eq:algeq}) has no multiple root.
\qed 
\par

\bigskip

Let $\boldsymbol{v}$ be a row vector of size $1\times (m+n-1)$,
decomposed as 
\[
\boldsymbol{v} =
\begin{bmatrix}
\boldsymbol{e} & f & \boldsymbol{g} \\
\end{bmatrix},
\]
where $f$ is a scalar and 
\[
\boldsymbol{e} = 
\begin{bmatrix}
e_1 & e_2 & \cdots & e_{n-1}
\end{bmatrix},
\quad 
\boldsymbol{g}
= \begin{bmatrix}
g_1 & g_2 & \cdots & g_{m-1}
\end{bmatrix}
\]
are row vectors of size $1\times (n-1)$ and $1\times (m-1)$,
respectively.
Then the matrix $M(m,n)$ in Lemma~\ref{lem:AmnAlexMat}
satisfies
%%%%%
\begin{equation}~\label{eq:prodMat}
\boldsymbol{v} \cdot M(m, n)=
\begin{bmatrix}
t^n \boldsymbol{g} 
& 
\left (
- \boldsymbol{e}\cdot \overrightarrow{T_{n-1}} - t^nf
-t^n\boldsymbol{g}\cdot \overrightarrow{i_{m-1}}
\right ) 
& 
\boldsymbol{e} \\
\end{bmatrix}.
\end{equation}
%%%%%
%
%
%%%%%%%%%%%%%%%%%%%%%%%%%%%%%%%%%%%%
\begin{lemma}~\label{lem:eigen}
If we take $e_i = t^{\sigma(i)m-i}x^{u_{\sigma(i)}} \ 
(i=1, 2, \ldots, n-1)$ as in Definition~\ref{defini:etx}, 
$f=xe_{n-m}$ and $g_i=t^{-n}x^{-1}e_i\ (i=1, 2, \ldots, m-1)$,
then we have
\[
\boldsymbol{v} \cdot M(m, n)=x^{-1}\boldsymbol{v}.
\]
\end{lemma}
%%%%%%%%%%%%%%%%%%%%%%%%%%%%%%%%%%%%
%
%
\medskip \noindent
{\bf Proof} \ 
By Lemma~\ref{lem:ex}, 
and the definition of $e_i$ and $g_i$, 
we can see immediately that
the $i$-th entry of $\boldsymbol{v} \cdot M(m, n)$
is equal to that of $M(m, n)$ multiplied by $x^{-1}$
except the case $i=m$.
The $m$-th entry of (\ref{eq:prodMat}) is computed as
\begin{eqnarray*}
-\boldsymbol{e}\cdot \overrightarrow{T_{n-1}}-t^nf
-t^n\boldsymbol{g}\cdot \overrightarrow{i_{m-1}} & = &
-\sum_{i=1}^{n-1}t^ie_i-t^nf-t^n\sum_{j=1}^{m-1}g_j
\medskip \\
& = &
-\sum_{i=1}^{n-1}t^{\sigma(i)m}x^{u_{\sigma(i)}}
-t^{mn}x^{m+n-1}
-\sum_{j=1}^{m-1}t^{\rho(j)n}x^{w_{\rho(j)}}
\medskip \\
& = & 1
\end{eqnarray*}
by Definition~\ref{defini:etx},
%by the definition~\ref{defini:etx} of $e_i$ and $g_i$,
Lemma~\ref{lem:ex},
bijectivities of $\sigma$ and $\rho$,
the equation (\ref{eq:algeq}) and 
Proposition~\ref{prop:Imn} (5). We have the lemma.
\qed 
\par

\medskip \noindent
{\bf Proof of Theorem~\ref{thm:AmnAlex}} \ 
By Lemma~\ref{lem:eigen},
$x^{-1}$ is an eigenvalue of the matrix $M(m, n)$
where $x$ is a root of (\ref{eq:algeq}).
Since the degree of the equation (\ref{eq:algeq})
and the size of $M(m, n)$ are identical ($=m+n-1$),
and the equation (\ref{eq:algeq}) have no multiple root
by Lemma~\ref{lem:root}, we have
\[
\sum_{i=0}^{m+n-1}t^{k_i}x^i
= x^{m+n-1} \cdot \det \left (
x^{-1} I_{m+n-1} - M(m, n)
\right ).
\]
By Lemma~\ref{lem:AmnAlexMat},
we have the result.
\qed

%%%%%%%%%%%%%%%%%%%%%%%%%%%%%%%%%%%%%%%%%%
%%%                   Section 4
%%%%%%%%%%%%%%%%%%%%%%%%%%%%%%%%%%%%%%%%%%
\section{Reidemeister torsions of $(A_{m,n}; r, 0)$ and key lemmas}~\label{sec:key}
We compute the Reidemeister torsions of $(A_{m,n}; r, 0)$.
The goal is Lemma~\ref{lem:Rmn}.

%%%%%%%%%%%%%%%%%%%%%%%%%%%%%%%%%%%%%%%%%%
%%%                   SubSection 4-1
%%%%%%%%%%%%%%%%%%%%%%%%%%%%%%%%%%%%%%%%%%
\subsection{The first homology of $(A_{m,n}; r,0)$}~\label{ssec:hom1}
%Conditions from the first homology of $(A_{m,n}; r,0)$
%
We calculate the first homology of $M = (A_{m,n}; \alpha/\beta,0)$.
Let $E$ denote the complement of $A_{m,n}$. We regard
$M = E\cup V_1\cup V_2$, see Notation in Subsection~\ref{ssec:surgery}
We set $M_1 := E\cup V_1 \subset M$.

\medskip

From now on, we always assume that 
$\gcd(\alpha, \beta)=1$, $\beta>0$, and 
%%%%%
\begin{eqnarray*}
\gcd(m+n, \alpha)=1
\end{eqnarray*}
%%%%%
which is equivalent to the condition for 
the first homology $H_1(M; \mathbb{Z})$ to be finite cyclic
by the elementary divisor theory.
Then the order is $(m+n)^2\beta$:
$H_1(M; \mathbb{Z})\cong \mathbb{Z}/(m+n)^2\beta \mathbb{Z}$.

\medskip

We determine the first homologies of $E$, $M_1$ and $M$,
define generators and study relations.
First, $H_1(E)$ is a free abelian group of rank 2
generated by $[m_1]$ and $[m_2]$:
%%%%%
\begin{equation*}
H_1(E)\cong \langle [m_1], [m_2]\ \vert \ -\rangle
\cong \mathbb{Z}^2.
\end{equation*}
We have
\begin{equation}~\label{eq:EA}
[l_1] = [m_2]^{m+n}\quad
\mbox{and}\quad
[l_2] = [m_1]^{m+n}.
\end{equation}
%%%%%
Next, we attach $V_1$ to $E$ to make $M_1$. 
We take integers $\gamma, \delta$ such that
$\alpha \delta-\beta \gamma = -1$, and fix the 
meridian-longitude system $m_1', l_1'$ of the solid torus $V_1$.
In $H_1(M_1)$, we have the relations (\ref{eq:EA}) and 
%%%%%
\begin{equation*}
[m'_1] = [m_1]^{\alpha}[l_1]^{\beta}=1
\quad \mbox{and} \quad
[l'_1]=[m_1]^{\gamma}[l_1]^{\delta}.
\end{equation*}
Thus, we have
\begin{eqnarray*}
H_1(M_1) & \cong & \langle [m_1], [m_2]\ \vert \ 
[m_1]^{\alpha}[m_2]^{(m+n)\beta}=1\rangle \nonumber\\
& \cong & \langle T\ \vert \ -\rangle \cong \mathbb{Z},
\end{eqnarray*}
%%%%%
where $T=[m_1]^{\gamma'}[m_2]^{\delta'}$\ 
by taking integers $\gamma', \delta'$ satisfying 
$\alpha \delta'-p\beta \gamma'=-1$.
By the relations above, we have
%%%%%
\begin{eqnarray}~\label{eq:M2A}
\begin{matrix}
[m_1] 
& = & [m_1]^{-\alpha \delta'+(m+n)\beta \gamma'}\hfill \\
& = & ([m_1]^{\alpha}[m_2]^{(m+n)\beta})^{-\delta'}
([m_1]^{\gamma'}[m_2]^{\delta'})^{(m+n)\beta}=T^{(m+n)\beta},\\[5pt]
%%%
[m_2] 
& = & [m_2]^{-\alpha \delta'+(m+n)\beta \gamma'}\hfill\\
& = & ([m_1]^{\alpha}[m_2]^{(m+n)\beta})^{\gamma'}
([m_1]^{\gamma'}[m_2]^{\delta'})^{-\alpha}=T^{-\alpha},\hfill\\[5pt]
%%%
[l'_1] & = & [m_1]^{\gamma}[l_1]^{\delta}\hfill\\
& = & [m_1]^{\gamma}[m_2]^{(m+n)\delta}
=T^{(m+n)\beta \gamma-(m+n)\alpha \delta}
=T^{m+n}.\hfill
\end{matrix}
\end{eqnarray}
%%%%%
Finally, we attach $V_2$ to $M_1$ to make $M$.
By (\ref{eq:EA}) and (\ref{eq:M2A}) in $H_1(M)$, we have 
%%%%%
\begin{equation*}
[m'_2]=[l_2]=[m_1]^{m+n}=T^{(m+n)^2\beta}=1,
\quad
[l'_2]=[m_2]=T^{-\alpha},
\end{equation*}
%%%%%
and
\begin{equation*}
H_1(M) \cong \langle T\ \vert \ 
T^{(m+n)^2\beta}=1\rangle \cong \mathbb{Z}/(m+n)^2\beta \mathbb{Z}.
\end{equation*}

%%%%%%%%%%%%%%%%%%%%%%%%%%%%%%%%%%%%%%%%%%
%%%                   SubSection 4-2
%%%%%%%%%%%%%%%%%%%%%%%%%%%%%%%%%%%%%%%%%%
\subsection{Reidemeister torsion of $(A_{m,n}; r, 0)$}~\label{ssec:R-tor}
In this subsection,
we compute the Reidemeister torsion of $(A_{m,n}; r, 0)$.
The goal is Lemma~\ref{lem:Rmn}.

\medskip

First, by Surgery formula~II (Lemma~\ref{lem:surgery2}) 
and the results on the first homology in the last subsection, 
we have
%%%%%
\begin{eqnarray}~\label{eq:torM1A}
\tau(M_1) & \doteq &
{\it \Delta}_{A_{m,n}}(T^{(m+n)\beta}, T^{-\alpha})
(T^{m+n}-1)^{-1}.
\end{eqnarray}
%%%%%
By the Alexander polynomial in Theorem~\ref{thm:AmnAlex} and (\ref{eq:torM1A}),
%%%%%
\begin{eqnarray}~\label{eq:torM2A}
%\begin{matrix}
\tau(M_1) 
& \doteq &
%\begin{matrix}
{\displaystyle
\left(\sum_{i=0}^{m+n-1}T^{k_i(m+n)\beta-i\alpha}\right)
(T^{m+n}-1)^{-1}} \nonumber \\
& = & 
{\displaystyle
\frac{
{\displaystyle \sum_{i=0}^{m+n-1}(T^{k_i(m+n)\beta-i\alpha}-T^{-i\alpha})}}
{T^{m+n}-1}+
\frac{{\displaystyle \sum_{i=0}^{m+n-1}T^{-i\alpha}}}
{T^{m+n}-1}}
\bigskip \hfill \nonumber \\
%%%
& = & 
{\displaystyle
\sum_{i=0}^{m+n-1}\left ( T^{-i\alpha}\cdot
\frac{T^{k_i(m+n)\beta}-1}
{T^{m+n}-1}\right ) +
\frac{T^{-(m+n)\alpha}-1}
{T^{m+n}-1} 
\cdot (T^{-\alpha}-1)^{-1}}.
%\end{matrix}
\end{eqnarray}
%%%%%
Note that $(T^{ab} -1)/(T^b-1)$ is a polynomial 
$1 + T^b + T^{2b} + \cdots + T^{(a-1)b}$.

Next, let $d\ge 2$ be a divisor of $(m+n)^2\beta$.
It holds that $\gcd(d, \alpha)=1$.
By Surgery formula~I (Lemma~\ref{lem:surgery1})
and the results on the first homology, 
the Reidemeister--Turaev torsion of $M$ is
%%%%%
\begin{equation}~\label{eq:torM3A}
\tau^{\psi_d}(M)\doteq
\tau^{\rho_d}(M_1)(\zeta_d^{\alpha}-1)^{-1},
\end{equation}
%%%%%
where  $\rho_d:=\psi_d\circ \iota$ is the composite of 
a ring homomorphism $\iota : \mathbb{Z}[H_1(M_1)]\to \mathbb{Z}[H_1(M)]$
induced from the natural inclusion, and a ring homomorphism 
$\psi_d : \mathbb{Z}[H_1(M)]\to \mathbb{Q}(\zeta_d)$
such that $\psi_d(T)=\zeta_d$.

\medskip

We take $d$ as a divisor of $(m+n)$ and
a ring homomorphism 
$\psi'_d : \mathbb{Z}[H_1(M)]\to \mathbb{Q}(\zeta_d)$
such that $\psi'_d(T)=\zeta_d^{-\overline{\alpha}}$
where $\alpha \overline{\alpha}\equiv 1\ (\mathrm{mod}\ \! m+n)$.
Then $d$ is still a divisor of $(m+n)^2\beta$,
and $\zeta_d^{-\overline{\alpha}}$ is still a primitive $d$-th root of unity,
since $\gcd(d, \alpha)=1$. By (\ref{eq:torM2A}) and (\ref{eq:torM3A}), we have
%%%%%
\begin{equation}~\label{eq:torM4A}
\tau^{\psi'_d}(M)\doteq
\left\{\beta(\zeta_d-1)\sum_{i=0}^{m+n-1}k_i\zeta_d^i-\alpha
\right\}
(\zeta_d-1)^{-2}.
\end{equation}
%
%
%%%%%%%%%%%%%%%%%%%%%%%%%%%%%%%%%%%%
\begin{definition}~\label{defini:Rmn}
For a divisor $d\ge 2$ of $(m+n)$, and the primitive $d$-th 
root ($\zeta = \zeta_d$) of unity, we define 
\[
R(m, n):=(\zeta-1)\sum_{i=0}^{m+n-1}k_i\zeta^i.
\]
\end{definition}
%%%%%%%%%%%%%%%%%%%%%%%%%%%%%%%%%%%%
%
%
By (\ref{eq:torM4A}), the Reidemeister--Turaev torsion of $M$
%the result $M$ of the Dehn surgery $(A_{m,n}; \alpha /\beta, 0)$ 
is expressed as 
%%%%%
\begin{equation}~\label{eq:Rmn}
\tau^{\psi'_d}(M)\doteq
\{\beta R(m, n)-\alpha \}(\zeta-1)^{-2}.
\end{equation}
%%%%%
%
%
%%%%%%%%%%%%%%%%%%%%%%%%%%%%%%%%%%%%
\begin{lemma}~\label{lem:R}
\begin{enumerate}
\item[(1)]
$R(m, n)$ is a real number.

\item[(2)]
$\displaystyle
R(m, n)=mn+\frac 12\sum_{i=1}^{m+n-1}(k_{i-1}-k_i)(\zeta^i+\zeta^{-i})$.

\item[(3)]
$\displaystyle
R(m, n)=m(n+1)+\sum_{j=1}^{m-1}(m-s_j)
\left( \zeta^{w_j}+\zeta^{-w_j}\right)$, 
\par

where $s_j$ is defined by 
$s_j := [jn]_m$ for an integer $j$,
see Notation in Subsection~\ref{ssec:Imn}.

\item[(4)]
$\displaystyle
R(m, n)=m(n+1)+\sum_{j=1}^{m-1}(m-j)
\left( \xi^{j}+\xi^{-j}\right)$.

\item[(5)]
$\displaystyle
R(m, n)=
\xi^{-(m-1)}\cdot\left(\frac{\xi^m-1}{\xi-1}\right)^2+mn
=\left|\frac{\xi^m-1}{\xi-1}\right|^2+mn$,
\par

where $\xi=\zeta^{\overline{m}}$ with 
$m\overline{m} \equiv 1 \mod (m+n)$.
\end{enumerate}
\end{lemma}
%%%%%%%%%%%%%%%%%%%%%%%%%%%%%%%%%%%%
%
%
\medskip \noindent
{\bf Proof} \ 
(1)\ The complex conjugate $\overline{R(m, n)}$ of $R(m, n)$ is
%%%%%
\begin{eqnarray*}
\overline{R(m, n)} & = & 
(\zeta^{-1}-1)\sum_{i=0}^{m+n-1}k_i\zeta^{-i}
=(\zeta-1)\sum_{i=0}^{m+n-1}-k_i\zeta^{-i-1}\\
& = & (\zeta-1)\sum_{i=0}^{m+n-1}(mn-k_i)\zeta^{m+n-i-1}\\
& = & (\zeta-1)\sum_{i=0}^{m+n-1}k_{m+n-i-1}\zeta^{m+n-i-1}\\
& = & (\zeta-1)\sum_{i=0}^{m+n-1}k_i\zeta^i=R(m, n)
\end{eqnarray*}
%%%%%
by Proposition~\ref{prop:Imn} (1) and 
the equality $1+\zeta+\cdots +\zeta^{m+n-1}=0$.
%Hence $R(m, n)$ is a real number.

\bigskip
\noindent
(2)\ First, we have
%%%%%
\begin{eqnarray*}
R(m, n) & = & \sum_{i=0}^{m+n-1}k_i(\zeta^{i+1}-\zeta^i)
=\sum_{i=1}^{m+n}k_{i-1}\zeta^i-\sum_{i=0}^{m+n-1}k_i\zeta^i\\
& = & mn+\sum_{i=1}^{m+n-1}(k_{i-1}-k_i)\zeta^i.
\end{eqnarray*}
%%%%%
By the proof of (1), we have
%%%%%
\begin{eqnarray*}
R(m, n) & = & \frac 12\left\{R(m, n)+\overline{R(m, n)}\right\}\\
& = & mn+\frac 12\sum_{i=1}^{m+n-1}(k_{i-1}-k_i)(\zeta^i+\zeta^{-i}).
\end{eqnarray*}
%%%%%

Before the proof of (3), (4) and (5), we prove the following: 
\par

\medskip \noindent
{\bf Claim} (Property on $s_j$) \ 
{\sl
\begin{enumerate}
\item[(a)]
The map from $j$ to $s_j$($= [jn]_m$) is a bijection on
$\{1, 2, \ldots, m-1 \}$ to itself.

\item[(b)] It holds that $s_j + s_{m-j} = m$.

\item[(c)]
There exists a unique element $h$ in $\{1,2, \ldots, m-1 \}$ such that 
$\gcd(h, m)=1$ and $s_h=1$.
It holds that $w_h\equiv -\overline{m}\ (\mathrm{mod}\ \! m+n)$.

\item[(d)]
For the same $h$ in (c) and each element $a$ in $\{1,2, \ldots, m-1 \}$,
we have 
\[
s_{ah}=a\quad
\mbox{and}\quad
w_{ah}\equiv aw_h\quad (\mathrm{mod}\ \! m+n),\]
where we regard $w_{ah}$ as $w_j$ with $j = [ah]_m$, precisely.
%where we regard $s_{ah}$ as $s_j$ ($w_{ah}$ as $w_j$, respectively)
%with $j = [ah]_m$, precisely.
\end{enumerate}
}

\noindent
{\bf Proof of Claim} \ 
(a)\ 
The map is induced by the multiplication of $n$
(i.e., $j \mapsto jn$) over $(\mathbb{Z}/m\mathbb{Z})\backslash \{0\}$.
It is a bijection, since $\gcd(m,n)=1$.

\medskip \noindent
(b)\ It is easy to see.

\medskip \noindent
(c)\ 
By (a), there exists a unique element $h$ in $\{1,2, \ldots, m-1 \}$ 
such that $s_h=1$.
In fact, it holds that $h\equiv \overline{n}\ (\mathrm{mod}\ \! m)$. 
We have $\gcd(h, m)=1$.
The second half is shown by $-mw_h=-m\left (h+\left[ \frac{hn}{m}\right] \right ) 
\equiv hn - m\left[ \frac{hn}{m}\right]
= [hn]_m = s_h =1\quad (\mathrm{mod}\ \! m+n)$,
see Proposition~\ref{prop:Imn} (2).

\medskip \noindent
(d)\ Since $s_{ah}\equiv ahn\equiv a\ (\mathrm{mod}\ \! m)$,
we have $s_{ah}=a$, for $a$ in $\{1,2, \ldots, m-1 \}$.
Since 
$mw_{ah}=mah+m\left[ \frac{ahn}{m}\right]
\equiv -\left ( ahn - m\left[ \frac{ahn}{m}\right] \right )
=-[ahn]_m =-s_{ah}= -a\quad (\mathrm{mod}\ \! m+n)$.
\qed

\bigskip \noindent
(3) We go back to the expression (2). 
%We divide the summation over $\{ 1,2, \ldots, m+n-1\}$
%into the sub-summations over $M, R$ and $L$
%as the set 
We divide the set $\{ 1,2, \ldots, m+n-1\}$
of indices of $k_i$, into $M, R$ and $L$
according to whether 
$k_{i-1}$ and/or $k_i$ belongs to $m\mathbb{Z} \cap {\frak I}(m, n)$.
\[
\{ 1,2, \ldots, m+n-1\} = M \cup R \cup L \quad 
\textrm{ (a disjoint union)}
\]
\begin{center}
\begin{tabular}{lll}
Definition of the subset & parameter $j$ & $k_{i-1} -k_i$  \\
\hline
$M := \{ i \vert k_{i-1} \in m\mathbb{Z}$ and $k_i \in m\mathbb{Z} \}$ 
& \quad --- & \quad $-m$ \\
$R := \{ i \vert k_{i-1} \in m\mathbb{Z}$ and $k_i \not \in m\mathbb{Z} \}$ 
& $k_i =jn$
& \quad $m \left[ \frac{jn}m \right ] -jn$ \\
$L := \{ i \vert k_{i-1} \not\in m\mathbb{Z}$ and $k_i \in m\mathbb{Z} \}$ 
& $k_{i-1} =jn$ 
& \quad $jn - m\left (\left[ \frac{jn}m \right ]+1\right )$ \\[5pt]
\hline
\end{tabular}
\end{center}
Note that the case both $k_{i-1} \not\in m\mathbb{Z}$ and $k_i \not\in m\mathbb{Z}$
(in other words, the case that both $k_{i-1}$ and $k_i$ belong to 
$n\mathbb{Z}$) never occur, since $m < n$ 
(Section~\ref{sec:intro}). 

For each $i \in R$, there exists $j$ with $1 \le j < m$
such that $k_i =jn$, equivalently $i = w_j $. 
Then, by Proposition~\ref{prop:Imn} (4),
$i' := w_{m-j}+1$ belongs to $L$ and it holds that
$i' = (m+n-1 - w_j)+1 = m+n-i$.
The correspondence between $i \in R$ and $i' \in L$ above is one to one. 
It also holds that
$\zeta^{i'} = \zeta^{w_{m-j}+1} = \zeta^{m+n-i} 
= \zeta^{-i}$ and 
%%%%%
\begin{align*}
k_{i-1} - k_i 
& =m\left[ \frac{jn}m \right ] - jn \ 
= -[jn]_m \ = -s_j, \\
k_{i'-1} - k_{i'} 
& = (m-j)n -m \left ( \left[ \frac{(m-j)n}m \right ] +1 \right) \\
& = [(m-j)n]_m - m = s_{m-j} - m = -s_j,
%& = [(m-j)n]_m - m = [-jn]_m - m = -[jn]_m.
\end{align*}
%%%%%
by Claim (b). Thus 
%%%%%
\begin{align*}
& (k_{i-1} - k_i) (\zeta^i + \zeta^{-i}) +
(k_{i'-1} - k_{i'}) (\zeta^{i'} + \zeta^{-i'}) \\
& 
= -2s_j(\zeta^i + \zeta^{-i}) \\
&
= -2s_j (\zeta^{w_j} + \zeta^{-w_j}),
\end{align*}
%%%%%
and 
%%%%%
\begin{align*}
& \sum_{i \in R}(k_{i-1}-k_i+m)(\zeta^i+\zeta^{-i}) \ + \ 
\sum_{i \in L}(k_{i-1}-k_i+m)(\zeta^i+\zeta^{-i}) \\
& = \sum_{i \in R} \left \{
(k_{i-1} - k_i + m) (\zeta^i + \zeta^{-i}) + 
(k_{i'-1} - k_{i'} + m) (\zeta^{i'} + \zeta^{-i'})\right \}  \\
& = 2 \sum_{j=1}^{m-1} (m-s_j) (\zeta^{w_j} + \zeta^{-w_j}).
\end{align*}
Thus, using $1 + \zeta + \zeta^2 + \cdots + \zeta^{m+n-1} =0$,
%%%%%
\begin{align*}
R(m, n)
& =  mn+\frac 12\sum_{i=1}^{m+n-1}(k_{i-1}-k_i)(\zeta^i+\zeta^{-i})\\
& =  mn+\frac 12\sum_{i=1}^{m+n-1}(k_{i-1}-k_i)(\zeta^i+\zeta^{-i})
+ \frac 12\sum_{i=0}^{m+n-1}m(\zeta^i+\zeta^{-i})\\
& =  mn + m +\frac 12\sum_{i=1}^{m+n-1}(k_{i-1}-k_i+m)(\zeta^i+\zeta^{-i}) \\
& =  m(n+1)+\frac 12\sum_{i \in M \cup R \cup L}(k_{i-1}-k_i+m)(\zeta^i+\zeta^{-i})\\
& =  m(n+1)+ \sum_{j=1}^{m-1} (m-s_j)(\zeta^{w_j}+\zeta^{-w_j}).
\end{align*}
%%%%%

\bigskip \noindent
(4) We take $h$($= \overline{n} \mod m$) in Claim (c). By Claim (d), we have
%%%%%
\begin{align*}
R(m, n)
& =  m(n+1)+ \sum_{j=1}^{m-1} (m-s_j)(\zeta^{w_j}+\zeta^{-w_j}) \\
& =  m(n+1)+ \sum_{a=1}^{m-1} (m-s_{ah})(\zeta^{w_{ah}}+\zeta^{-w_{ah}}) \\
& =  m(n+1)+ \sum_{a=1}^{m-1} (m-a)(\zeta^{-a\overline{m}}+\zeta^{a\overline{m}}) \\
& =  m(n+1)+ \sum_{a=1}^{m-1} (m-a)(\xi^a+\xi^{-a}),
\end{align*}
%%%%%
where we set $\xi=\zeta^{\overline{m}}$,
which is also a $d$-th primitive root of unity.

\bigskip \noindent
(5) Using elementary calculus
\[
\sum_{a=1}^{m-1} (m-a) x^a 
\ = \ \dfrac{x^m - mx -1}{x-1} + \dfrac{x^m-1}{(x-1)^2}
\ = \ x \dfrac{x^m-1}{(x-1)^2} - m \dfrac{x}{x-1},
%= \dfrac{x^m - mx -1}{x-1} + \dfrac{x^m-1}{(x-1)^2}
%\sum_{i=1}^{m-1} ix^i = \dfrac{(m-1)x^m +1}{x-1} - \dfrac{x^m-1}{(x-1)^2},
\]
we have
\[
R(m, n)=\xi^{-(m-1)}\cdot\left(\frac{\xi^m-1}{\xi-1}\right)^2+mn.
\]
The proof of Lemma~\ref{lem:R} is completed. 
\qed 
\par

\medskip
The result of this subsection is summarized as:
%
%
%%%%%%%%%%%%%%%%%%%%%%%%%%%%%%%%%%%%
\begin{lemma}~\label{lem:Rmn}
Let $M = (A_{m,n}; \alpha/\beta,0)$,
and $d\ge 2$ a divisor of $m+n$.
Then the Reidemeister--Turaev torsion of $M$ related with 
$\psi'_d(T) = \zeta_d^{-\bar{\alpha}}$ is 
%%%%%
\begin{eqnarray*}
\tau^{\psi'_d}(M)
& \doteq &
\xi^{-m} \left\{\beta R(m,n) - \alpha \right\}
(\xi^m-1)^{-2} \\
& = &
\xi^{-m} 
\left\{\beta \left| \frac{\xi^m-1}{\xi-1}\right|^2 -(\alpha-mn\beta) \right\}
(\xi^m-1)^{-2},
\end{eqnarray*}
%%%%%
where $\xi=\zeta^{\overline{m}}$ (thus $\xi^m=\zeta$) is a primitive $d$-th root
of unity.
\end{lemma}
%%%%%%%%%%%%%%%%%%%%%%%%%%%%%%%%%%%%
%
%

%%%%%%%%%%%%%%%%%%%%%%%%%%%%%%%%%%%%%%%%%%
%%%                   SubSection 4-3
%%%%%%%%%%%%%%%%%%%%%%%%%%%%%%%%%%%%%%%%%%
\subsection{Necessary conditions}~\label{ssec:key}
%{Key lemmas}
We study some necessary conditions for
$(A_{m,n}; \alpha/\beta, 0)$ to be a lens space 
by the Reidemeister--Turaev torsions.
%
%
%%%%%%%%%%%%%%%%%%%%%%%%%%%%%%%%%%%%
\begin{lemma}~\label{lem:Equality}
Suppose that $(A_{m,n}; \alpha/\beta, 0)$ is a lens space.
Then there exists integers $i$ and $j$ such that
$\gcd(i, m+n)=\gcd(j, m+n)=1$ and
%%%%%
\begin{equation}~\label{eq:torR}
\left\{\beta R(m, n)-\alpha \right\}
(\xi^m-1)^{-2}\doteq (\xi^i-1)^{-1}(\xi^j-1)^{-1},
\end{equation}
%%%%%
equivalently, 
%%%%% i.e., by Lemma~\ref{lem:R} (4),
\begin{equation}~\label{eq:Rijm}
\left\{\beta \left| \frac{\xi^m-1}{\xi-1}\right|^2 - \alpha' \right\}
(\xi^m-1)^{-2}
\doteq
\frac{1}{(\xi^i-1)(\xi^j-1)},
\end{equation}
%%%%%
where $\alpha' = \alpha - mn\beta$, and  
$\xi$ is a primitive $d$-th root of unity.
\end{lemma}
%%%%%%%%%%%%%%%%%%%%%%%%%%%%%%%%%%%%
%
%
The equalities (\ref{eq:torR}), (\ref{eq:Rijm}) correspond to 
two expressions of $\tau^{\psi'_d}(M)$ in 
Lemma~\ref{lem:Rmn}.
%
%
%%%%%%%%%%%%%%%%%%%%%%%%%%%%%%%%%%%%
\begin{lemma}~\label{lem:normR}
Let $d\ge 2$ be a divisor of $m+n$.
Suppose that $(A_{m,n}; \alpha/\beta, 0)$ is a lens space.
Then we have
\begin{enumerate}
\item[(1)] 
$|N_d\left( \beta R(m, n)-\alpha \right)|=1$, 
where $N_d$ is the $d$-norm, see Subsection~\ref{ssec:norm}.

\item[(2)] $\alpha'=\alpha - mn\beta \ge 0$.
\end{enumerate}
\end{lemma}
%%%%%%%%%%%%%%%%%%%%%%%%%%%%%%%%%%%%
%
%
\medskip \noindent
{\bf Proof} \ 
(1) We take the $d$-norm of the equality (\ref{eq:torR}).
Since $N_d(\xi^m-1)=N_d(\xi^i-1)=N_d(\xi^j-1)\ne 0$ 
by Proposition~\ref{pr:norm} (1) and Lemma~\ref{lem:cyclotomic},
we have the result. 

\par \medskip
(2) Suppose that the integer $\alpha'=\alpha-mn\beta<0$.
Then we have 
\[
\beta R(m, n)-\alpha
=\beta \left| \frac{\xi^m-1}{\xi-1}\right|^2-\alpha' > 1,
\]
hence $|N_d\left( \beta R(m, n)-\alpha \right)|>1$.
By (1), we have the result.
\qed 
\par

\medskip

Fixing the combinatorial Euler structure,
we will regard (\ref{eq:Rijm}) as a control equivalence
of the value sequences of degree $m+n$,
in the sense of Lemma \ref{lem:lensEx} (1).
Note that the first factor in the left-hand side is a real value.
On the right-hand side, we have to control $(i,j)$ to use 
Lemma~\ref{lem:Euler2} or Lemma~\ref{lem:Euler}.
\par

\medskip
\noindent
{\bf Conditions on $(i,j)$ and $(e,f)$} \ 
We can take $i$ and $j$ satisfying $1\le i \le j\le (m+n-1)/2$.
%%%%%
If $i+j$ is odd (then $m+n$ is odd),
then we replace $j$ with $m+n-j$ and denote it by $j$ again.
Then, as a condition of $(i,j)$, we may assume 
%%%%%
\begin{eqnarray}~\label{eq:ij-cond}
1\le i\le j\le m+n-1, \ 
2\le i+j\le m+n-1 \textrm{ and $i+j$ is even.}
\end{eqnarray}
%%%%%
From now on, we regard the equality (\ref{eq:Rijm})
as a controll equivalence between the 
\underline{real} value sequences
\begin{equation}~\label{eq:Rijm2}
u \xi^{-\frac{i+j}2}
\left\{\beta \left| \frac{\xi^m-1}{\xi-1}\right|^2 - \alpha' \right\}
(\xi^i-1)(\xi^j-1)
\ = \ 
\xi^{-m} (\xi^m-1)^2, 
\end{equation}
where $u=\pm 1$, or $\pm \xi^{\frac{m+n}{2}}$ only if $m+n$ is even,
by Lemma~\ref{lem:control}.
We define the integers
%%%%%
\begin{equation}~\label{eq:ef}
e :=\frac{j-i}{2} \ \mbox{and} \ 
f :=\frac{i+j}{2}. \quad
\mbox{They satisfies }\ 0\le e<f\le (m+n-1)/2.
\end{equation}
%%%%%
Using $(e,f)$, we can deform (\ref{eq:Rijm2})
as
%%%%%
\begin{eqnarray}~\label{eq:Fz=Gz}
\begin{matrix}
u \left\{ \beta \left( \xi^m+\xi^{-m}\right) -\alpha' \left(
\xi+\xi^{-1}\right)
+2\left( \alpha' -\beta \right) \right\} % corrected
\left\{ \left( \xi^f+\xi^{-f}\right)- \left( \xi^e+\xi^{-e}\right) \right\}
\hfill \medskip\\
\ = 
\left( \xi^{m+1}+\xi^{-(m+1)}\right)
-2\left( \xi^m+\xi^{-m}\right)
+\left( \xi^{m-1}+\xi^{-(m-1)}\right)
-2 \left( \xi+\xi^{-1}\right) +4. \hfill
\end{matrix}
\end{eqnarray}
We define two symmetric Laurent polynomials
%%%%%
\begin{eqnarray}~\label{eq:FG}
\begin{aligned}
F(t) & = &
\left\{ \beta \left( t^m+t^{-m}\right) -\alpha' \left( t+t^{-1}\right)
+2\left( \alpha' -\beta \right) \right\} 
\left\{ \left( t^f+t^{-f}\right)- \left( t^e+t^{-e}\right) \right\},
\medskip\\
G(t) & = & \left( t^{m+1}+t^{-(m+1)}\right)
-2\left( t^m+t^{-m}\right)
+\left( t^{m-1}+t^{-(m-1)}\right)
-2 \left( t+t^{-1}\right) + 4,
\end{aligned}
\end{eqnarray}
%%%%%
then (\ref{eq:Rijm2}) means that 
two real value sequences (of degree $m+n$)
induced by $F(t)$ and $G(t)$
are control equivalent, see Subsection~\ref{ssec:Euler}.
Note that $F(1)=G(1)=0$.
By Lemma~\ref{lem:Euler}, (\ref{eq:Fz=Gz}) lifts to a congruence
of the symmetric Laurent polynomials.
%%%%%
%
%
%%%%%%%%%%%%%%%%%%%%%%%%%%%%%%%%%%%%
\begin{lemma}[Necessary condition]~\label{lem:NesC}
Suppose that $(A_{m,n}; \alpha/\beta, 0)$ is a lens space.
We set $\alpha'=\alpha-mn\beta$.
Then there exist integers $e$ and $f$ such that
$0\le e< f \le (m+n-1)/2$, and the following congruence holds:
\begin{enumerate}
\item[(a)]
If $m+n$ is odd,
we have
$F(t)\equiv \pm G(t)\ (\mathrm{mod}\ \! t^{m+n}-1)$.

\item[(b)]
If $m+n$ is even,
we have
$F(t)\equiv \pm G(t)$ or
$F(t)\equiv \pm t^{\frac{m+n}{2}}G(t)\ (\mathrm{mod}\ \! t^{m+n}-1)$.
\end{enumerate}
\end{lemma}
%%%%%%%%%%%%%%%%%%%%%%%%%%%%%%%%%%%%
%
%
If $m+n$ is even, then $\textrm{span}(G(t)) = 2(m+1) \le m+n$, 
since the pair $m, n$ is coprime, thus both are odd and $m+2 \le n$.
Furthermore,
the congruence also induces an identity
%%%%%
\begin{equation}~\label{eq:F=G}
\textrm{(i) } \textrm{red}(F(t)) = \pm G(t) 
\quad \textrm{ or \quad
(ii) } \textrm{red}( t^{\frac{m+n}{2}}F(t)) = \pm G(t),
\end{equation}
%%%%%
as in the second half of Lemma~\ref{lem:Euler}.
We will regard it as \underline{an equation} of $(e,f)$
on the surgery coefficient $\alpha/\beta$ for $M = (A_{m,n}; \alpha/\beta, 0)$ to be a lens space:
Suppose that $M$ is a lens space, 
then there exists a solution $(e,f)$ of the equation.
We mainly use its contraposition:
If the equation has no solution $(e,f)$,
then $M$ is not a lens space. 
The case (b) looks troublesome. To prove that $M$ is not a lens space,
we have to show that 
neither (i) nor (ii)
has a solution.
Fortunately, we only have to show one of them.
%
%
%%%%%%%%%%%%%%%%%%%%%%%%%%%%%%%%%%%%
\begin{remark}~\label{rem:Nes(b)}
{\rm 
In either case $m+n$ is odd or even,
to prove that
$(A_{m,n}; \alpha/\beta, 0)$ is not a lens space,
it is sufficient to show that 
$\textrm{red}(F(t)) = \pm G(t)$ has no solution $(e,f)$,
because we can prove the following.
}
\end{remark}
%%%%%%%%%%%%%%%%%%%%%%%%%%%%%%%%%%%%
%
%
%%%%%%%%%%%%%%%%%%%%%%%%%%%%%%%%%%%%
\begin{lemma}~\label{lem:Nes(b)}
In the case (b) $m+n$ is even in Lemma~\ref{lem:NesC},
if the equation (i) $\textrm{red}(F(t)) = \pm G(t)$ in (\ref{eq:F=G})
has a solution, 
the other equation (ii) $\textrm{red}( t^{\frac{m+n}{2}}F(t)) = \mp G(t)$
has a solution, and vice versa.
\end{lemma}
%%%%%%%%%%%%%%%%%%%%%%%%%%%%%%%%%%%%
%
%
\medskip \noindent
{\bf Proof} \ 
We concentrate on the factor 
$\left\{ \left( t^f+t^{-f}\right)- \left( t^e+t^{-e}\right) \right\}$
of $F(t)$. We transform $(e,f)$ to $(e',f')$ by
\begin{equation*}~\label{eq:ij'}
e' =\frac{m+n}{2}-f\quad \mbox{and}\quad
f' =\frac{m+n}{2}-e,
\end{equation*}
which satisfies the same condition $0\le e' < f' \le (m+n-1)/2$ with (\ref{eq:ef}).
For a solution $(e,f)$ of the equation $\textrm{red}(F(t)) = \pm G(t)$,
its transformation $(e',f')$ is a solution of
$\textrm{red}( t^{\frac{m+n}{2}}F(t)) = \mp G(t)$, and vice versa.
\qed 
\par

\medskip

In Subsection 5.3, we will prove that
$(e, f)=(0, 1)$ with $\alpha' = 0$ is the only solution for the equation
$\textrm{red}(F(t)) = \pm G(t)$
in general cases (see Lemma \ref{lem:special} (1) below).
Note that $\alpha' =0$ implies $\alpha/\beta=mn$, which is related 
to the lens space surgery in Theorem~\ref{thm:(mn)}.

Using the expression of $R(m,n)$ in Lemma~\ref{lem:R} (4),
we can prove the following.
%
%
%%%%%%%%%%%%%%%%%%%%%%%%%%%%%%%%%%%%
\begin{lemma}~\label{lem:special}
\begin{enumerate}
\item[(1)]
%The equation (\ref{eq:Rijm}) has a root $\alpha/\beta=mn$
%(i.e.\ $\alpha'=0$) if and only if $(e, f)=(0, 1)$ (i.e.\ $i=j=1$).
The condition $\alpha/\beta=mn$ (i.e.\ $\alpha'=0$) is equivalent to $(e, f)=(0, 1)$ (i.e.\ $i=j=1$).

\item[(2)]
In (\ref{eq:Rijm}), if $(e, f)=(0, m)$ (i.e.\ $i=j=m$), then
there is no root for $\alpha/\beta$.

\end{enumerate}
\end{lemma}
%%%%%%%%%%%%%%%%%%%%%%%%%%%%%%%%%%%%
%
%
\medskip \noindent
{\bf Proof} \ 
(1)\ Suppose that $\alpha/\beta=mn$.
Then we have $\tau^{\psi_d'}(M)\doteq (\xi-1)^{-2}$
by Lemma~\ref{lem:Rmn}, and the equality (\ref{eq:Rijm})
becomes 
\[
\dfrac{1}{(\xi -1)^2} \ \dot= \ \dfrac1{(\xi^i -1)(\xi^j -1)}.
\]
We have $i=j=1$ by the Franz lemma (Lemma~\ref{lem:Franz}).

\medskip

Conversely suppose that $i=j=1$.
Then the equality (\ref{eq:torR}) can be deformed into
\[
\beta R(m,n) - \alpha 
= u
\xi^{-(m-1)} \left ( \frac{\xi^m-1}{\xi-1}\right )^2
\]
where $u=\pm 1$ or $\pm \xi^{\frac{m+n}{2}}$.
Using the expression of $R(m,n)$ in Lemma~\ref{lem:R} (4)
and the calculus in the proof of Lemma~\ref{lem:R} (5), we have
%%%%%
\begin{equation*}
\beta m(n+1)-\alpha 
+\sum_{j=1}^{m-1}\beta (m-j)\left( \xi^j+\xi^{-j}\right)
=u \left\{
m+\sum_{j=1}^{m-1}(m-j)\left( \xi^j+\xi^{-j}\right)
\right\}.
\end{equation*}
%%%%%
By taking the symmetric polynomial lift and 
Lemma~\ref{lem:Euler2} (Note that the span is 
$2(m-1) \le 2([(m+n)/2]-1)$), 
the case $u=\pm \xi^{\frac{m+n}{2}}$ does not occur,
by Lemma~\ref{lem:symeven}.
We have $\beta=1$ and $\alpha=mn$.
\bigskip

\noindent
(2)\ Suppose that $i=j=m$.
Then the equality (\ref{eq:torR}) is deformed into
\[
\beta R(m,n) - \alpha 
= u
\]
where $u=\pm 1$ or $\pm \xi^{\frac{m+n}{2}}$.
By the same method with above, we have
%%%%%
\begin{equation*}
\beta m(n+1)-\alpha 
+\sum_{j=1}^{m-1}\beta (m-j)\left( \xi^j+\xi^{-j}\right)=u.
\end{equation*}
%%%%%
By Lemma~\ref{lem:Euler2} and Lemma~\ref{lem:symeven},
we have $\beta =0$. Hence there is no root for $\alpha/\beta$.
\qed
%
%
%%%%%%%%%%%%%%%%%%%%%%%%%%%%%%%%%%%%
\begin{lemma}~\label{lem:gcd}
Suppose that $\alpha' >0$.
\begin{enumerate}
\item[(1)]
$\gcd (\alpha', \beta) =1$.
\item[(2)] If $\alpha' = \beta$, then the congruence $F(t) \equiv \pm G(t)$ (mod $t^{m+n} -1$)
has a unique solution $(m,n)=(2,3)$.
\end{enumerate}
\end{lemma}
%%%%%%%%%%%%%%%%%%%%%%%%%%%%%%%%%%%%
%
%%%%%%%%%%%%%%%%%%%%%%%%%%%%%%%%%%%%
%
%
\medskip \noindent
{\bf Proof} \ 
(1)\ 
By the Euclidean algorithm, we have 
\[
\gcd (\alpha', \beta) 
=
\gcd (\alpha - mn \beta, \beta)
= 
\gcd (\alpha, \beta)
=1.
\]

\noindent
(2)\ 
Suppose $\alpha'=\beta$.
Then, by (1), we have $\alpha'=\beta=1$ and the polynomials (\ref{eq:FG})
are
\begin{eqnarray*}
F(t) & = & t^{-m-f}(t^{m+1}-1)(t^{m-1}-1)(t^i-1)(t^j-1), \\
G(t) & = & t^{-m-1}(t^m-1)^2(t-1)^2.
\end{eqnarray*}
The congruence $F(t) \equiv \pm G(t)$ (mod $t^{m+n} -1$) implies
\[
F(\zeta) \doteq G(\zeta),
\]
where $\zeta$ is a primitive $(m+n)$-th root of unity.
Suppose $\gcd(m-1, m+n)\ge 2$ or $\gcd(m+1, m+n)\ge 2$.
Then the left-hand side of the equation above is $0$ for some $d$.
Hence we have $\gcd(m-1, m+n)=1$ and $\gcd(m+1, m+n)=1$.
By the Franz lemma \cite{Fz} (Lemma~\ref{lem:Franz}), 
we have $$\{\pm (m-1), \pm (m+1), \pm i, \pm j\ (\mathrm{mod}\ \! m+n)\}
=\{\pm 1, \pm 1, \pm m, \pm m\ (\mathrm{mod}\ \! m+n)\}$$
as multiple sets.
It has a unique solution $(m, n)=(2, 3)$ with $(i,j)=(1,3)$.
\qed

%%%%%%%%%%%%%%%%%%%%%%%%%%%%%%%%%%%%%%%%%%
%%%                   Section 5
%%%%%%%%%%%%%%%%%%%%%%%%%%%%%%%%%%%%%%%%%%
\section{Proof of Theorem \ref{thm:MT1}\ 
(Lens space surgeries along $A_{m,n}$)}~\label{sec:proofA}
The ``if part" of Theorem~\ref{thm:MT1} (1) 
follows from Theorem~\ref{thm:(mn)},
thus our main purpose is
to prove the ``only if part".
We study the condition on $\alpha/\beta$ for the 
equations (\ref{eq:Rijm2}) or (\ref{eq:F=G})(i)
has a solution $(i,j)$ or $(e,f)$, respectively.
Our proof is divided into three cases:
$m=2$ (Subsection~\ref{ssec:m=2}), 
$n=m+1$ (Subsection~\ref{ssec:n=m+1}),
and the general case 
where $m\ge 3$ and $n\ge m+2$ (Subsection~\ref{ssec:ngem+2}).
Note that the first two cases contains the exceptional case $(A_{2,3};7,0)$.
In Subsection~\ref{ssec:type} and~\ref{ssec:move},
we prove Theorem~\ref{thm:MT1} (2)
by using the values of the Reidemeister torsions.
We also use Kirby moves. 
\par

\medskip

To make expressions of symmetric Laurent polynomials short,
we use the notation $\langle t^x \rangle = t^x + t^{-x}$
for any integer $x$.
We regard $\langle t^x \rangle$ as $\langle t^{-x} \rangle$
if $x <0$, and $\langle t^0 \rangle= 2$.
For the terminologies ^^ ^^ reduce, reduction (denoted by $\textrm{red}(P(t))$)"
of symmetric Laurent polynomials, see 
Definition~\ref{defini:symLPoly}.

%%%%%%%%%%%%%%%%%%%%%%%%%%%%%%%%%%%%%%%%%%
%%%                   SubSection 5-1
%%%%%%%%%%%%%%%%%%%%%%%%%%%%%%%%%%%%%%%%%%
\subsection{The case $m=2$}~\label{ssec:m=2}
In this case, $n$ and $m+n =n+2$ are odd.
Let $\xi$ denote any $d$-th root of unity, where
$d$ is a divisor of $n+2$ with $d \ge 2$.
We have 
$R(2,n) = \vert \xi +1 \vert^2 +mn = \xi + \xi^{-1} + mn + 2$
by Lemma~\ref{lem:R} (5),
thus the equation (\ref{eq:Rijm2}), divided by $\xi^{-1}(\xi -1)^2$
as a value sequence, becomes
%%%%%
\begin{equation}~\label{eq:m=2tor}
\xi^{-\frac{i+j-2}{2}}\cdot
\left\{ \beta \left( \xi+\xi^{-1}\right)-\alpha'' \right\}
\frac{(\xi^i-1)(\xi^j-1)}{(\xi-1)^2}
=\eta \xi^{-1}\cdot \frac{(\xi^2-1)^2}{(\xi-1)^2}
\end{equation}
%%%%%
where $\alpha'' = \alpha' - 2\beta = \alpha - 2(n+1)\beta$, $\eta=\pm 1$
and
$(i,j)$ satisfies the condition (\ref{eq:ij-cond}) in the last section. 
We regard (\ref{eq:m=2tor}) as an equality between 
real value sequences,
defined in Subsection~\ref{ssec:Euler}.

\medskip
\noindent
(1) \ The case $2\le i+j < n+1$. Note that $i+j$ is even, see (\ref{eq:ij-cond}).

By Lemma~\ref{lem:Euler2}, the equalities (\ref{eq:m=2tor}) 
lift to a congruence and 
%%%%%
\begin{equation*}
t^{-\frac{i+j-2}{2}}\cdot
\left\{ \beta \left( t+t^{-1}\right)-\alpha'' \right\}
\frac{(t^i-1)(t^j-1)}{(t-1)^2}
=\eta t^{-1}\cdot \frac{(t^2-1)^2}{(t-1)^2}.
\end{equation*}
%%%%%
Note that $(t^x -1)/(t-1)$ is a polynomial for an integer $x$
and that both hand sides are symmetric Laurent polynomials.
The span of the left-hand side 
is $i+j \le 2([(m+n)/2]-1) = 2n$.
From $\beta >0$, we have $i = j =1$,
$\beta \left( t+t^{-1}\right)-\alpha''= t^{-1}(t+1)^2$,
$\beta=1$, $\alpha'' = -2$, and $\alpha=2n$.

\medskip
\noindent
(2) \ The case $i+j=n+1$ ($=2f$, see (\ref{eq:ef}) in Subsection~\ref{ssec:key}).

Then $e = (j-i)/2$ is an integer satisfying $1\le e<\frac{n+1}{2}$.
The equation (\ref{eq:Rijm2}) is
%%%%%
\begin{equation*}
\xi^{-\frac{n+1}{2}}\cdot
\left\{ \beta \left( \xi+\xi^{-1}\right)-\alpha'' \right\}
( \xi^i-1)(\xi^j-1)
=\eta \xi^{-2}\cdot (\xi^2-1)^2.
\end{equation*}
%%%%%
By Lemma~\ref{lem:Euler}, it lifts to
%%%%%
\begin{eqnarray}~\label{eq:m=2poly}
\begin{matrix}
& & (\beta-\alpha'')\left( t^{\frac{n+1}{2}}+t^{-\frac{n+1}{2}}\right)
+\beta \left( t^{\frac{n-1}{2}}+t^{-\frac{n-1}{2}}\right)\hfill\medskip\\
& & -\beta \left( t^{e+1}+t^{-(e+1)}\right)
+\alpha''\left( t^e+t^{-e}\right)
-\beta \left( t^{e-1}+t^{-(e-1)}\right)\hfill\medskip\\
& = & \eta \left\{ \left( t^2+t^{-2}\right)-2\right\}, \hfill
\end{matrix}
\end{eqnarray}
which is, using notations $\langle t^x \rangle = t^x + t^{-x}$,
\[
 (\beta-\alpha'') \langle t^{\frac{n+1}{2}} \rangle
+\beta \langle t^{\frac{n-1}{2}}\rangle
- \beta \langle t^{e+1} \rangle
+ \alpha''\langle  t^e \rangle
- \beta \langle t^{e-1} \rangle
= \eta (\langle  t^2 \rangle-2). \hfill
\]
%%%%%
Here we used $\langle t^{\frac{n+3}{2}} \rangle 
\equiv \langle t^{\frac{n+1}{2}} \rangle$
(mod $t^{n+2}-1$).
Note that the span of the left-hand side is at most $n+1 = 2[(m+n)/2]$.
We have $e=1$, $\eta = \beta =1$ and $n=3$.
Then (\ref{eq:m=2poly}) is deformed into
%%%%%
\begin{equation*}
-\alpha'' \left( t^2+t^{-2}\right)
+(1+\alpha'')\left( t+t^{-1}\right)-2=\left( t^2+t^{-2}\right)-2.
\end{equation*}
%%%%%
Hence we have $\alpha''=-1$ and $\alpha=7$.
In this case, $(i,j)=(1,3)$.
This corresponds to the case $(A_{2,3}; 7, 0)$, and it is a lens space,
see Subsection~\ref{ssec:move}.

%%%%%%%%%%%%%%%%%%%%%%%%%%%%%%%%%%%%%%%%%%
%%%                   SubSection 5-2
%%%%%%%%%%%%%%%%%%%%%%%%%%%%%%%%%%%%%%%%%%
\subsection{The case $n=m+1$}~\label{ssec:n=m+1}
In this case, $m+n = 2m+1$ is odd.
Let $\xi$ denote any $d$-th root of unity, where
$d$ is a divisor of $2m+1$ with $d \ge 2$.
We use $\zeta=\xi^m$ as in Lemma~\ref{lem:Rmn},
then $\xi=\zeta^{-2}$. 
The Reidemeister torsion in Lemma~\ref{lem:Rmn} 
is  deformed to
%%%%%
\begin{eqnarray*}%~\label{eq:trans}
\tau^{\psi_d'}(M)
& = &
\xi^{-m}
\left\{\beta \left| \frac{\xi^m-1}{\xi-1}\right|^2 - \alpha' \right\}
(\xi^m-1)^2 \
= \zeta 
\left\{ \beta \zeta
\left( \frac{\zeta-1}{\zeta^2-1}\right)^2-\alpha' \right\} 
(\zeta-1)^{-2} \\
& = &
-\zeta^2
\left\{ \alpha' \zeta^{-1}
\left( \frac{\zeta^2-1}{\zeta-1}\right)^2-\beta \right\} 
(\zeta^2-1)^{-2} \
= - \zeta^2
\left\{ \alpha' (\zeta + \zeta^{-1}) -\beta \right\} 
(\zeta^2-1)^{-2}.
\end{eqnarray*}
%%%%%
We apply the same argument on this equality as in Subsection~\ref{ssec:key}
and retake $(i,j)$ satisfying the condition (\ref{eq:ij-cond}).
The equality of the real value sequence (\ref{eq:Rijm2}) is
%%%%%
\begin{equation*}~\label{eq:ABC}
\zeta^{-\frac{i+j}2}
\left\{ \alpha' \zeta^{-1}
\left( \frac{\zeta^2-1}{\zeta-1}\right)^2-\beta \right\} 
(\zeta^i-1)(\zeta^j-1)
\ = \ 
\eta \zeta^{-2} (\zeta^2-1)^2.
\end{equation*}
%%%%%
Divided by $\zeta^{-1}(\zeta -1)^2$, 
it induces similar equation to (\ref{eq:m=2tor}).
We can apply the same argument as in Subsection~\ref{ssec:m=2}.
Instead of (\ref{eq:m=2tor}), we study
%%%%%
\begin{equation}~\label{eq:m+1tor}
\zeta^{-\frac{i+j-2}{2}}\cdot
\left\{ \alpha' (\zeta + \zeta^{-1}) -\beta' \right\} 
\frac{(\zeta^i-1)(\zeta^j-1)}{(\zeta-1)^2}
=\eta \zeta^{-1}\cdot \frac{(\zeta^2-1)^2}{(\zeta-1)^2},
\end{equation}
%%%%%
where $\beta' = \beta - 2 \alpha'$($= \beta - 2(\alpha -mn\beta)$).
Equation (\ref{eq:m+1tor}) is obtained from (\ref{eq:m=2tor})
by changing ($\xi$ to $\zeta$), $\beta$ to $\alpha'$
and $\alpha''$ to $\beta$, respectively.
Thus, using the correspondence,
we can study their roots by the same argument.
In the last subsection, $\beta=0$ was not allowed, but 
here $\alpha' = 0$ is allowed (Lemma~\ref{lem:normR}).

\medskip
\noindent
(1) \ The case $\alpha'=0$\ ($\alpha/\beta=m(m+1)=mn$).

In this case, $\alpha/\beta=m(m+1)$ is a root
by Lemma~\ref{lem:special} (1).

\medskip
\noindent
(2) \ The case $\alpha' \ge 1$.
The argument is similar to the case $m=2$
in Subsection~\ref{ssec:m=2}.
\begin{enumerate}
\item[(i)]
The case $2\le i+j\le 2m-2$.

There is no root for $\alpha', \beta$,
by the argument with (1) in Subsection~\ref{ssec:m=2},
since the corresponding root $(\alpha', \beta') = (1,-2)$ 
implies $\beta = 0$, which is not allowed.

\medskip
\item[(ii)]
The case $i+j=2m$.

Then $e = (j-i)/2$ is an integer satisfying $1\le e < m$.
The argument is similar to that of (2) in Subsection~\ref{ssec:m=2}.
We have
\[
(\alpha' - \beta') \langle t^{m} \rangle
+ \alpha' \langle t^{m-1} \rangle
- \alpha'  \langle t^{e+1} \rangle
+ \beta' \langle  t^e \rangle
- \alpha' \langle t^{e-1} \rangle
= \eta (\langle  t^2 \rangle-2) \hfill
\]
Here we used $\langle t^{m+1} \rangle 
\equiv \langle t^m \rangle$ (mod $t^{2m+1}-1$).
The span of the left-hand side is at most $2m = 2[(m+n)/2]$.
Corresponding to that
the equation (\ref{eq:m=2tor}) has a 
root $(\beta, \alpha) = (-1,1)$, 
this equation has a root $(\alpha', \beta') = (1, -1)$ only if $m=2$, 
which implies $(\alpha, \beta) = (7, 1)$.
This corresponds to the case $(A_{2,3}; 7, 0)$, and it is a lens space,
see Subsection~\ref{ssec:move}.
\end{enumerate}

%%%%%%%%%%%%%%%%%%%%%%%%%%%%%%%%%%%%%%%%%%
%%%                   SubSection 5-3
%%%%%%%%%%%%%%%%%%%%%%%%%%%%%%%%%%%%%%%%%%
\subsection{The case $m\ge 3$ and $n\ge m+2$}~\label{ssec:ngem+2}
Let $P(t)$ be a symmetric Laurent polynomial in Definition~\ref{defini:symLPoly} of the form:
\[
P (t) = 
a_0 + \sum_{i=1}^{\infty} \ a_i 
\left (
t^ i + t^{-i}
\right )
= a_0 + \sum_{i=1}^{\infty} \ 
a_i \langle t^ i \rangle.
\]
We call 
$a_0$ the {\it constant term},
$a_i \langle t^ i \rangle$ the {\it $i$-th term},
and $a_i$ the {\it $i$-th coefficient}. 
When $P(t)$ is considered in 
$\mathbb{Q}[t, t^{-1}]/(t^N-1)$, 
$P(t)$ (as above) is {\it reduced} if 
$a_i =0$ for all $i > [N/2]$. 
We denote the reduction of $P(t)$ by $\textrm{red}(P(t))$.

\medskip \noindent
%%%%%%%%%%%%%%%%%%%%%%%%%%%%%%%%%
{\bf Assumption} (of this subsection) \ 
Let $(m,n)$ be a fixed pair of positive coprime integers 
with $m \geq 3$ and $n \geq m+2$.
We assume that $\beta \geq 1, \alpha' = \alpha - mn\beta \geq 0$,
see Lemma~\ref{lem:normR}.
We set $F(t)$ and $G(t)$ as:
%%%%%
\begin{eqnarray*}
F(t) & = &
\{
\beta \langle t^{m} \rangle
- \alpha' \langle t^{1} \rangle
+2 (\alpha' -\beta)
\}
\{
\langle t^{f} \rangle - \langle t^{e} \rangle
\},
\\
%%%
G(t) & = &
\langle t^{m+1} \rangle
- 2 \langle t^{m} \rangle
+ \langle t^{m-1} \rangle
- 2 \langle t^1 \rangle
+ 4
\end{eqnarray*}
%%%%%
in (\ref{eq:FG}) in the last section.
From now on, we fix $N=m+n$. Note that $G(t)$ is already reduced.
Let 
%%%%%
\begin{eqnarray}~\label{eq:Equation}
\textrm{red}(F(t)) = \pm G(t),
%F(t) \equiv \pm G(t) \mod (t^{m+n}-1).
\end{eqnarray}
%%%%%
be an equation on $(e, f)$. Existence 
of an integral solution is a necessary condition
for $(A_{m,n};\alpha/\beta,0)$ to be a lens space.
We show that there exist no integral solution of 
(\ref{eq:Equation}) with 
$0 \le e < f \le (m+n-1)/2$ under $\alpha' >0$,
see Lemma~\ref{lem:special} (1) for the case $\alpha'=0$.
By Lemma~\ref{lem:gcd}, we assume that $\alpha' \not=\beta$.
\par 

\medskip 
 
Using $\langle t^{a} \rangle \cdot \langle t^{b} \rangle
=  \langle t^{a+b} \rangle + \langle t^{a-b} \rangle$, we have
%%%%%
\begin{eqnarray*}
F(t) & = &
\beta \langle t^{m+f} \rangle 
+ \beta \langle t^{m-f} \rangle
- \alpha' \langle t^{f+1} \rangle 
- \alpha' \langle t^{f-1} \rangle
+ 2 (\alpha' -\beta)\langle t^{f} \rangle
\\
&  &
- \beta \langle t^{m+e} \rangle 
- \beta \langle t^{m-e} \rangle
+ \alpha' \langle t^{e+1} \rangle 
+ \alpha' \langle t^{e-1} \rangle
- 2 (\alpha' -\beta)\langle t^{e} \rangle.
\end{eqnarray*}
%%%%%
This consists of non-zero ten terms, but may not be reduced.
On the other hand, $G(t)$ has five terms 
and is already reduced.
Thus our problem is 
^^ ^^ Which term in $F(t)$ is reduced to a term in $G(t)$?"
and ^^ ^^ Which terms in $F(t)$ are cancelled with other terms?"
On the terms
$\langle t^{m+f} \rangle, \langle t^{m-f} \rangle, 
\langle t^{m+e} \rangle$ and $\langle t^{m-e} \rangle$, 
we have
\[
\textrm{red} (\langle t^{m+x} \rangle )
= \begin{cases}
\langle t^{m+x} \rangle & \textrm{ if } x \le (n-m)/2 \\
\langle t^{n-x} \rangle & \textrm{ if } x > (n-m)/2 \\
\end{cases},
\quad
\textrm{red} (\langle t^{m-x} \rangle )
= \langle t^{\vert m-x \vert } \rangle,
\]
for $x = e$ or $f$ (thus $0 \le x \le (m+n-1)/2$). 
The term $\langle t^{f+1} \rangle$ is already reduced 
except only one case:
\par \medskip
%%%%%%%%%%%%%%%%%%%%%%%%%%%%%%%%%%%%
{\bf Case~F} \
If $(m+n)$ is odd and $f = (m+n-1)/2$, it holds that
$\textrm{red} (\langle t^{f+1} \rangle ) = \langle t^{f} \rangle$.
%%%%%%%%%%%%%%%%%%%%%%%%%%%%%%%%%%%%
\par
\begin{figure}[h]
\begin{center}
\includegraphics[scale=0.4]{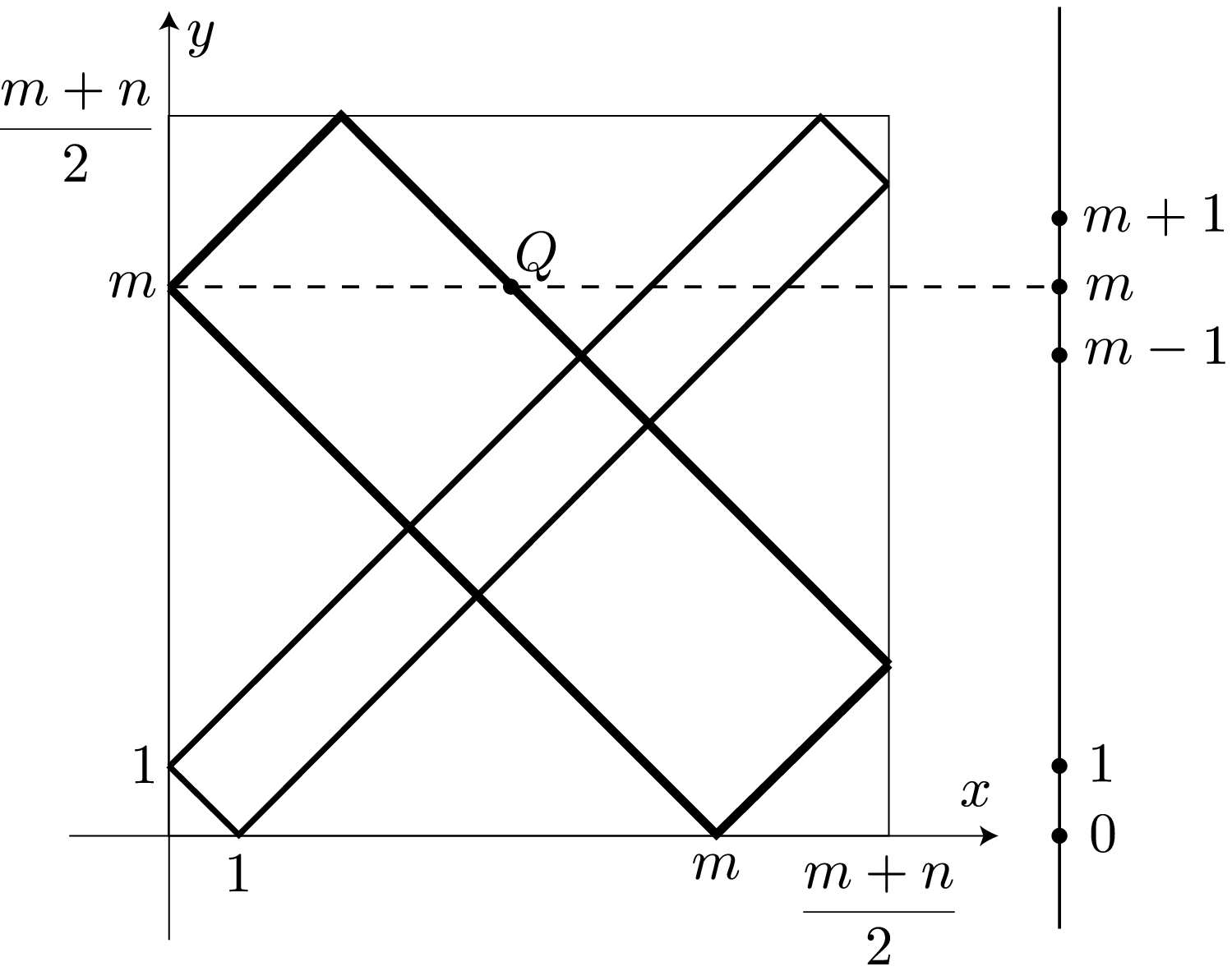}
\quad
\includegraphics[scale=0.4]{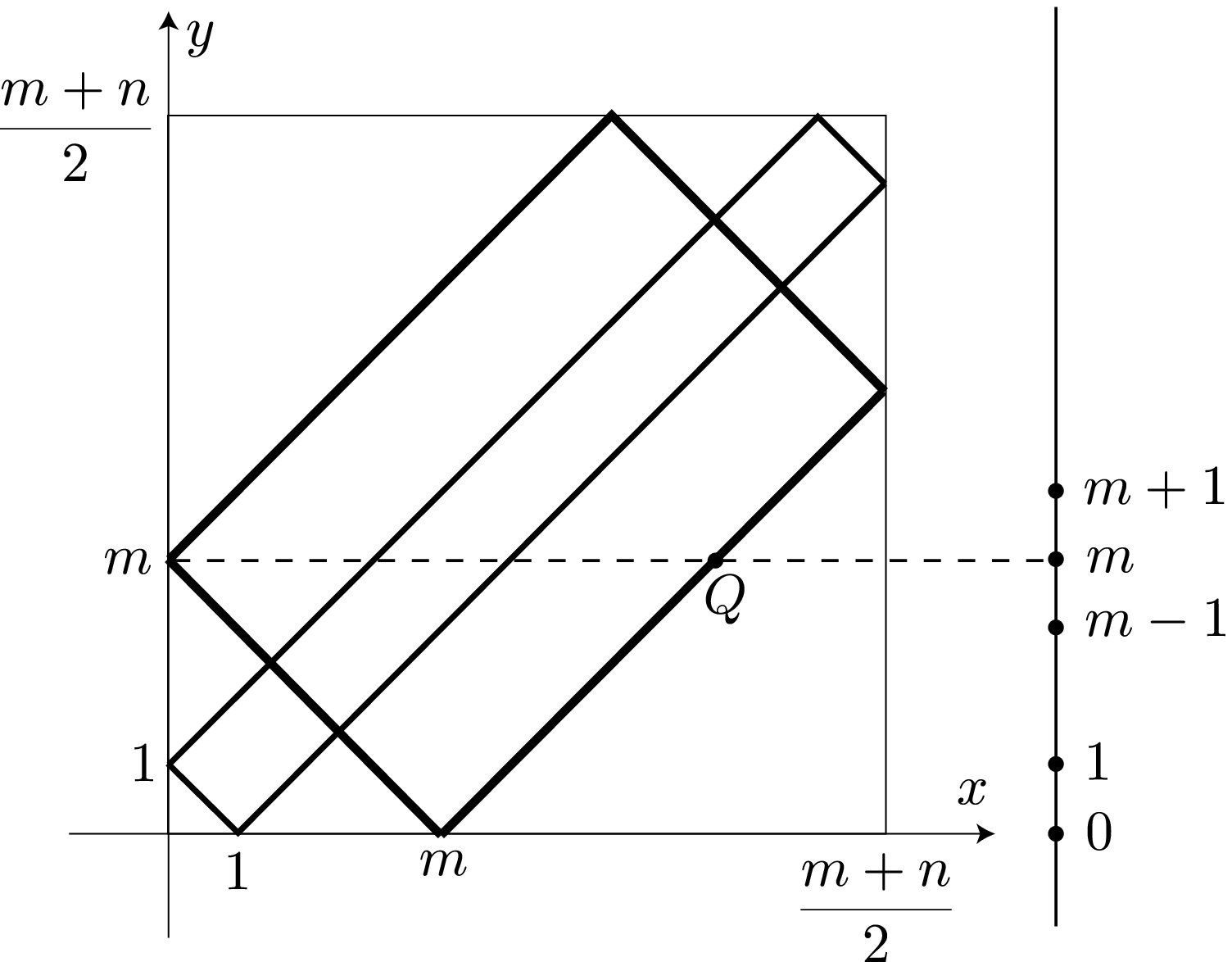} \\
(a) $n<3m$  \hskip5cm (b) $n>3m$ \qquad \qquad
\caption{Graph of the degrees I}
\label{fig:graph1}
\end{center}
\end{figure}

 \medskip \noindent
In Case~F, it holds that 
$f -m = (n-m-1)/2 \geq 1/2$, thus we have:
%
%
%%%%%%%%%%%%%%%%%%%%%%%%%%%%%%%%%%%%
\begin{lemma}~\label{lem:CaseF}
In Case~F, it holds that $f \geq m+1$. Furthermore, for an integer $a \geq 1$,
$f = m+ a$ is equivalent to $n = m+2a+1$.
\end{lemma}
%%%%%%%%%%%%%%%%%%%%%%%%%%%%%%%%%%%%
%
%
%
We will often take care of this exceptional case.
The other five terms are already reduced. It is easy to see:
%
%
%%%%%%%%%%%%%%%%%%%%%%%%%%%%%%%%%%%%
\begin{lemma}~\label{lem:Obs1}
\begin{enumerate}
\item[(1)] Neither $\langle t^{m+x} \rangle$ nor $\langle t^{x+1} \rangle$ 
can be reduced to the constant term.
\item[(2)] The term $\langle t^{2m} \rangle$ ($\langle t^{m+x} \rangle$ with $x=m$) can be reduced to 
neither the constant term nor the $1$-st term.
\end{enumerate}
\end{lemma}
%%%%%%%%%%%%%%%%%%%%%%%%%%%%%%%%%%%%
%
%
The graph of the degrees
of $\textrm{red}(\langle t^{m+x} \rangle)$,
$\textrm{red}(\langle t^{m-x} \rangle)$,
$\textrm{red}(\langle t^{x+1} \rangle)$ and
$\textrm{red}(\langle t^{x-1} \rangle)$
are useful, see Figure~\ref{fig:graph1}.
We are interested in only the points whose coordinates are integers.
We set the rectangle formed by
$y=\deg(\mathrm{red}(\langle t^{m+x}\rangle))$ 
and
$y=\deg(\mathrm{red}(\langle t^{m-x}\rangle))$
as $R$.
%
%
%%%%%%%%%%%%%%%%%
%%%%%%%%%%%%%%%%%
%
%
The intersection points of $R$ and the line $y = k$
correspond to solve the equation
$\mathrm{red}(\langle t^{m\pm x}\rangle)
=\langle t^k\rangle$.
The point $Q$ is one of the intersections for $k=m$
other than $(0, m)$.
% If $n>3m$, then $Q(2m, m)$. If $n<3m$, then $Q(n-m, m)$.
%In the graph, the point $Q,,$ on the line $y=m$ corresponds to the case 
%$\textrm{red}(\langle t^{m \pm x} \rangle) = \langle t^{m} \rangle$.
The $(x,y)$-coordinate of $Q$ is as follows: 
\par

If $n < 3m$, then $Q(n-m,m)$, which corresponds to 
$\textrm{red}(\langle t^{m + x} \rangle) = \langle t^{m} \rangle$, 
\par

If $n > 3m$, then $Q(2m,m)$, which corresponds to 
$\textrm{red}(\langle t^{m - x} \rangle) = \langle t^{m} \rangle$.
\par 

\noindent
Here, note that $n=3m$ contradicts coprimeness of $m$ and $n$.
%
%
%%%%%%%%%%%%%%%%%%%%%%%%%%%%%%%%%%%%
\begin{lemma}~\label{lem:Obs2}
{\rm (on the case 
$\textrm{red}(\langle t^{m \pm x} \rangle) = \langle t^m \rangle$) \\
}
By $x_Q$, we denote the $x$-coordinate ($= n-m$ or $2m$) of $Q$. Then we have
\begin{enumerate}
\item[(A)] $x_Q \not= m$, 
\item[(B)] $x_Q =m + a$ only if $n=2m+a$, for $a = -2, -1, 1$ or $2$,
\item[(C)] $x_Q =a$ only if $n=m+a$, for $a = 2$ or $3$.
\end{enumerate}
\end{lemma}
%%%%%%%%%%%%%%%%%%%%%%%%%%%%%%%%%%%%
%
%

%From now on, \underline{
%under the assumption 
%$\alpha' > 0$ and $\alpha' \not= \beta$}
%(i.e., every coefficient of $F(t)$ is non-zero),
%we show that any solution $(e,f)$ of 
%$\textrm{red}(F(t)) = \pm G(t)$ implies a contradiction.
Since $G(t)$ has the non-zero constant term,  
by Lemma~\ref{lem:Obs1} (1), we have
%at least one of $e$ and $f$ is equal to $0,1$ or $m$
\[
\{ e, f \} \cap \{ 0,1, m \} \not= \emptyset. 
\]
The proof is divided into the five cases:
(1) $e=0$, 
(2) $e=1$ with $m \geq 4$, 
(3) $e=m$, 
(4) ^^ ^^ $e \not=0,1$ and $f=m$", and 
(5) $e=1$ with $m=3$. 
Note that $f=0$ is impossible
and that the case $f=1$ is included by the case $e=0$,
because of the assumption $e<f$.
\par

\bigskip \noindent
%%%%%%%%%%%%        
%%%%%%%%%%%%
{\bf Case~1} ($e=0$) \ 
\begin{eqnarray*}
F(t) & = & 
\beta \langle t^{m+f} \rangle 
+ \beta \langle t^{m-f} \rangle
- \alpha' \langle t^{f+1} \rangle 
- \alpha' \langle t^{f-1} \rangle
+2 (\alpha' -\beta)\langle t^{f} \rangle 
\\
& &
- 2\beta \langle t^{m} \rangle
+ 2\alpha' \langle t^{1} \rangle
- 4 (\alpha' -\beta).
\end{eqnarray*}
Since $G(t)$ has a non-trivial $(m+1)$- and $(m-1)$-term, by Lemma~\ref{lem:Obs2},
we have (caring Case~F)
\[
f \in \{1, x_Q + \epsilon, m, m+1, m+2\} \cap \{1, x_Q - \epsilon, m-2, m-1, m \},
\]
where $\epsilon = -1$ if $n <3m$ (and $\epsilon = +1$ if $n >3m$, respectively).
%
%
%%%%%%%%%%%%%%%%%
\begin{figure}[h]
\begin{center}
\includegraphics[scale=0.4]{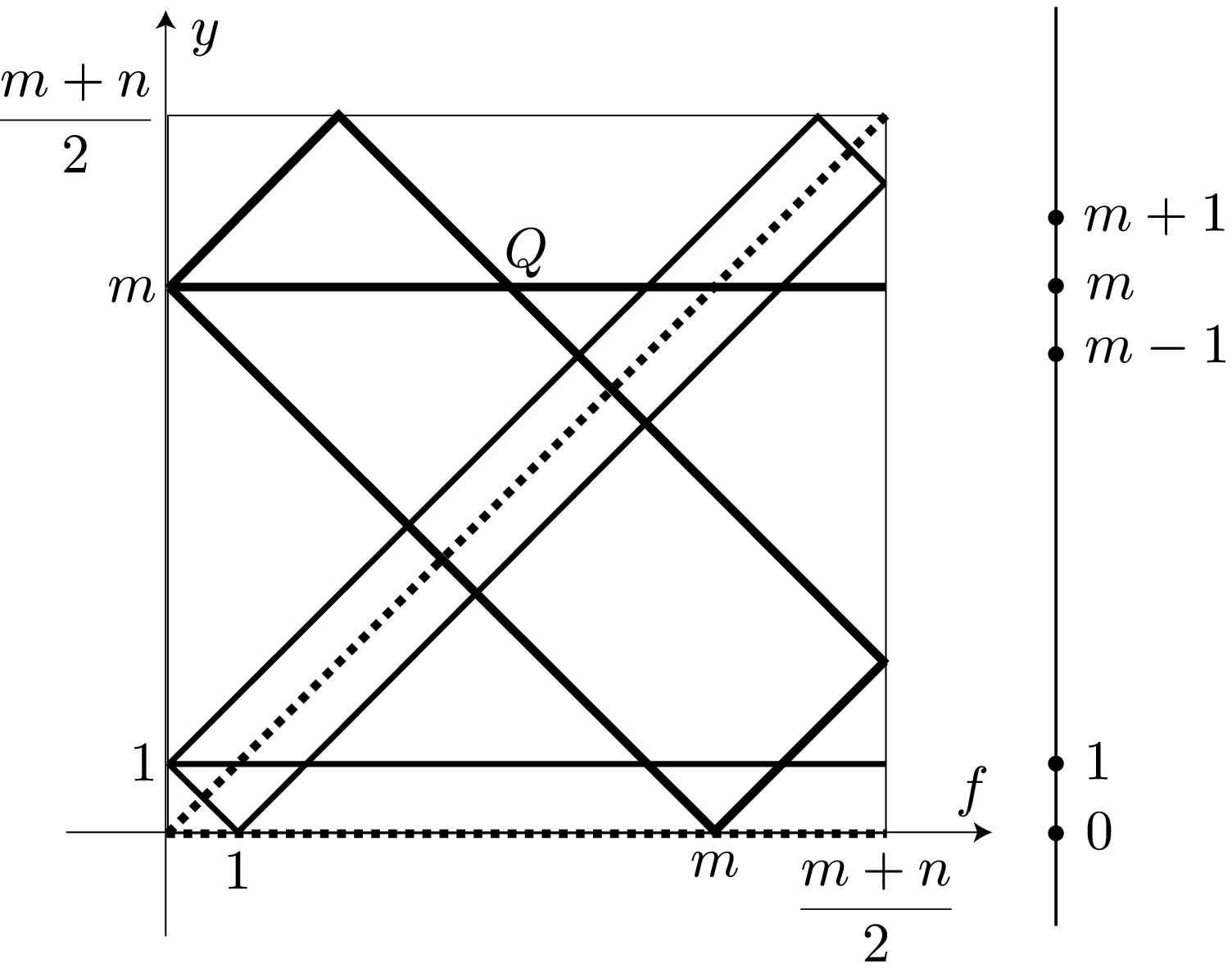}
\quad
\includegraphics[scale=0.4]{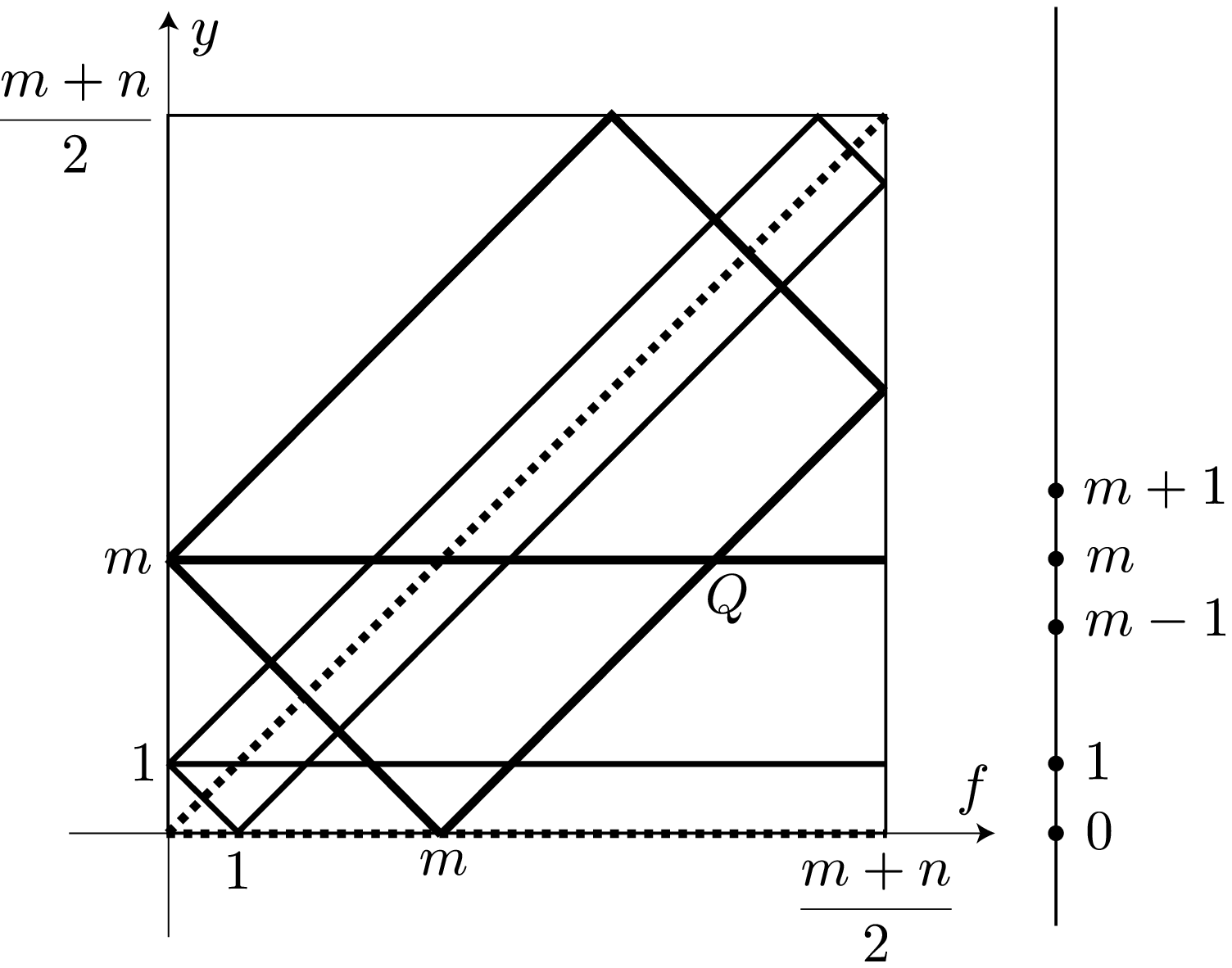} \\
(a) $n<3m$  \hskip5cm (b) $n>3m$ \qquad \qquad
\caption{Graph of the degrees II}
\label{fig:graphe2}
\end{center}
\end{figure} 
%%%%%%%%%%%%%%%%%
%
%
Figure~\ref{fig:graphe2} helps to understand it.
Thus we have four cases:
(i) $f = 1$,
(ii) $f = m$,
(iii) $f = m+2$ (if $n = 2m+1$ ($x_Q =m+1$)) or 
(iv) $f = m-2$ (if $n = 2m-1$ ($x_Q =m-1$))
by Lemma~\ref{lem:Obs2} (B).

In each case, Case~F 
(i.e., $\textrm{red}(\langle t^{f+1} \rangle) = \langle t^{f} \rangle$)
does not occur because $f \le m$
in (i), (ii) and (iv), and because $f \not= (m+n-1)/2$ in (iii),
by Lemma~\ref{lem:CaseF}.

\medskip \noindent %%%%%
Subcase(1-i): $(e,f)=(0,1)$. \ 
\[
F(t) = 
\beta \langle t^{m+1} \rangle 
- 2\beta \langle t^{m} \rangle
+ \beta \langle t^{m-1} \rangle
- \alpha' \langle t^{2} \rangle 
+2 (2\alpha' -\beta)\langle t^{1} \rangle 
- 2 (3\alpha' -2\beta).
\]
Since $\alpha' \not= 0$,
to cancel $-\alpha' \langle t^{2} \rangle$,
we need $m=3$ and $\alpha' = \beta$,
which contradicts the assumption.
(If we admit $\alpha' = 0$, then $(e,f)=(0,1)$ is a solution,
see Lemma~\ref{lem:special} (1).)

\medskip \noindent %%%%%
Subcase(1-ii): $(e,f)=(0,m)$. \ 
\[
F(t) = 
\beta \langle t^{2m} \rangle 
- \alpha' \langle t^{m+1} \rangle 
+2 (\alpha' -2\beta)\langle t^{m} \rangle 
- \alpha' \langle t^{m-1} \rangle
+ 2\alpha' \langle t^{1} \rangle
- 2 (2\alpha' - 3\beta).
\]
Since $2m>m+1$, the $2m$-th term is not reduced.
By Lemma~\ref{lem:Obs1} (2), and the ratio of the $(m+1)$-th and 
the $(m-1)$-th coefficients ($=1$), we have 
$\mathrm{red}(\langle t^{2m}\rangle)=\langle t^m\rangle$ and $n-m=m$.
It contradicts coprimeness of $m$ and $n$.
%If the first term $\beta \langle t^{2m} \rangle$ 
%is already reduced, we have a contradiction.
%(In the rest of the proof, we will omit this kind of sentence.)
%Considering the ratio of the coefficients of $G(t)$,
%we need (at least) 
%$\textrm{red}( \langle t^{2m} \rangle ) = \langle t^{m} \rangle $.
%But, it holds only if $n -m =m$, which contradicts to the coprimeness of $(m,n)$.

\medskip \noindent %%%%%
Subcase(1-iii): $(e,f)=(0,m+2)$ with $n=2m+1$. \ 
It holds that 
$\textrm{red}( \langle t^{2m+2} \rangle ) = \langle t^{m-1} \rangle$.
%%%%%
\begin{eqnarray*}
F(t) & = & 
- \alpha' \langle t^{m+3} \rangle 
+2 (\alpha' -\beta)\langle t^{m+2} \rangle 
\\
& &
- \alpha' \langle t^{m+1} \rangle
- 2\beta \langle t^{m} \rangle
+ \beta \langle t^{m-1} \rangle 
+ \beta \langle t^{2} \rangle
+ 2\alpha' \langle t^{1} \rangle
- 4 (\alpha' -\beta).
\end{eqnarray*}
%%%%%
Since $m\ge 3$ (thus $m+2 \le (m+n)/2$) and $\alpha' \not= \beta$,
the $(m+2)$-th term cannot be canceled. We have a contradiction.
%One of the first two terms has to be reduced and to cancel $\beta \langle t^{2} \rangle$,
%but it is impossible, because $\textrm{red}( \langle t^{m+a} \rangle )
%= \langle t^{2m+1 -a} \rangle $
%and $2m + 1 - a \geq 4$, for $a=2,3$.

\medskip \noindent %%%%%
Subcase(1-iv): $(e,f)=(0,m-2)$ with $n=2m-1$. \ 
If $m=3$, we go back to (1-i), thus we assume $m\geq 4$. 
It holds that $\textrm{red}( \langle t^{2m-2} \rangle ) = \langle t^{m+1} \rangle $.
%%%%%
\begin{eqnarray*}
F(t) & = & 
\beta \langle t^{m+1} \rangle 
- 2\beta \langle t^{m} \rangle
- \alpha' \langle t^{m-1} \rangle \\
& &
+2 (\alpha' -\beta)\langle t^{m-2} \rangle 
- \alpha' \langle t^{m-3} \rangle
+ \beta \langle t^{2} \rangle
+ 2\alpha' \langle t^{1} \rangle
- 4 (\alpha' -\beta).
\end{eqnarray*}
%%%%%
All terms are already reduced. 
By the signs of the $(m+1)$-th and the $(m-1)$-th coefficients,
we have a contradiction.
%Since the $(m-2)$-th coefficient $2(\alpha' -\beta) \not=0$,
%we need $m=4$ and $2 (\alpha' -\beta) + \beta =0$, i.e., 
%$\beta = 2\alpha'$.
%Then we have 
%\[
%\textrm{red}(F(t)) = \alpha' \cdot 
%(2\langle t^{5} \rangle 
%- 4\langle t^{4} \rangle
%- \langle t^{3} \rangle 
%+ \langle t^{1} \rangle
%+ 2),
%%2\alpha'\langle t^{5} \rangle 
%%- 4\alpha' \langle t^{4} \rangle
%%- \alpha' \langle t^{3} \rangle 
%%+ \alpha' \langle t^{1} \rangle
%%+ 2 \alpha',
%\]
%which is not equal to $G(t)$. We have a contradiction.
%\par
%
% We have proved that in the case $e=0$, there exists no solution.
%

\bigskip \noindent
%%%%%%%%%%%%
%%%%%%%%%%%%
{\bf Case~2} ($e=1$ with $m\geq 4$) \ 
%%%%%
\begin{eqnarray*}
F(t) & = & 
\beta \langle t^{m+f} \rangle 
+ \beta \langle t^{m-f} \rangle
- \alpha' \langle t^{f+1} \rangle 
- \alpha' \langle t^{f-1} \rangle
+ 2 (\alpha' -\beta)\langle t^{f} \rangle 
\\
& &
- \beta \langle t^{m+1} \rangle 
- \beta \langle t^{m-1} \rangle
+  \alpha' \langle t^{2} \rangle
- 2 (\alpha' -\beta) \langle t^{1} \rangle
+ 2\alpha'.
\end{eqnarray*}
%%%%%
%
%
%%%%%%%%%%%%%%%%%
\begin{figure}[h]
\begin{center}
\includegraphics[scale=0.4]{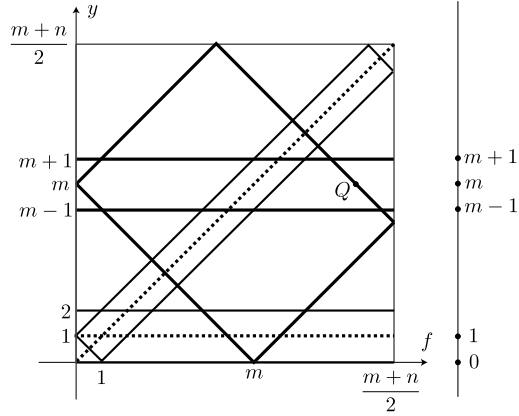}
\caption{Graph of the degrees III}
\label{fig:graph3}
\end{center}
\end{figure} 
%%%%%%%%%%%%%%%%%
%
%
Recall the assumption $f > e=1$.
Comparing $F(t)$ and $G(t)$, we have that 
at least one term in $F(t)$ is equal to or reduced to the $m$-th term, 
and that 
the term $\alpha' \langle t^{2} \rangle$ is canceled with another term
whose coefficient is negative.
Here we used $m \geq 4$
(i.e., $\langle t^{m-1} \rangle \not= \langle t^{2} \rangle$).
We have (caring Case~F)
\[
f \in \{ m-1, m, m+1, x_Q \} \cap \{ 2, 3 \},
\]
where $x_Q$ is the value defined in Lemma~\ref{lem:Obs2}. 
Figure~\ref{fig:graph3} helps to understand it.
Since $f \le 3$, Case~F does not occur by Lemma~\ref{lem:CaseF}.
Thus we have three cases:
(i) $m = 4$ and $f = 3$,
(ii) $f = 2$ (if $n = m+2$ ($x_Q =2$)),
(iii) $f = 3$ (if $n = m+3$ ($x_Q =3$)),
see Lemma~\ref{lem:Obs2} (C).

\medskip \noindent %%%%%
Subcase(2-i): $m=4$ and $(e,f)=(1,3)$.
\[
F(t) = 
\beta \langle t^{7} \rangle 
- \beta \langle t^{5} \rangle 
- \alpha' \langle t^{4} \rangle 
+ (2\alpha' -3\beta)\langle t^{3} \rangle 
- (2\alpha' -3\beta) \langle t^{1} \rangle
+ 2\alpha'.
\]
The term $\beta \langle t^{7} \rangle$ has to be reduced:
$\textrm{red} (\langle t^{7} \rangle) = \langle t^{n-3} \rangle$
and $3 \le n-3 \le 5$. 
By coprimeness of $m$ and $n$, we have $n=7$.
By the $5$-th and the $4$-th coefficients, 
we have a contradiction.
%
%Considering the ratio of the $1$-st coefficient and the constant term
%$- (2\alpha' - 3\beta) : 2 \alpha' = -2 : 4$,
%we have $\alpha' = 3\beta$
%and 
%\[
%F(t) = 
%\beta \cdot 
%(\langle t^{7} \rangle 
%- \langle t^{5} \rangle 
%- 3\langle t^{4} \rangle 
%+ 3\langle t^{3} \rangle 
%-3\langle t^{1} \rangle
%+ 6). 
%%\beta \langle t^{7} \rangle 
%%- \beta \langle t^{5} \rangle 
%%- 3\beta \langle t^{4} \rangle 
%%+ 3\beta \langle t^{3} \rangle 
%%-3\beta \langle t^{1} \rangle
%%+ 6\beta. 
%\]
%In any cases 
%$\textrm{red} (\langle t^{7} \rangle) = 
%\langle t^{5} \rangle$,
%$\langle t^{4} \rangle$ or $\langle t^{3} \rangle$, 
%it does not hold that $\textrm{red} (F(t)) = \pm G(t)$.
%We have a contradiction.

\medskip \noindent %%%%%
Subcase(2-ii): $m\geq 4, f = 2$ (if $n = m+2$ ($x_Q =2$)). \ 
It holds that $\textrm{red}( \langle t^{m+2} \rangle ) = \langle t^{m} \rangle $.
\[
F(t)  = 
- \beta \langle t^{m+1} \rangle 
+\beta \langle t^{m} \rangle 
- \beta \langle t^{m-1} \rangle
+ \beta \langle t^{m-2} \rangle
- \alpha' \langle t^{3} \rangle 
+ (3\alpha' -2\beta)\langle t^{2} \rangle 
- (3\alpha' -2\beta) \langle t^{1} \rangle
+ 2\alpha'.
\]
All terms are already reduced. 
Considering the ratio of the $(m+1)$-st and the $m$-th
coefficients, we have a contradiction.

\medskip \noindent %%%%%
Subcase(2-iii): $m\geq 4, f = 3$ (if $n = m+3$ ($x_Q =3$)). \ 
It holds that $\textrm{red}( \langle t^{m+3} \rangle ) = \langle t^{m} \rangle $.
\[
F(t)  = 
- \beta \langle t^{m+1} \rangle 
+ \beta \langle t^{m} \rangle 
- \beta \langle t^{m-1} \rangle
+ \beta \langle t^{m-3} \rangle
- \alpha' \langle t^{4} \rangle 
+ 2 (\alpha' -\beta)\langle t^{3} \rangle 
- 2 (\alpha' -\beta) \langle t^{1} \rangle
+ 2\alpha' 
\]
All terms are already reduced. 
If $m \geq 5$, then we have a contradiction
by the same method with the last case.
If $m = 4$, 
\[
\textrm{red}(F(t))  = 
- \beta \langle t^{5} \rangle 
- (\alpha' - \beta) \langle t^{4} \rangle 
+ (2\alpha' -3\beta)\langle t^{3} \rangle 
- (2\alpha' -3\beta) \langle t^{1} \rangle
+ 2\alpha'.
\]
Considering the ratio of the first two coefficients,
$-\beta : -(\alpha' - \beta ) = 1 : -2$, i.e., $\beta = - \alpha' < 0$.
We have a contradiction. 
\par

\medskip

In the rest of the proof, we will often use the following:
%
%
%%%%%%%%%%%%%%%%%%%%%%%%%%%%%%%%%%%%
\begin{lemma}~\label{lem:Obs3}
Even if 
$\textrm{red} (\langle t^{a} \rangle ) = \langle t^{b} \rangle$, 
the sum $\pm \beta \langle t^{a}\rangle \pm  \alpha' \langle t^{b} \rangle$
of any signs never cancel by the reduction.
\end{lemma}
%%%%%%%%%%%%%%%%%%%%%%%%%%%%%%%%%%%%
%
%
This lemma is easily shown by the assumption
$\alpha' > 0,  \beta > 0$ and $\alpha' \not= \beta$. 

\par \bigskip \noindent
%%%%%%%%%%%%
%%%%%%%%%%%%
{\bf Case~3} ($e=m$) \ 
%%%%%
\begin{eqnarray*}
F(t) & = &
\beta \langle t^{m+f} \rangle 
+ \beta \langle t^{m-f} \rangle
- \alpha' \langle t^{f+1} \rangle 
- \alpha' \langle t^{f-1} \rangle
+2 (\alpha' -\beta)\langle t^{f} \rangle
\\
&  &
- \beta \langle t^{2m} \rangle 
- 2\beta 
+ \alpha' \langle t^{m+1} \rangle 
+ \alpha' \langle t^{m-1} \rangle
-2 (\alpha' -\beta)\langle t^{m} \rangle.
\end{eqnarray*}
%%%%%
Recall the assumption $f > e=m$.
Comparing $F(t)$ and $G(t)$, 
at least one term in $F(t)$ is equal to or reduced to the $1$-st term.
%Since $\langle t^{2m} \rangle$ is not 
%reduced to the $1$-st term (by Lemma~\ref{lem:Obs1} (2)), 
By Lemma~\ref{lem:Obs1} (2),
we have two cases:
(i) $f = m+1$ (It can be in Case~F) and (ii)
$\textrm{red}(\langle t^{m+f} \rangle )= \langle t^{1} \rangle$.
But in the latter case, $f = n-1 \le (m+n)/2$
(see the graph in Figure~\ref{fig:graph1}),
thus $n =m+2$ and $f = m+1$. 
We only have to study the case (i).

\medskip \noindent %%%%%
Subcase(3-i): $(e,f) = (m,m+1)$. \
%%%%%
\begin{eqnarray*}
F(t) & = &
\beta \langle t^{2m+1} \rangle 
- \beta \langle t^{2m} \rangle 
- \alpha' \langle t^{m+2} \rangle 
\\
&  &
+ (3\alpha' -2\beta)\langle t^{m+1} \rangle
- (3\alpha' - 2\beta)\langle t^{m} \rangle
+ \alpha' \langle t^{m-1} \rangle
+ \beta \langle t^{1} \rangle
- 2\beta.
\end{eqnarray*}
%%%%%
We focus on $- \alpha' \langle t^{m+2} \rangle$. 
Neither 
$\beta \langle t^{2m+1} \rangle- \alpha' \langle t^{m+2} \rangle$
nor 
$-\beta \langle t^{2m} \rangle- \alpha' \langle t^{m+2} \rangle$
cancel by the reduction by Lemma~\ref{lem:Obs3}.
Thus, if $m+2 \le (m+n)/2$ (i.e., $n \geq m+4$), 
the term $- \alpha' \langle t^{m+2} \rangle$ is left after the reduction,
which is a contradiction. 
We have $n = m+2$ or $n=m+3$.
\par

First we assume $n = m+2$,
then 
$\textrm{red} (
\beta \langle t^{2m+1} \rangle 
- \beta \langle t^{2m} \rangle 
- \alpha' \langle t^{m+2} \rangle 
) = 
\beta \langle t^{1} \rangle 
- \beta \langle t^{2} \rangle 
- \alpha' \langle t^{m} \rangle$.
Considering the ratio of the $(m+1)$-th and the $(m-1)$-th coefficient, 
we need $\alpha' = \beta$. We have a contradiction.
Next we assume $n = m+3$, which implies $m \geq 4$.
Then 
$\textrm{red} (
\beta \langle t^{2m+1} \rangle 
- \beta \langle t^{2m} \rangle 
- \alpha' \langle t^{m+2} \rangle 
) = 
\beta \langle t^{2} \rangle 
- \beta \langle t^{3} \rangle 
- \alpha' \langle t^{m+1} \rangle$. (It is in Case~F.)
The second coefficient $\beta \not= 0$, we have a contradiction.
\par
%
%We have proved that in the case $e=m$, there exists no solution.
%

\bigskip \noindent
%%%%%%%%%%%%
%%%%%%%%%%%%
{\bf Case~4} ($f=m$) \ It is not in Case~F by Lemma~\ref{lem:CaseF}.
%%%%%
\begin{eqnarray*}
F(t) & = &
- \beta \langle t^{m+e} \rangle 
- \beta \langle t^{m-e} \rangle
+ \alpha' \langle t^{e+1} \rangle 
+ \alpha' \langle t^{e-1} \rangle
-2 (\alpha' -\beta)\langle t^{e} \rangle
\\
&  &
+ \beta \langle t^{2m} \rangle 
+ 2\beta 
- \alpha' \langle t^{m+1} \rangle 
- \alpha' \langle t^{m-1} \rangle
+ 2 (\alpha' -\beta)\langle t^{m} \rangle.
\end{eqnarray*}
%%%%%
Recall the assumption $e < f =m$.
Comparing $F(t)$ and $G(t)$, we have that
at least one term in $F(t)$ is equal to or reduced to the $1$-st term.
By Lemma~\ref{lem:Obs1} (2), $\langle t^{2m} \rangle$ is not reduced to the $1$-st term.
It does not hold that 
$\textrm{red}(\langle t^{m+e} \rangle )= \langle t^{1} \rangle$,
because it implies $e = n-1$, which contradicts to $e < f = m$.
Thus we need $\textrm{red}(\langle t^{m-e} \rangle )= \langle t^{1} \rangle$, i.e., $e = m-1$.
%%%%%
\begin{eqnarray*}
F(t) & = &
\beta \langle t^{2m} \rangle 
- \beta \langle t^{2m-1} \rangle 
\\
&  &
- \alpha' \langle t^{m+1} \rangle 
+ (3\alpha' -2\beta)\langle t^{m} \rangle
- (3\alpha' -2\beta)\langle t^{m-1} \rangle
+ \alpha' \langle t^{m-2} \rangle
- \beta \langle t^{1} \rangle
+ 2\beta.
\end{eqnarray*}
%%%%%
One of the first two terms has to cancel 
$\alpha' \langle t^{m-2} \rangle$, but it is impossible,
because neither 
$\beta \langle t^{2m} \rangle + \alpha' \langle t^{m-2} \rangle$
nor 
$- \beta \langle t^{2m-1} \rangle + \alpha' \langle t^{m-2} \rangle$
cancel by Lemma~\ref{lem:Obs3}.
%
%We have proved that in the case $f=m$, there exists no solution.
%

\bigskip \noindent
%%%%%%%%%%%%        
%%%%%%%%%%%%
{\bf Case~5} ($e=1$ with $m=3$) \ 
Note that $G(t) = 
\langle t^{4} \rangle
- 2 \langle t^{3} \rangle
+ \langle t^{2} \rangle
- 2 \langle t^1 \rangle
+ 4$.
%%%%%
\begin{eqnarray*}
F(t) & = & 
\beta \langle t^{3+f} \rangle 
+ \beta \langle t^{3-f} \rangle
- \alpha' \langle t^{f+1} \rangle 
- \alpha' \langle t^{f-1} \rangle
+ 2 (\alpha' -\beta)\langle t^{f} \rangle 
\\
& &
- \beta \langle t^{4} \rangle 
+ (\alpha' -\beta) \langle t^{2} \rangle
- 2 (\alpha' -\beta) \langle t^{1} \rangle
+ 2\alpha'.
\end{eqnarray*}
%%%%%
The proof is divided into four cases:
(i) $f = 2$,
(ii) $f = 3$,
(iii) $f = 4$ or
(iv) $f \geq 5$. 
\par

\medskip \noindent %%%%%
Subcase(5-i):  $(e,f)=(1,2)$.
\[
F(t) = 
\beta \langle t^{5} \rangle 
- \beta \langle t^{4} \rangle 
- \alpha' \langle t^{3} \rangle 
+ 3 (\alpha' -\beta)\langle t^{2} \rangle 
- 3 (\alpha' -\beta) \langle t^{1} \rangle
+ 2\alpha'.
\]
The term $\beta \langle t^{5} \rangle$ has to be reduced:
$\textrm{red} (\langle t^{5} \rangle) = \langle t^{n-2} \rangle$
and $n-2 \geq m \geq 3$. 
Considering the ratio of the $1$-st coefficient and the constant term,
we have $- 3 (\alpha' -\beta) : 2 \alpha' = -2 : 4$, i.e.,
$2\alpha' = 3\beta$.
Considering the ratio of the $2$-nd and the $1$-st coefficients,
we have a contradiction.

\medskip \noindent %%%%%
Subcase(5-ii): $(e,f)=(1,3)$.
\[
F(t) = 
\beta \langle t^{6} \rangle 
- (\alpha' + \beta) \langle t^{4} \rangle 
+ 2 (\alpha' -\beta)\langle t^{3} \rangle 
-\beta \langle t^{2} \rangle
- 2 (\alpha' -\beta) \langle t^{1} \rangle
+ 2(\alpha' + \beta) 
\]
The term $\beta \langle t^{6} \rangle$ has to be reduced:
$\textrm{red} (\langle t^{6} \rangle) = \langle t^{n-3} \rangle$
and $n-3 \geq 2$. 
Considering the ratio of the $1$-st coefficient and the constant term,
we have $-2 (\alpha' -\beta)  : 2(\alpha' + \beta) = -2 : 4$, i.e.,
$\alpha' = 3\beta$.
\[
F(t) = \beta \cdot (
\langle t^{6} \rangle 
- 4\langle t^{4} \rangle 
+ 4\langle t^{3} \rangle 
-\langle t^{2} \rangle
- 4\langle t^{1} \rangle
+ 8).
%\beta \langle t^{6} \rangle 
%- 4\beta \langle t^{4} \rangle 
%+ 4\beta \langle t^{3} \rangle 
%-\beta \langle t^{2} \rangle
%- 4 \beta \langle t^{1} \rangle
%+ 8\beta.
\]
In any cases 
$\textrm{red} (\langle t^{6} \rangle) = 
\langle t^{4} \rangle$,
$\langle t^{3} \rangle$ or $\langle t^{2} \rangle$, 
it does not hold that $\textrm{red} (F(t)) = \pm G(t)$.
We have a contradiction.

\medskip \noindent %%%%%
Subcase(5-iii): $(e,f)=(1,4)$.
\[
F(t) =
\beta \langle t^{7} \rangle 
- \alpha' \langle t^{5} \rangle 
+ (2\alpha' -3\beta)\langle t^{4} \rangle 
- \alpha' \langle t^{3} \rangle
+ (\alpha' -\beta) \langle t^{2} \rangle
- (2\alpha' -3\beta) \langle t^{1} \rangle
+ 2\alpha'.
\]
By Lemma~\ref{lem:Obs3}, both terms 
$\beta \langle t^{7} \rangle$ and $- \alpha' \langle t^{5} \rangle$
have to be reduced.
It is possible only if $(m+n)/2 < 5$, i.e., $n < 7$.
Since $n \geq m+2 =5$, and $n$ is coprime to $m=3$, we have $n=5$.
Then 
$\textrm{red} (\langle t^{7} \rangle) = \langle t^{1} \rangle$ and
$\textrm{red} (\langle t^{5} \rangle) = \langle t^{3} \rangle$.
Considering the ratio of the $4$-th and the $2$-nd coefficients, 
we have 
$(2\alpha' -3\beta) = (\alpha' -\beta)$, i.e., $\alpha' = 2\beta$.
In this case, we lose the $1$-st term:
$\textrm{red} (
\beta \langle t^{7} \rangle
- (2\alpha' -3\beta) \langle t^{1}\rangle ) =0$. 
We have a contradiction.

\medskip \noindent %%%%%
Subcase(5-iv): $e=1$ and $f \geq 5$.
%%%%%
\begin{eqnarray*}
F(t) 
& = & 
\beta \langle t^{f+3} \rangle 
- \alpha' \langle t^{f+1} \rangle 
+ 2 (\alpha' -\beta)\langle t^{f} \rangle 
\\
& &
- \alpha' \langle t^{f-1} \rangle
+ \beta \langle t^{f-3} \rangle
- \beta \langle t^{4} \rangle 
+ (\alpha' -\beta) \langle t^{2} \rangle
- 2 (\alpha' -\beta) \langle t^{1} \rangle
+ 2\alpha'.
\end{eqnarray*}
%%%%%
The $f$-th term has to be cancelled. 
If it is cancelled by the $(f+3)$-th term,
then the $(f+1)$-th term is left after the reduction. 
If it is cancelled by the $(f+1)$-th term, which is in Case~F,
we have $\alpha' = 2 \beta$, $n= 2f-2$, 
and $\textrm{red}(\langle t^{f+3} \rangle ) 
= \langle t^{f-2} \rangle$ with $f-2 \geq 3$, thus
\[
\textrm{red}(F(t)) =  \beta \cdot (
- 2 \langle t^{f-1} \rangle
+ \langle t^{f-2} \rangle 
+ \langle t^{f-3} \rangle
- \langle t^{4} \rangle 
+ \langle t^{2} \rangle
- 2 \langle t^{1} \rangle
+ 4).
%- 2 \beta \langle t^{f-1} \rangle
%+ \beta \langle t^{f-2} \rangle 
%+ \beta \langle t^{f-3} \rangle
%- \beta \langle t^{4} \rangle 
%+ \beta \langle t^{2} \rangle
%- 2 \beta \langle t^{1} \rangle
%+ 4\beta.
\]
If $f > 5$, then the top degree is $f-1>4$, we have a contradiction.
Otherwise $f =5$, we have
\[
\textrm{red}(F(t)) =  \beta  \cdot
(- 3 \langle t^{4} \rangle
+ \langle t^{3} \rangle 
+ 2 \langle t^{2} \rangle
- 2 \langle t^{1} \rangle
+ 4).
%- 3 \beta \langle t^{4} \rangle
%+ \beta \langle t^{3} \rangle 
%+ 2\beta \langle t^{2} \rangle
%- 2 \beta \langle t^{1} \rangle
%+ 4\beta.
\]
The ratio of the coefficients is different from that of $G(t)$.
We have a contradiction.
\par
 
\bigskip
In any cases, we have a contradiction.
The proof of Theorem~\ref{thm:MT1} (1) is completed. 
\qed 
\par

%%%%%%%%%%%%%%%%%%%%%%%%%%%%%%%%%%%%%%%%%%
%%%                   SubSection 5-4
%%%%%%%%%%%%%%%%%%%%%%%%%%%%%%%%%%%%%%%%%%
%\subsection{Proof of $(A_{m,n}; mn, 0)\cong L((m+n)^2, m\overline{n})$}
\subsection{Type of the lens space $(A_{m,n}; mn, 0)$}
~\label{ssec:type}
We verify the second term ($= m\overline{n}$) of the lens space
$(A_{m,n}; mn, 0)\cong L((m+n)^2, m\overline{n})$
by the Reidemeister torsion.
Here we are also interested in possibility of the 
transformation between the parameters
$(m,n)$ and $(p,q)$. 
\par

\medskip

We use the notations ${\frak I}(m, n)=\{k_i\ |\ i=0, 1, \ldots, m+n-1\}$,
$u_j\ (j=0, 1, \ldots, n)$ and $w_j\ (j=0, 1, \ldots, m)$
defined in Section~\ref{sec:Alex}.
%
%
%%%%%%%%%%%%%%%%%%%%%%%%%%%%%%%%%%%%
\begin{lemma}~\label{lem:ki'}
\begin{enumerate}
\item[(1)]
We set $k'_i=k_i(m+n)-imn\ (i=0, 1, \ldots, m+n-1)$.
Then we have $k'_i+k'_{m+n-1-i}=mn$ and
$k'_0=0\le k'_i\le k'_{m+n-1}=mn$.

\item[(2)]
${\frak I}(m, n)=\{k'_i\ |\ i=0, 1, \ldots, m+n-1\}$.
\end{enumerate}
\end{lemma}
%%%%%%%%%%%%%%%%%%%%%%%%%%%%%%%%%%%%
%
%
\medskip \noindent
{\bf Proof} \ 
(1)\ By the definition of $k'_i$ and 
Proposition~\ref{prop:Imn} (1), we have
%%%%%
\begin{eqnarray*}
k'_i+k'_{m+n-1} & = & k_i(m+n)-imn+k_{m+n-1-i}(m+n)-(m+n-1-i)mn\\
& = & (k_i+k_{m+n-1-i})(m+n)-(m+n-1)mn\\
& = & mn(m+n)-(m+n-1)mn=mn.
\end{eqnarray*}
%%%%%
By Proposition~\ref{prop:Imn} (2),
every element $k_i\in {\frak I}(m, n)$ is 
uniquely expressed of the form $k_{u_j}$ or $k_{w_j}$.
It is easy to see that $k'_0=0$ and $k'_{m+n-1}=mn$.

\medskip \noindent
(i)\ The case $i=u_j\ (j=1, 2, \ldots, n-1)$. \ 
Since $k_i=k_{u_j}=jm$ by Proposition~\ref{prop:Imn} (2)
and $\gcd(m, n)=1$, we have
\[
k'_i=jm(m+n)-\left( \left[ \frac{jm}{n}\right]+j\right)mn
=mn\left( \frac{jm}{n}-\left[ \frac{jm}{n}\right]\right)>0.
\]
Since $k'_i+k'_{m+n-1}=mn$, we have $0<k'_i<mn$.

\medskip \noindent
(ii)\ The case $i=w_j\ (j=1, 2, \ldots, m-1)$. \ 
Since $k_i=k_{w_j}=jn$ by Proposition~\ref{prop:Imn} (2) 
and $\gcd(m, n)=1$, we have
\[
k'_i=jn(m+n)-\left( \left[ \frac{jn}{m}\right]+j\right)mn
=mn\left( \frac{jn}{m}-\left[ \frac{jn}{m}\right]\right)>0.
\]
Since $k'_i+k'_{m+n-1}=mn$, we have $0<k'_i<mn$.
Therefore we have the result.

\bigskip 
\noindent
(2)\ (i)\ The case $i=u_j\ (j=1, 2, \ldots, n-1)$. \ 
By the proof of (1) (i), 
$k'_i$ is of the form $k'_i=j'm$ for some $j'\in \mathbb{Z}$ and
$j'\equiv jm\ (\mathrm{mod}\ \! n)$ is uniquely determined.

\medskip \noindent
(ii)\ The case $i=w_j\ (j=1, 2, \ldots, m-1)$. \ 
By the proof of (1) (ii), 
$k'_i$ is of the form $k'_i=j'n$ for some $j'\in \mathbb{Z}$ and
$j'\equiv jn\ (\mathrm{mod}\ \! m)$ is uniquely determined.

\medskip

From (i), (ii) and (1),
the set $\{k'_i\ |\ i=0, 1, \ldots, m+n-1\}$ consists of
distinct $m+n$ elements, and we have the result.
\qed 
\par
%
%
%%%%%%%%%%%%%%%%%%%%%%%%%%%%%%%%%%%%
\begin{lemma}~\label{lem:Amn(T)}
For a coprime positive pair $(m, n)$,
The Alexander polynomial ${\it \Delta}_{A_{m,n}}(t, x)$ of 
$A_{m,n}$ satisfies that
\[
{\it \Delta}_{A_{m,n}}(T^{m+n}, T^{-mn})
\ \doteq \ 
\dfrac{(T^{mn}-1)(T^{m+n}-1)}{(T^m-1)(T^n-1)}.
\]
\end{lemma}
%%%%%%%%%%%%%%%%%%%%%%%%%%%%%%%%%%%%
%
%
\medskip \noindent
{\bf Proof} \ 
As we remarked in Section~\ref{sec:intro},
the first component $K_1$ of $A_{m,n}$ is
the $(m, n)$-torus knot $T_{m,n}$.
By Theorem~\ref{thm:AmnAlex} and 
the Torres formula (Lemma~\ref{lem:Torres}), 
we have
\[
{\it \Delta}_{A_{m,n}}(t, 1)
\doteq \frac{t^{m+n}-1}{t-1}\cdot {\it \Delta}_{T_{m,n}}(t)
\doteq \frac{(t^{mn}-1)(t^{m+n}-1)}{(t^m-1)(t^n-1)}
=\sum_{i=0}^{m+n-1} t^{k_i}.
\]
By Theorem~\ref{thm:AmnAlex}, it holds that
\[
{\it \Delta}_{A_{m,n}}(T^{m+n}, T^{-mn})
\ \doteq \ 
\sum_{i=0}^{m+n-1} T^{k_i(m+n)-imn}
\quad
= \sum_{i=0}^{m+n-1} T^{k'_i},
\]
and hence
we have 
${\it \Delta}_{A_{m,n}}(T^{m+n}, T^{-mn})
\doteq 
{\it \Delta}_{A_{m,n}}(T, 1)$
by Lemma~\ref{lem:ki'} (2).
\qed 
\par

\bigskip
\noindent
{\bf Proof} of $(A_{m,n}; mn, 0)\cong L((m+n)^2, m\overline{n})$ 
\par 
\noindent
Assuming that $M = (A_{m,n}; mn, 0)$ is a lens space $L(P,Q)$, 
it is clear that $P = (m+n)^2$ by the first homology.
By Lemma \ref{lem:Amn(T)} and the formulae (\ref{eq:torM1A}) and (\ref{eq:torM2A}), we have
%%%%%
\begin{equation*}
\tau (M_1)\doteq
{\it \Delta}_{A_{m,n}}(T^{m+n}, T^{-mn})(T^{m+n}-1)^{-1}
= \dfrac{T^{mn}-1}{(T^{m}-1)(T^{n}-1)},
\end{equation*}
%%%%%
\begin{equation*}
\tau^{\psi_{(m+n)^2}}(M)\doteq
(\zeta_{(m+n)^2}^m-1)^{-1}(\zeta_{(m+n)^2}^n-1)^{-1}.
\end{equation*}
%%%%%
Because of the ambiguity of the Reidemeister torsion, 
the parameters of the lens space is almost decided but
the orientation (at the choice $\pm$) is left undecided.
\[
M\cong L((m+n)^2, \pm m\overline{n})
%\ \textrm{or} \ L((m+n)^2, \pm \overline{m}n),
\]
where $\overline{n}$ is taken modulo $(m+n)^2$.

\medskip

Using the knowledge of existence of a lens space surgery
 $(A_{m,n}; mn, 0)\cong L(p^2, pq-1)$
in Theorem~\ref{thm:(mn)}, 
we can determine the sign of $\pm m\overline{n}$ as follows:
Suppose that it is $-m\overline{n}$.
Then we have 
$pq-1\equiv -m\overline{n}\ (\mathrm{mod}\ \! p^2)$, and 
%or $-\overline{m}n\ (\mathrm{mod}\ \! p^2)$. 
%If $pq-1\equiv -m\overline{n}$, 
\[
pq\equiv 1-m\overline{n}\equiv \overline{n}(n-m)
\quad (\mathrm{mod}\ \! p^2).
\]
Since $\gcd(pq, p^2)=p=m+n$ is changed to
$\gcd(\overline{n}(n-m), (m+n)^2) = m+n$
(we regard $\overline{n}$ as an integral lift),
which contradicts that $\gcd(\overline{n}, m+n)=1$ and
$\gcd(n-m, m+n) = \gcd(n-m, 2m)\le 2m < m+n$.
\qed

%%%%%%%%%%%%%%%%%%%%%%%%%%%%%%%%%%%%%%%%%%
%%%                   SubSection 5-5
%%%%%%%%%%%%%%%%%%%%%%%%%%%%%%%%%%%%%%%%%%
%\subsection{Proof of $(A_{2,3}; 7, 0)\cong L(25, 7)$}~\label{ssec:move}
\subsection{Type of the lens space $(A_{2,3}; 7, 0)$}~\label{ssec:move}
Suppose that $(m, n)=(2, 3)$ and $\alpha/\beta=7$.
Then we have
%%%%%
\begin{equation*}
{\it \Delta}_{A_{2,3}}(t, x)
=1+t^2x+t^3x^2+t^4x^3+t^6x^4
\end{equation*}
%%%%%
by Theorem~\ref{thm:AmnAlex}, and
(\ref{eq:torM1A}) is computed as
%%%%%
\begin{equation*}
\tau(M_1)\doteq
{\it \Delta}_{A_{2,3}}(T^5, T^{-7})(T^5-1)^{-1}
=\frac{1+T^3+T+T^{-1}+T^2}{T^5-1}\doteq \frac{1}{T-1}. 
\end{equation*}
%%%%%
Hence we have
%%%%%
\begin{equation*}
\tau^{\psi_{25}}(M)\doteq 
(\zeta_{25}-1)^{-1}(\zeta_{25}^7-1)^{-1}
\end{equation*}
%%%%%
by (\ref{eq:torM2A}).
Assuming that
$(A_{2,3}; 7, 0)$ is a lens space, 
only by our method of
Reidemeister--Turaev torsion,
we can say that the lens space is homeomorphic to 
$L(25, 7)\cong L(25, -7)$.

\begin{figure}[h]
\begin{center}
\includegraphics[scale=0.5]{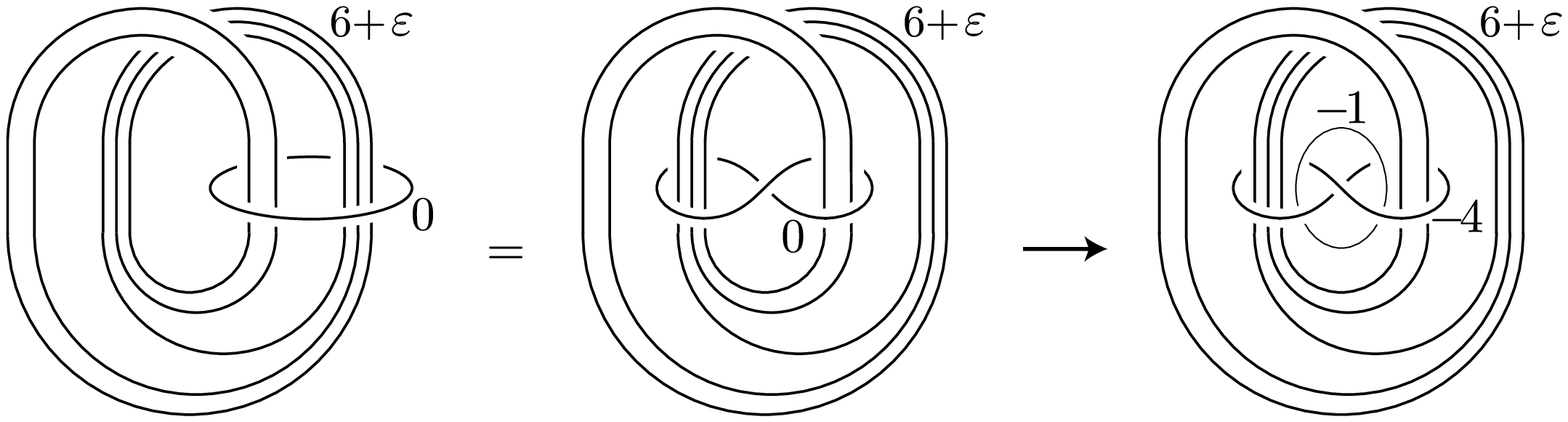}\\
\includegraphics[scale=0.5]{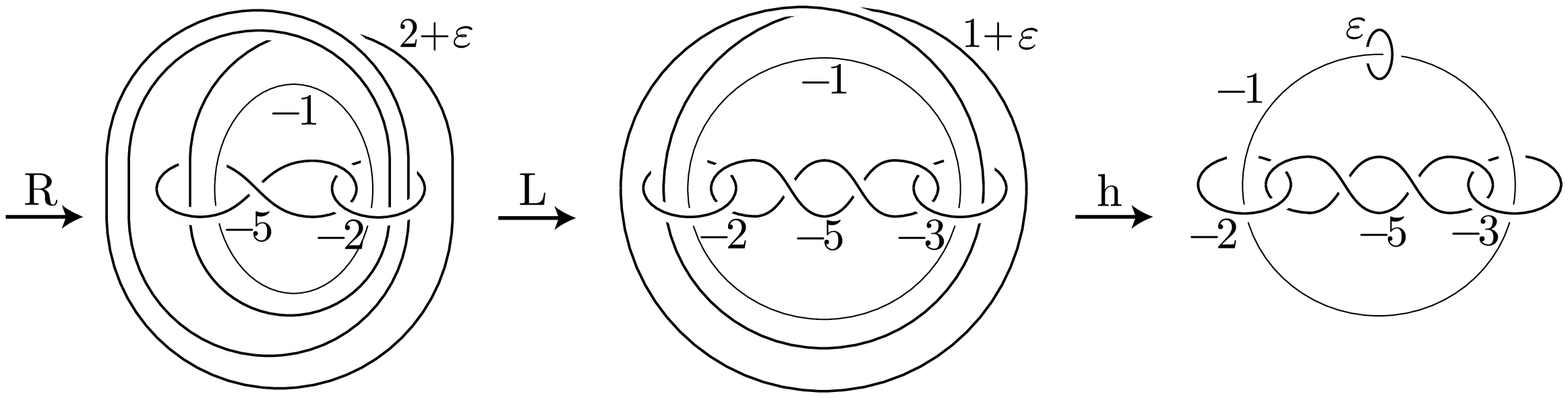}
%\caption{$(A_{2,3}; 7, 0)\cong L(25, 7)$ I}
\caption{$(A_{2,3}; 7, 0)$ is $L(25, 7)$ I}~\label{fig:A23_1}
\end{center}
\end{figure}

\begin{figure}[h]
\begin{center}
\includegraphics[scale=0.5]{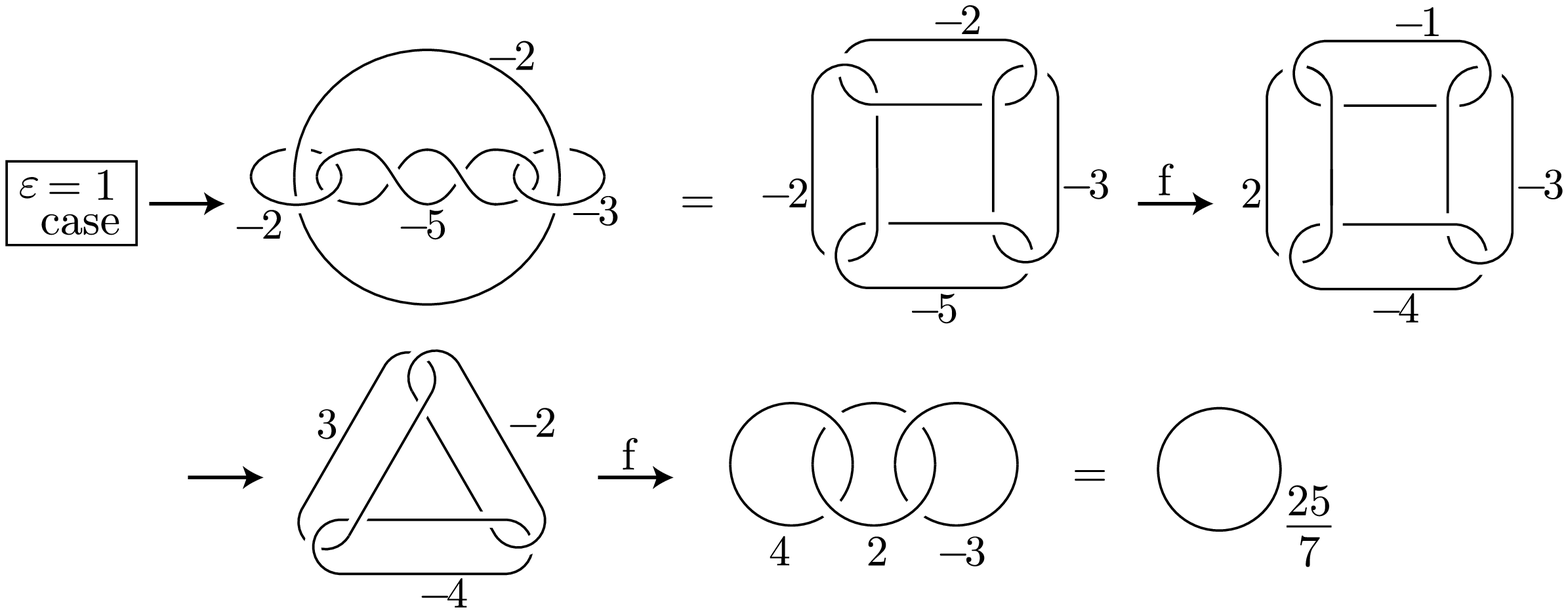}
%\caption{$(A_{2,3}; 7, 0)\cong L(25, 7)$ III}
\caption{$(A_{2,3}; 7, 0)$ is $L(25, 7)$ II}\label{fig:A23_2}
\end{center}
\end{figure}

\begin{figure}[h]
\begin{center}
\includegraphics[scale=0.5]{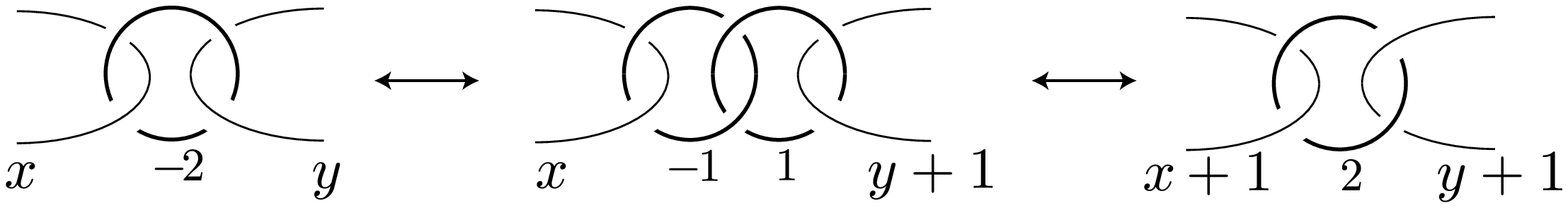}
\caption{Kirby move f}
\label{fig:KirbyF}
\end{center}
\end{figure}

%We also give a geometric proof of $(A_{2,3}; 7, 0)\cong L(25, 7)$
By Kirby calculus in Figure \ref{fig:A23_1} and \ref{fig:A23_2},
we have $(A_{2,3}; 7, 0)\cong L(25, 7)$.
For the moves in Figure~\ref{fig:A23_1}, see \cite{Yam3}.
By setting $\varepsilon=0$, we have $(A_{2,3}; 6, 0)\cong L(25, 9)$.
In the moves \lq\lq f\rq\rq in Figure~\ref{fig:A23_2},
we used the formula in Figure~\ref{fig:KirbyF}
on a $-2$-framed unknot.
In the second f move, the union of
two ($3$- and $-4$-framed) components
are regarded $x$, and get unlinked by the positive full twist.

%%%%%%%%%%%%%%%%%
%%%%%%%%%%%%%%%%%

%%%%%%%%%%%%%%%%%
%%%%%%%%%%%%%%%%%

%%%%%%%%%%%%%%%%%
%%%%%%%%%%%%%%%%%

%%%%%%%%%%%%%%%%%%%%%%%%%%%%%%%%%%%%%%%%%%
%%%                   Section 6
%%%%%%%%%%%%%%%%%%%%%%%%%%%%%%%%%%%%%%%%%%
\section{Proof of Theorem \ref{thm:MT2}\ 
(Lens space surgeries along $B_{p,q}$)}~\label{sec:Bpq}
In this section, we prove Theorem~\ref{thm:MT2} on lens space surgeries
along $B_{p,q}$.
%%%%%%%%%%%%%%%%%%%%%%%%%%%%%%%%%%%%%%%%%%
%%%                   SubSection 6-1
%%%%%%%%%%%%%%%%%%%%%%%%%%%%%%%%%%%%%%%%%%
\subsection{Alexander polynomial of $B_{p,q}$}~\label{ssec:BpqAlex}
\begin{figure}[h]
\begin{center}
\includegraphics[scale=0.4]{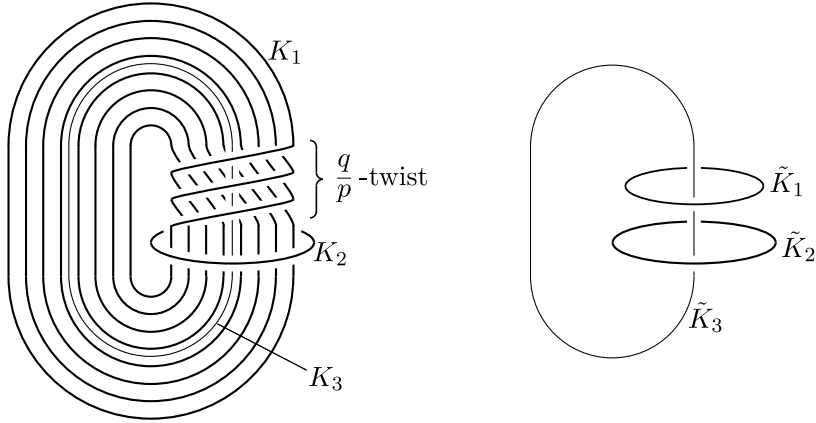}
\caption{Link $L$ (ex. $(p,q)=(8,3)$) and ${\tilde L}$}
\label{fig:LLt}
\end{center}
\end{figure}

The goal of this subsection is: 
%
%
%%%%%%%%%%%%%%%%%%%%%%%%%%%%%%%%%%%%
\begin{lemma}~\label{lem:BpqAlex}
The Alexander polynomial of the link $B_{p,q}$ is 
\[
{\it \Delta}_{B_{p,q}}(t, x)
\doteq \frac{t^{pq}x^p -1}{t^q x -1}.
\]
where $t$ (and $x$, respectively) is represented by
the meridian of the torus knot component $K_1$
(that of the unknotted component $K_2$).
\end{lemma}
%%%%%%%%%%%%%%%%%%%%%%%%%%%%%%%%%%%%
%
%
\medskip \noindent
{\bf Proof} \ 
We add the third component $K_3$ to $B_{p,q} = K_1 \cup K_2$
such that $H := K_2 \cup K_3$ is a Hopf link and 
$K_1$ is isotopic to a simple closed loop in the level torus 
under the identification of the complement of $H$ and 
$T^2 \times (-1,1)$.
We set $L := B_{p,q}\cup K_3$ and 
use the notations defined in Subsection~\ref{ssec:surgery}.
The complement $E_L$ of $L$ is homeomorphic to 
that $E_{{\tilde L}}$ of the connected sum
${\tilde L} = {\tilde K}_1 \cup {\tilde K}_2 \cup {\tilde K}_3$
of two Hopf links 
(${\tilde L}$ is \lq\lq$B_{0,1}\cup K_3$\rq\rq\, in a sense),
see Figure~\ref{fig:LLt}.
%
%
%%%%%%%%%%%%%%%%%
%%%%%%%%%%%%%%%%%
%
%

\medskip

Our strategy of the proof is as follows:
First we prove the Alexander polynomial of $L$
by studying the transformation of the homologies 
from the easier ${\tilde L}$ to $L$, and
use the Torres formula (Lemma~\ref{lem:Torres}) on 
the sublink $B_{p,q}$ of $L$. 

\medskip

We take a homeomorphism $f : E_{{\tilde L}} \rightarrow E_L$
that carries each regular neighborhood of ${\tilde K}_i$
to that of $K_i$ ($i = 1,2,3$).
We will use a tilde for the notation of each component of ${\tilde L}$:
We denote the meridian-longitude system of ${\tilde K}_i$
by ${\tilde m}_i, {\tilde l}_i$, and
its Alexander polynomial of ${\tilde L}$ by
${\it \Delta}_{{\tilde L}}({\tilde t}_1, {\tilde t}_2, {\tilde t}_3)$.
In $H_1(E_L)$ and $H_1(E_{{\tilde L}})$,
the longitudes are presented by meridians:
%%%%%
\begin{eqnarray}~\label{eq:rel1}
\begin{matrix}
[l_1]=[m_2]^p[m_3]^q, \hfill 
& [l_2]=[m_1]^p[m_3], \hfill 
& [l_3]=[m_1]^q[m_2], \hfill \\[5pt]
[{\tilde l}_1]=[{\tilde m}_3], \hfill 
& [{\tilde l}_2]=[{\tilde m}_3], \hfill 
& [{\tilde l}_3]=[{\tilde m}_1][{\tilde m}_2].\hfill
\end{matrix}
\end{eqnarray}
%%%%%
We note that
%%%%%
\begin{equation}~\label{eq:rel2}
[m_2]^p[l_2]^q=[m_2]^p([m_1]^p[m_3])^q=[m_1]^{pq}[l_1].
\end{equation}
%%%%%
We take integers $a, b$ satisfying $aq-bp=1$.
The isomorphism $f_{\ast} : H_1(E_{{\tilde L}})\to H_1(E_L)$
induced by $f$ satisfies:
%%%%%
\begin{eqnarray}~\label{eq:map1}
\begin{matrix}
f_{\ast}([{\tilde m}_1])=[m_1], \hfill 
& f_{\ast}([{\tilde l}_1])=[m_2]^p[l_2]^q=[m_1]^{pq}[l_1], \hfill\\
f_{\ast}([{\tilde m}_2])=[m_2]^a[l_2]^b,\hfill 
& f_{\ast}([{\tilde l}_2])=[m_2]^p[l_2]^q, \hfill\\
f_{\ast}([{\tilde m}_3])=[m_3]^q[l_3]^p, \hfill
& f_{\ast}([{\tilde l}_3])=[m_3]^b[l_3]^a. \hfill& 
\end{matrix}
\end{eqnarray}
%%%%%
%and, by (\ref{eq:rel2}), we have
%%%%%%
%\begin{eqnarray*}
%\begin{matrix}
%(f_{\ast})^{-1}([m_1])=[{\tilde m}_1], \hfill & 
%(f_{\ast})^{-1}([l_1])=[{\tilde m}_1]^{-pq}[{\tilde l}_1], \hfill\\
%%
%(f_{\ast})^{-1}([m_2])=[{\tilde m}_2]^q[{\tilde l}_2]^{-b},\hfill & 
%(f_{\ast})^{-1}([l_2])=[{\tilde m}_2]^{-p}[{\tilde l}_2]^a, \hfill\\
%%
%(f_{\ast})^{-1}([m_3])=[{\tilde m}_3]^a[{\tilde l}_3]^{-p}, \hfill & 
%(f_{\ast})^{-1}([l_3])=[{\tilde m}_3]^{-b}[{\tilde l}_3]^q. \hfill& 
%\end{matrix}
%\end{eqnarray*}
%%%%%
Thus, setting 
$t_i:=[m_i]$ and ${\tilde t}_i:=[{\tilde m}_i]\ (i=1, 2, 3)$, 
by (\ref{eq:rel1}), (\ref{eq:rel2}) and (\ref{eq:map1}), we have
%%%%%
\begin{equation}~\label{eq:rel3}
f_{\ast}({\tilde t}_1)=t_1,\quad
f_{\ast}({\tilde t}_2)=t_1^{pb}t_2^at_3^b,\quad
f_{\ast}({\tilde t}_3)=t_1^{pq}t_2^pt_3^q
\end{equation}
%%%%%
The Alexander polynomial of ${\tilde L}$ is known.
%%%%%
\begin{equation*}
{\it \Delta}_{{\tilde L}}({\tilde t}_1, {\tilde t}_2, {\tilde t}_3)
\doteq {\tilde t}_3-1.
\end{equation*}
%%%%%
Substituting (\ref{eq:rel3}) to this, 
we have the Alexander polynomial of $L$.
%%%%%
\begin{equation*}
{\it \Delta}_L(t_1, t_2, t_3)\doteq t_1^{pq}t_2^pt_3^q-1.
\end{equation*}
%%%%%
Hence, by Torres formula (Lemma~\ref{lem:Torres}),
we have the Alexander polynomial of $B_{p,q}$:
\[
{\it \Delta}_{B_{p,q}}(t, x)
\doteq \frac{{\it \Delta}_L(t, x, 1)}{t^q x -1}
\doteq \frac{t^{pq}x^p-1}{t^q x-1}.
\]
\qed

%%%%%%%%%%%%%%%%%%%%%%%%%%%%%%%%%%%%%%%%%%
%%%                   SubSection 6-2
%%%%%%%%%%%%%%%%%%%%%%%%%%%%%%%%%%%%%%%%%%
\subsection{Conditions from the first homology of $(B_{p,q}; r,0)$}~\label{ssec:hom2}
We set $M = (B_{p,q}; \alpha/\beta,0)$.
The arguments are parallel to Subsection~\ref{ssec:hom1}.
We let $E$ denote the complement of $B_{p,q}$, and regard
$M = E\cup V_1\cup V_2$, see the notations in Subsection~\ref{ssec:surgery}.
We set $M_1 = E\cup V_1$.

\medskip

We assume that $\gcd(p, \alpha)=1$,
which is equivalent to the condition 
for $H_1(M; \mathbb{Z})$ to be finite cyclic. 
Then we have
$H_1(M; \mathbb{Z})\cong \mathbb{Z}/p^2\beta \mathbb{Z}$.

\medskip

We determine the first homologies of $E$, $M_1$ and $M$,
and represent some elements by the generators.
First, $H_1(E)$ is a free abelian group of rank 2
generated by $[m_1]$ and $[m_2]$:
%%%%%
\begin{equation*}
H_1(E)\cong \langle [m_1], [m_2]\ \vert \ -\rangle
\cong \mathbb{Z}^2.
\end{equation*}
%%%%%
We have
%%%%%
\begin{equation}~\label{eq:EB}
[l_1] = [m_2]^p \quad\mbox{and}\quad [l_2] = [m_1]^p.
\end{equation}
%%%%%
We attach $V_1$ to $E$ to make $M_1$. 
We take integers $\gamma, \delta$ such that
$\alpha \delta-\beta \gamma = -1$, and fix the 
meridian-longitude system $m_1', l_1'$ of the solid torus $V_1$.
In $H_1(M_1)$, 
%%%%%
\begin{equation*}
[m'_1] = [m_1]^{\alpha}[l_1]^{\beta}=1, \quad
[l'_1]=[m_1]^{\gamma}[l_1]^{\delta}\quad
\end{equation*}
%%%%%
and (\ref{eq:EB}) also hold. Thus, we have
%%%%%
\begin{eqnarray*}
H_1(M_1) & \cong & \langle [m_1], [m_2]\ \vert \ 
[m_1]^{\alpha}[m_2]^{p\beta}=1\rangle \nonumber\\
& \cong & \langle T\ \vert \ -\rangle \cong \mathbb{Z},
\end{eqnarray*}
%%%%%
where $T=[m_1]^{\gamma'}[m_2]^{\delta'}$
by taking integers $\gamma', \delta'$ satisfying 
$\alpha \delta'-p\beta \gamma'=-1$.
By the relations above, we have
%%%%%%%%%%
\begin{eqnarray}~\label{eq:M2B}
\begin{matrix}
[m_1] 
& = & [m_1]^{-\alpha \delta'+p\beta \gamma'}\hfill \\
& = & ([m_1]^{\alpha}[m_2]^{p\beta})^{-\delta'}
([m_1]^{\gamma'}[m_2]^{\delta'})^{p\beta}=T^{p\beta},\\[5pt]
%%%
[m_2] 
& = & [m_2]^{-\alpha \delta'+p\beta \gamma'}\hfill\\
& = & ([m_1]^{\alpha}[m_2]^{p\beta})^{\gamma'}
([m_1]^{\gamma'}[m_2]^{\delta'})^{-\alpha}=T^{-\alpha}, \hfill\\[5pt]
%%%
[l'_1] 
& = & [m_1]^{\gamma}[l_1]^{\delta}\hfill\\
& = & [m_1]^{\gamma}[m_2]^{p\delta}
=T^{p\beta \gamma-p\alpha \delta}=T^p. \hfill
\end{matrix}
\end{eqnarray}
%%%%%%%%%%
Next, we attach $V_2$ to $M_1$ to make $M$.
By (\ref{eq:EB}) and (\ref{eq:M2B}) in $H_1(M)$, we have 
%%%%%
\begin{equation*}
[m'_2]=[l_2]=[m_1]^p=T^{p^2\beta}=1,
\quad
[l'_2]=[m_2]=T^{-\alpha},
\end{equation*}
and
%%%%%
\begin{equation*}
H_1(M) \ \cong \ 
\langle T\ \vert \ 
T^{p^2\beta}=1\rangle 
\ \cong \ 
\mathbb{Z}/p^2\beta \mathbb{Z}.
\end{equation*}
%%%%%

%%%%%%%%%%%%%%%%%%%%%%%%%%%%%%%%%%%%%%%%%%
%%%                   SubSection 6-3
%%%%%%%%%%%%%%%%%%%%%%%%%%%%%%%%%%%%%%%%%%
\subsection{Proof of Theorem \ref{thm:MT2}}~\label{ssec:ProofTh1.4}
%
%\begin{Proof}
By Surgery formula~II (Lemma~\ref{lem:surgery2})(2), 
Lemma~\ref{lem:BpqAlex}
and (\ref{eq:M2B}), we have
%%%%%
\begin{equation}~\label{eq:torM1B}
\tau(M_1)\doteq
{\it \Delta}_{B_{p,q}}(T^{p\beta}, T^{-\alpha})
(T^p-1)^{-1}
\doteq
\frac{T^{p(-\alpha+pq\beta)}-1}{(T^{-\alpha+pq\beta}-1)(T^p-1)}
\end{equation}
%%%%%
Let $\iota : \mathbb{Z}[H_1(M_1)]\to \mathbb{Z}[H_1(M)]$
be a ring homomorphism induced from the natural inclusion, 
$\psi_d : \mathbb{Z}[H_1(M)]\to \mathbb{Q}(\zeta_d)$
a ring homomorphism such that $\psi_d(T)=\zeta_d$,
and set $\rho_d:=\psi_d\circ \iota$,
where $d$ is a divisor of $p$. Note that $\gcd(d, \alpha)=1$.
Then by Surgery formula~I (Lemma~\ref{lem:surgery1}), 
(\ref{eq:M2B}) and (\ref{eq:torM1B}), we have
%%%%%
\begin{equation*}
\tau^{\psi_d}(M)\doteq
\tau^{\rho_d}(M_1)(\zeta_d^{\alpha}-1)^{-1}
\doteq
(\alpha-pq\beta)(\zeta_d^{\alpha}-1)^{-2}.
\end{equation*}
%%%%%

Suppose that $M$ is a lens space.
Then its Reidemeister torsion is equal to
that of a lens space (Example \ref{ex:lens}), i.e.\ 
there exist integers $i$ and $j$ such that
$\gcd(p, i)=\gcd(p, j)=1$ and
%%%%%
\begin{equation*}
(\alpha-pq\beta)(\zeta_d^{\alpha}-1)^{-2}
\doteq
(\zeta_d^i-1)^{-1}(\zeta_d^j-1)^{-1}.
\end{equation*}
%%%%%
By taking $d$-norm (see Subsection~\ref{ssec:norm}) 
of both hand-sides, 
we have a necessary condition
%%%%%
\begin{equation*}
\vert \alpha-pq\beta \vert = 1
\end{equation*}
%%%%%
by Lemma~\ref{lem:cyclotomic}.

\medskip

Conversely, suppose $\vert \alpha-pq\beta \vert = 1$ and 
set $\varepsilon=\alpha-pq\beta (=\pm 1)$.
We can prove that 
$(B_{p,q}; \alpha/\beta, 0)$ is homeomorphic to the lens space
by Kirby--Rolfsen moves \cite{Ro} as follows 
(cf. \cite{Mos}, see also \cite[Proposition~4.1]{IS}):
$M$ is the result of $(\alpha/\beta, 0, \infty)$-surgery 
along $L = B_{p,q} \cup K_3$. 
It is homeomorphic to
that of $(\varepsilon/\beta, -p/a, -a/p)$-surgery ($aq-bp=1$) along ${\tilde L}$
because of the identification 
between $E_L$ and $E_{{\tilde L}}$ in (\ref{eq:map1}),
and that of $(-p/a, -(a+\varepsilon p\beta)/p)$-surgery along the Hopf link 
by $(-\varepsilon \beta)$-twist on ${\tilde K}_1$.
Therefore we have $M\cong L(p^2\beta, \alpha).$
\qed
\par

%%%%%%%%%%%%%%%%%%%%%%%%%%%%%%%%%%%%%%%%%%
%%%                   SubSection 6-4
%%%%%%%%%%%%%%%%%%%%%%%%%%%%%%%%%%%%%%%%%%
\subsection{Other surgeries along $B_{p,q}$}~\label{ssec:otherB}
The link $B_{p,q}$ is included in a family of 
{\it Burde-Murasugi's links} \cite{BM}, 
whose complements admit structures of Seifert fiber spaces.
The complement of $B_{p,q}$ is homeomorphic 
to that of the $(3, 3)$-torus link (and to that of the link in Figure 10 also). %\ref{fig:LLt}
Thus we can determine all Dehn surgeries along $B_{p,q}$
by a result of the first author and M.~Shimozawa \cite{KS}
on Dehn surgeries along torus links.

\medskip

Let $(B_{p,q}; \alpha_1/\beta_1, \alpha_2/\beta_2)$
be the result of 
$(\alpha_1/\beta_1, \alpha_2/\beta_2)$-surgery
along $B_{p,q}$ where 
$\alpha_i$ and $\beta_i$ ($i=1, 2$) are coprime integers 
with $\beta_i\ge 1$.
%
%
%%%%%%%%%%%%%%%%%%%%%%%%%%%%%%%%%%%%
\begin{theorem}~\label{th:otherB}
We set
$\varepsilon_1=\alpha_1-\beta_1pq$, 
$\varepsilon_2=\alpha_2q-\beta_2p$, 
$P=\alpha_1\alpha_2-\beta_1\beta_2p^2$, and
$a$ and $b$ as coprime integers such that $aq-bp=1$.
Then the result of surgery
$M=(B_{p,q}; \alpha_1/\beta_1, \alpha_2/\beta_2)$
with $p\ge 2$ is as follows:
\begin{enumerate}
\item[(1)]
If $\varepsilon_1\varepsilon_2\ne 0$, then
$M$ is a Seifert fiber space over $S^2$
with at most three singular fibers
whose multiplicities are 
$|\varepsilon_1|$, $|\varepsilon_2|$ and $p$.

\item[(2)]
Suppose $\varepsilon_1\varepsilon_2\ne 0$.
Then $M$ is a lens space if and only if
$|\varepsilon_1|=1$ or $|\varepsilon_2|=1$.
The resulting lens spaces are
\[
L(P, -\beta_2 \varepsilon_1-\beta_1 \varepsilon_2q)\cong
L(P, -\alpha_1\beta_2-\alpha_2\beta_1q^2+2\beta_1\beta_2pq)
\]
in both cases.

\item[(3)]
\begin{enumerate}
\item[(i)]
If $\varepsilon_1=0$ (i.e.\ $\alpha_1/\beta_1=pq$), 
then $M$ is a connected sum of lens spaces
\[
L(\varepsilon_2, \alpha_2b-\beta_2a)\sharp L(p, -q).
\]

\item[(ii)]
If $\varepsilon_2=0$ (i.e.\ $\alpha_2/\beta_2=p/q$), 
then $M$ is a connected sum of lens spaces
\[
L(\varepsilon_1, -\beta_1)\sharp L(p, -q).
\]
\end{enumerate}
\end{enumerate}
Here we regard $L(\pm 1, Q)$ ($\cong S^3$)
and $L(0,Q)$($\cong S^1 \times S^2$) for any $Q$
as lens spaces.
\end{theorem}
%%%%%%%%%%%%%%%%%%%%%%%%%%%%%%%%%%%%
%
%
\medskip \noindent
{\bf Proof} \ 
Let $L_{p,q}$ be the link $L$ in the proof of Lemma 6.1. %\ref{lem:BpqAlex}
Then $L_{1,1}$ is the $(3, 3)$-torus link.
We set $L_{1,1}=\hat{L}=\hat{K}_1\cup \hat{K}_2\cup \hat{K}_3$.
Let $\hat{m}_i$ and $\hat{l}_i$ be a meridian and a longitude of 
$\hat{K}_i$ ($i=1, 2, 3$), respectively.

\medskip

\begin{figure}[h]
\begin{center}
\includegraphics[scale=0.5]{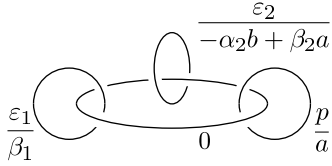}
\caption{Seifert fiber space $(B_{p,q}; \alpha_1/\beta_1, \alpha_2/\beta_2)$}
\label{fig:SeifertB}
\end{center}
\end{figure}

Let $h : E_L\to E_{{\hat L}}$ be a homeomorphism
inducing an isomorphism 
$h_{\ast} : H_1(E_L)\to H_1(E_{{\hat L}})$
such that
%%%%%
\begin{eqnarray*}
\begin{matrix}
h_{\ast}([m_1])=[{\hat m}_1], \hfill 
& h_{\ast}([l_1])=[{\hat m}_1]^{-pq+1}[{\hat l}_1], \hfill\\
h_{\ast}([m_2])=[{\hat m}_2]^{q-b}[{\hat l}_2]^{-b}, \hfill 
& h_{\ast}([l_2])=[{\hat m}_2]^{-p+a}[{\hat l}_2]^a, \hfill\\
h_{\ast}([m_3])=[{\hat m}_3]^a[{\hat l}_3]^{-p+a}, \hfill 
& h_{\ast}([l_3])=[{\hat m}_3]^{-b}[{\hat l}_3]^{q-b}.\hfill& 
\end{matrix}
\end{eqnarray*}
%%%%%
By the relations, we have
\begin{eqnarray*}
\begin{matrix}
h_{\ast}([m_1]^{\alpha_1}[l_1]^{\beta_1})
& = & [{\hat m}_1]^{\alpha_1-\beta_1(pq-1)}
[{\hat l}_1]^{\beta_1}, \hfill\\
h_{\ast}([m_2]^{\alpha_2}[l_2]^{\beta_2})
& = & [{\hat m}_2]^{\alpha_2(q-b)-\beta_2(p-a)}
[{\hat l}_2]^{-\alpha_2b+\beta_2a}, \hfill\\
h_{\ast}([m_3])
& = & [{\hat m}_3]^a[{\hat l}_3]^{-p+a},\hfill& 
\end{matrix}
\end{eqnarray*}
and hence $M$ is the result of 
$(\{\alpha_1-\beta_1(pq-1)\}/\beta_1, 
\{\alpha_2(q-b)-\beta_2(p-a)\}/(-\alpha_2b+\beta_2a), 
a/(-p+a))$-surgery along the $(3, 3)$-torus link.
By \cite[Theorem 3.3 (3)]{KS}, 
if $\varepsilon_1\varepsilon_2\ne 0$, then
$M$ is a Seifert fiber space over $S^2$
with at most three singular fibers
whose indices are
\[
\left(-1; \frac{\beta_1}{\varepsilon_1}, 
\frac{-\alpha_2b+\beta_2a}{\varepsilon_2}, 
\frac{-p+a}{p}\right)
=\left(0; \frac{\beta_1}{\varepsilon_1}, 
\frac{-\alpha_2b+\beta_2a}{\varepsilon_2}, 
\frac{a}{p}\right),
\]
see Figure~\ref{fig:SeifertB}.
%
%
%%%%%%%%%%%%%%%%%
%%%%%%%%%%%%%%%%%
%
%
Other cases are obtained from it.
\qed

%
%
%%%%%%%%%%%%%%%%%%%%%%%%%%%%%%%%%%%%
\begin{remark}
{\rm 
The convention on lens spaces in \cite{KS}
is different from that in the present paper.
Theorem 1.5 corresponds to %\ref{thm:MT2}
the case $\alpha_2=0$ and $\beta_2=1$
in Theorem \ref{th:otherB} (2), and
we can deduce Theorem 1.5 (2) %\ref{thm:MT2}
from Theorem \ref{th:otherB} (2) after some deformations.
}
\end{remark}
%%%%%%%%%%%%%%%%%%%%%%%%%%%%%%%%%%%%
%
%

%%%%%%%%%%%%%%%%%%%%%%%%%%%%%%%%%%%%%%%%%%
%%%                   Section 7
%%%%%%%%%%%%%%%%%%%%%%%%%%%%%%%%%%%%%%%%%%
\section{Lens space surgeries along $A_{m,n}$
other than $(r, 0)$-surgery}~\label{sec:otherA}
Contrasted to the case of $B_{p,q}$ in Subsection~\ref{ssec:otherB},
to determine all Dehn surgeries along $A_{m,n}$
is a hard problem. Related to Conjecture~\ref{conj:Amn}, 
in Subsection~\ref{ssec:mnr} and Subsection~\ref{ssec:7r},
we compute the Reidemeister torsions of $M$
under the case (1) $\alpha_1/\beta_1=mn$ and (2) $\alpha_1/\beta_1=7$ for $(m, n)=(2, 3)$, respectively.

%%%%%%%%%%%%%%%%%%%%%%%%%%%%%%%%%%%%%%%%%%
%%%                   SubSection 7-1
%%%%%%%%%%%%%%%%%%%%%%%%%%%%%%%%%%%%%%%%%%
\subsection{The case $(A_{m,n};mn, r)$}~\label{ssec:mnr}
For $A_{m,n}=K_1\cup K_2$,
we compute the Reidemeister torsions of $M=(A_{m,n};mn, r)$.
We set $r=\alpha/\beta$, where
$\alpha$ and $\beta$ are coprime integers with $\beta \ge 1$,
and $P= (m+n)^2\beta-mn \alpha$.
Let $d (\ge 2)$ be a divisor of $P$.

\medskip

We use the notations $E$, $M_1$ and $M$ defined in Subsection~\ref{ssec:surgery}.
By the similar way to Subsection~\ref{ssec:hom1}, 
we have two lemmas.
%
%
%%%%%%%%%%%%%%%%%%%%%%%%%%%%%%%%%%%%
\begin{lemma}~\label{lem:hom71}
We take integers $u$ and $v$ as $(m+n)u-mnv=1$, and
set $T=[m_1]^u[m_2]^v$ in $H_1(M_1)$.
Then 
\begin{enumerate}
\item[(1)] $T$ is a generator of $H_1(M_1)$. 
\item[(2)] It holds that $[m_1]=T^{m+n}$ and $[m_2]=T^{-mn}$.
\end{enumerate}
\end{lemma}
%%%%%%%%%%%%%%%%%%%%%%%%%%%%%%%%%%%%%
%
%
%%%%%%%%%%%%%%%%%%%%%%%%%%%%%%%%%%%%%
\begin{lemma}~\label{lem:core71}
Let $T$ be the generator of $H_1(M)$ which is induced
by that of $H_1(M_1)$ in Lemma~\ref{lem:hom71}.
Let $\psi_d : \mathbb{Z}[H_1(M)]\to \mathbb{Q}(\zeta_d)$
be a ring homomorphism such that $\psi_d(T)=\zeta_d$.
Then the core of $V_2$ (i.e., $[l_2]$) is mapped to
$\zeta_d^{(m+n)^2\delta-mn\gamma}$ by $\psi_d$,
where $\gamma$ and $\delta$ are integers
such that $\alpha \delta-\beta \gamma=-1$.
\end{lemma}
%%%%%%%%%%%%%%%%%%%%%%%%%%%%%%%%%%%%
%
%
%%%%%%%%%%%%%%%%%%%%%%%%%%%%%%%%%%%%
\begin{lemma}~\label{lem:gamma71}
%We take integers $\gamma$ and $\delta$ 
%as $\alpha \delta-\beta \gamma=-1$.
The integer $(m+n)^2\delta-mn \gamma$ is coprime to $P$.
\end{lemma}
%%%%%%%%%%%%%%%%%%%%%%%%%%%%%%%%%%%%
%
%
\medskip \noindent
{\bf Proof} \ 
We have
%%%%%
\begin{eqnarray}~\label{eq:matrix}
\begin{pmatrix}
\alpha & \beta \\
\gamma & \delta
\end{pmatrix}
\begin{pmatrix}
-mn \\
(m+n)^2
\end{pmatrix}
=
\begin{pmatrix}
(m+n)^2\beta-mn\alpha \\
(m+n)^2\delta-mn\gamma
\end{pmatrix}.
\end{eqnarray}
%%%%%
Since $\alpha \delta-\beta \gamma=-1$, the matrix is 
invertible over $\mathbb{Z}$, thus we have
\[
\gcd((m+n)^2\beta-mn\alpha,
(m+n)^2\delta-mn\gamma)
= \gcd((m+n)^2, mn)=1.
\]
\qed
%
%
%%%%%%%%%%%%%%%%%%%%%%%%%%%%%%%%%%%%%
\begin{theorem}~\label{thm:tor71}
\[
\tau^{\psi_d}(M)
\doteq
\frac{\zeta_d^{mn}-1}
{(\zeta_d^m-1)(\zeta_d^n-1)(\zeta_d^{(m+n)^2\delta-mn \gamma}-1)}.
\]
\end{theorem}
%%%%%%%%%%%%%%%%%%%%%%%%%%%%%%%%%%%%
%
%
\medskip \noindent
{\bf Proof} \ 
We use the surgery formula II (Lemma~\ref{lem:surgery2}).
By Lemma~\ref{lem:core71} and Lemma~\ref{lem:gamma71},
$\zeta_d^{(m+n)^2\delta-mn\gamma}$ is 
a primitive $d$-th root of unity.
By the Alexander polynomial of $A_{m,n}$ 
in Theorem~\ref{thm:AmnAlex}, Lemma~\ref{lem:hom71}, 
Lemma~\ref{lem:core71} and Lemma~\ref{lem:Amn(T)}, we have the result.
\qed

\bigskip

By the Franz lemma (Lemma~\ref{lem:Franz}),
$\tau^{\psi_d}(M)$ is same with
the Reidemeister torsion of a lens space if and only if
\[
m\equiv \pm 1, \quad
n\equiv \pm 1\quad \mbox{or}\quad
(m+n)^2\delta-mn \gamma\equiv \pm mn\quad
(\mathrm{mod}\! \ P).
\]
%
%
%%%%%%%%%%%%%%%%%%%%%%%%%%%%%%%%%%%%
\begin{lemma}~\label{cor:ByTor71}
\begin{enumerate}
\item[(i)]
If $\beta=1$, then $M$ has the same Reidemeister torsions
with a lens space $L(P, \pm m\overline{n})$
where $n\overline{n}\equiv 1\ (\mathrm{mod}\! \ P)$.
\item[(ii)]
$(m+n)^2\delta-mn \gamma \equiv \pm mn\ 
(\mathrm{mod}\! \ P)$
is equivalent to
$\beta \equiv \pm 1\ (\mathrm{mod}\! \ P)$.
\item[(iii)]
Suppose that $\alpha<0$, and
$M$ has the same Reidemeister torsions
with a lens space.
Then we have $\beta=1$.
\end{enumerate}
\end{lemma}
%%%%%%%%%%%%%%%%%%%%%%%%%%%%%%%%%%%%
%
%
\medskip \noindent
{\bf Proof} \ 
(i) If $\beta=1$, then we can take $\gamma=1, \delta=0$, thus
$(m+n)^2\delta-mn \gamma \equiv \pm mn\ 
(\mathrm{mod}\! \ P)$.

(ii) Suppose that $(m+n)^2\delta-mn \gamma \equiv \pm mn\ 
(\mathrm{mod}\! \ P)$.
Then by (\ref{eq:matrix}), it holds that
\[
\begin{pmatrix}
-mn \\
(m+n)^2
\end{pmatrix}
\equiv
\begin{pmatrix}
-\delta & \beta \\
\gamma & -\alpha
\end{pmatrix}
\begin{pmatrix}
0 \\
\pm mn
\end{pmatrix}
= 
\begin{pmatrix}
\pm \beta mn\\
\mp \alpha mn
\end{pmatrix}
\quad (\mathrm{mod}\! \ P),
\]
hence $\beta \equiv \pm 1\ (\mathrm{mod}\! \ P)$.
Conversely suppose that $\beta \equiv \pm 1\ (\mathrm{mod}\! \ P)$.
Then we can take
$\gamma \equiv \pm 1, \delta \equiv 0\ (\mathrm{mod}\! \ P)$, 
thus $(m+n)^2\delta-mn \gamma \equiv \pm mn\ 
(\mathrm{mod}\! \ P)$.

(iii) $P = (m+n)^2\beta-mn\alpha = (m+n)^2\beta + mn\vert \alpha \vert$.
By the assumption $2 \le m<n$, 
it holds that $m\not \equiv \pm 1$ and $n\not\equiv \pm 1\ (\mathrm{mod}\! \ P)$.
Since $0< \beta < P- mn \vert \alpha \vert$, we have  $\beta =1$.
\qed

\medskip

\begin{figure}[h]
\begin{center}
\includegraphics[scale=0.5]{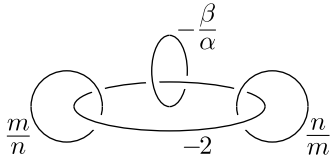}
\caption{Seifert manifold $(A_{m,n}; mn, \alpha/\beta)$}
\label{fig:SeifertA}
\end{center}
\end{figure}

Extending the Kirby calculus in \cite{Yam3},
one can prove the following:
%
%
%%%%%%%%%%%%%%%%%%%%%%%%%%%%%%%%%%%%%
\begin{lemma}[{An extension of \cite{Yam3}}]~\label{lem:AmnSeif}
The surgered manifold $(A_{m,n}; mn, \alpha/\beta)$ is a Seifert fiber space over $S^2$
with at most three singular fibers
whose indices are
\[
\left(
-2; 
\frac{n}{m},
\frac{m}{n},
-\frac{\alpha}{\beta}
\right),
\]
see Figure~\ref{fig:SeifertA}.
\end{lemma}
%%%%%%%%%%%%%%%%%%%%%%%%%%%%%%%%%%%%
%
%
%%%%%%%%%%%%%%%%%
%%%%%%%%%%%%%%%%%
%
%
As a corollary, we can prove some lens space surgeries.
%
%%%%%%%%%%%%%%%%%%%%%%%%%%%%%%%%%%%%%
\begin{corollary}~\label{lem:Amnlens}
The surgered manifold $(A_{m,n}; mn, \alpha/\beta)$ is a lens space if and only if $\beta =1$
(thus $P = (m+n)^2-\alpha mn$) and the lens space is $L(P, m\overline{n})$, where 
$n\overline{n} \equiv 1\ (\mathrm{mod}\! \ P)$.
\end{corollary}
%%%%%%%%%%%%%%%%%%%%%%%%%%%%%%%%%%%%
%
%
%
%
\medskip \noindent
{\bf Proof} \ (Using Lemma~\ref{lem:AmnSeif}) \
We apply \cite[Lemma 2.2]{KS} to the description in Figure~\ref{fig:SeifertA}
in the case $\beta =1$.
Then the manifold is $L(P,Q)$ with 
\[
P = (m+n)^2-\alpha mn, \quad
Q= mx + \{m(\alpha-2)-n\}y,
\]
where $x, y$ are coprime integers satisfying $nx+my=1$, and it holds that
%%%%%
\begin{eqnarray*}
nQ
& =& mnx + \{mn(\alpha-2)-n^2\}y \\
& \equiv & mnx + \{(m+n)^2-2mn-n^2\}y\ (\mathrm{mod}\ \! P) \\
& = & mnx + m^2y \\ 
& = & m(nx+my) \\
& =& m.
\end{eqnarray*}
%%%%%
\qed 

%%%%%%%%%%%%%%%%%%%%%%%%%%%%%%%%%%%%%%%%%%
%%%                   SubSection 7-2
%%%%%%%%%%%%%%%%%%%%%%%%%%%%%%%%%%%%%%%%%%
\subsection{The case $(A_{2,3};7, r)$}~\label{ssec:7r}
\begin{figure}[h]
\begin{center}
\includegraphics[scale=0.5]{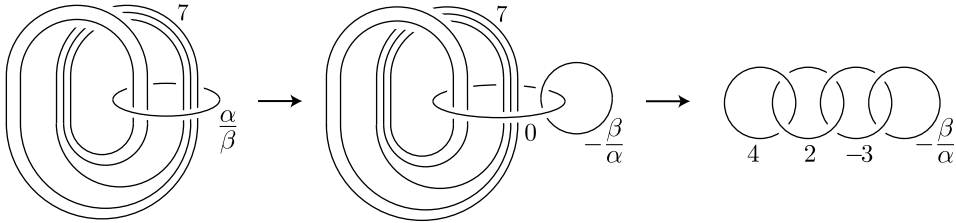}
\caption{Lens space $(A_{2,3}; 7, \alpha/\beta)$}
\label{fig:Lens}
\end{center}
\end{figure}

We set $(m, n)=(2, 3)$, $r=\alpha/\beta$,  where
$\alpha$ and $\beta$ are coprime integers
with $\beta \ge 1$ and $P= 25\beta-7\alpha$.
Let $d (\ge 2)$ be a divisor of $P$.
We compute the Reidemeister torsions of $M=(A_{2,3};7, r)$.
%
%
%%%%%%%%%%%%%%%%%%%%%%%%%%%%%%%%%%%%
\begin{theorem}\label{thm:tor72}
Let $\psi_d : \mathbb{Z}[H_1(M)]\to \mathbb{Q}(\zeta_d)$
be a ring homomorphism
which maps a generator of $H_1(M)$ to $\zeta_d$.
Then we have
$$\tau^{\psi_d}(M)
\doteq
(\zeta_d-1)^{-1}(\zeta_d^{7\gamma-25\delta}-1)^{-1}$$
where integers $\gamma$ and $\delta$ 
satisfy $\alpha \delta-\beta \gamma=-1$.

In particular, $M$ has the same Reidemeister torsions
with a lens space $L(P, \pm (7\gamma-25\delta))$.
\end{theorem}
%%%%%%%%%%%%%%%%%%%%%%%%%%%%%%%%%%%%
%
%
We can verify the lens space surgeries.
%
%
%%%%%%%%%%%%%%%%%%%%%%%%%%%%%%%%%%%%
\begin{lemma}~\label{thm:L25}
The surgered manifold $(A_{2,3};7, \alpha/\beta )$ is a lens space $L(25 \beta - 7\alpha, 2\alpha - 7 \beta)$.
\end{lemma}
%%%%%%%%%%%%%%%%%%%%%%%%%%%%%%%%%%%%
%
%
\medskip \noindent
{\bf Proof} \ 
We modify the Kirby calculus in 
Figure~\ref{fig:A23_1} and Figure~\ref{fig:A23_2}
as in Figure~\ref{fig:Lens}.
%
%
%%%%%%%%%%%%%%%%%
%%%%%%%%%%%%%%%%%
%
%
\qed

%%%%%%%%%%%%%%%%%%%%%%%%%%%%%%%%%%%%%%%%%%
%%%                   Section 8
%%%%%%%%%%%%%%%%%%%%%%%%%%%%%%%%%%%%%%%%%%
\section{Final remarks}~\label{sec:rem}
In this section, we give some remarks concerning our results.

%%%%%%%%%%%%%%%%%%%%%%%%%%%%%%%%%%%%%%%%%%
%%%                   SubSection 8-1 
%%%%%%%%%%%%%%%%%%%%%%%%%%%%%%%%%%%%%%%%%%
\subsection{Alternative proof for Theorem \ref{thm:MT1} (1)}
~\label{ssec:CGLS}
At the beginning of our study, the authors did not know whether
$(A_{m,n}; \emptyset, 0)$ is proved to be non-Seifert,
and our method by Reidemeister--Turaev torsion works 
without the knowledge. In fact, we are interested in 
the condition on the Alexander polynomial 
for a link to admit a lens space surgery.

In the preparation of this paper, K.~Motegi informed 
us that  $(A_{m,n}; \emptyset, 0)$ is proved to be non-Seifert.
We give an alternative rough proof of Theorem \ref{thm:MT1} (1)
by assuming $(A_{m,n}; \emptyset, 0)$ is non-Seifert
and by using the Cyclic Surgery Theorem \cite{CGLS}. 
\par

\medskip \noindent
{\bf Proof of Theorem \ref{thm:MT1} (1)} 
(Using the results in \cite{DMM1} and \cite{CGLS}) \ 
We set $M=(A_{m,n}; r, 0)$ ($r\in \mathbb{Q}$).
By the Cyclic Surgery Theorem \cite{CGLS}, 
only $r=mn-1$ or $mn+1$ can be a solution of Problem~\ref{prob1} other than $r=mn$.

\medskip
\noindent
{\bf Case 1} \ $r=mn-1$\\
This case does not occur by Lemma~\ref{lem:normR} (2).

\medskip
\noindent
{\bf Case 2} \ $r=mn+1$\\
By Lemma~\ref{lem:Rmn}, we have
\begin{equation*}
\tau^{\psi'_d}(M) \doteq 
(\xi^{m-1}-1)(\xi^{m+1}-1)(\xi-1)^{-2}(\xi^m-1)^{-2}
\end{equation*}
where $\xi$ is a primitive $d$-th root of unity ($d|m+n$ and $d\ge 2$).
If $M=(A_{m,n}; mn+1, 0)$ is a lens space, 
there exists integers $i$ and $j$ such that
$\gcd(i, m+n)=\gcd(j, m+n)=1$, as in (\ref{eq:torR}):
\begin{equation*}
(\xi^{m-1}-1)(\xi^{m+1}-1)(\xi^i-1)(\xi^j-1)
\doteq (\xi-1)^2(\xi^m-1)^2.
\end{equation*}
%Suppose $\gcd(m-1, m+n)\ge 2$ or $\gcd(m+1, m+n)\ge 2$.
%Then the left-hand side of the equation above is $0$ for some $d$.
%Hence we have $\gcd(m-1, m+n)=1$ and $\gcd(m+1, m+n)=1$.
By the similar way to the proof of Lemma~\ref{lem:gcd}, 
using Franz lemma \cite{Fz} (Lemma~\ref{lem:Franz}), 
we have $$\{\pm (m-1), \pm (m+1), \pm i, \pm j\ (\mathrm{mod}\ \! m+n)\}
=\{\pm 1, \pm 1, \pm m, \pm m\ (\mathrm{mod}\ \! m+n)\}$$
as multiple sets.
It has a unique solution $(m, n)=(2, 3)$ with $(i,j)=(1,3)$.
%and (as in Subsection~\ref{ssec:m=2}). 
Thus we have $r=7$.
\qed

%%%%%%%%%%%%%%%%%%%%%%%%%%%%%%%%%%%%%%%%%%
%%%                   SubSection 8-2 
%%%%%%%%%%%%%%%%%%%%%%%%%%%%%%%%%%%%%%%%%%
\subsection{Algebraic generalization}~\label{ssec:alg-gene}
Our main theorems can be extended to
the cases of 2-component links in homology 3-spheres
with the same Alexander polynomials as $A_{m,n}$ and $B_{p,q}$.

%
%
%%%%%%%%%%%%%%%%%%%%%%%%%%%%%%%%%%%%
\begin{theorem}~\label{thm:extMT1}
Let $L_{m,n}$ be a $2$-component link in a homology $3$-sphere
with the same Alexander polynomial as $A_{m,n}$
(Theorem~\ref{thm:AmnAlex}).
If $(L_{m,n}; r, 0)$ is a lens space, 
then we have $r=mn$, or $r=7$ for $(m, n)=(2, 3)$.
Moreover we have 
$(L_{m,n}; mn, 0)\cong L((m+n)^2, \pm m\overline{n})$
or $(L_{2,3}; 7, 0)\cong L(25, 7)$ respectively,
where $n\overline{n}\equiv 1\ (\mathrm{mod}\ \! (m+n)^2)$.
\end{theorem}
%%%%%%%%%%%%%%%%%%%%%%%%%%%%%%%%%%%%
%
%
%%%%%%%%%%%%%%%%%%%%%%%%%%%%%%%%%%%%
\begin{theorem}~\label{thm:extMT2}
Let $L_{p,q}'$ be a $2$-component link in a homology $3$-sphere
with the same Alexander polynomial as $B_{p,q}$
(Lemma~\ref{lem:BpqAlex}).
If $(L_{p,q}'; \alpha/\beta, 0)$ is a lens space, then we have $|\alpha-pq\beta|=1$.
Furthermore, in this case, we have
$(L_{p,q}'; \alpha/\beta, 0)\cong L(p^2\beta, \pm \alpha)$.
\end{theorem}
%%%%%%%%%%%%%%%%%%%%%%%%%%%%%%%%%%%%
%
%

We remark that the converses of the theorems above
do not always hold in general.

\bigskip \noindent
%%%%%%%%%%%%%%%%%%%%%%%%%%%%%%%%%%%%%%%%%%
{\bf Acknowledgement.}
The authors would like to thank to Professor Kimihiko Motegi 
and Professor Takayuki Morifuji for thier useful advices.
\bigskip

%%%%%%%%%%%%%%%%%%%%%%%%%%%%%%%%%%%%%%%%%%
%%%                   References
%%%%%%%%%%%%%%%%%%%%%%%%%%%%%%%%%%%%%%%%%%
%%%%%%%%%%%%%%%%%%%%%%%%%%%%%%%%%%%%%%%%%%
%%% END of Footnotesize

{\small
\par
Teruhisa KADOKAMI\par
Department of Mathematics,
East China Normal University,\par
Dongchuan-lu 500, Shanghai, 200241, China \par
{\tt kadokami2007@yahoo.co.jp}, 
{\tt mshj@math.ecnu.edu.cn}
%Department of Applied Mathematics,\par
%Dalian University of Technology\par
%Dalian-si, Liaoning-sheng, 116024, China \par
%{\tt kadokami@dlut.edu.cn}\par
%{\tt kadokami2007@yahoo.co.jp} \par
\medskip
\par
Yuichi YAMADA\par
Department of Mathematics,
The University of Electro-Communications \par
1-5-1,Chofugaoka, Chofu, Tokyo, 182-8585, JAPAN \par
{\tt yyyamada@sugaku.e-one.uec.ac.jp} \par
}

\end{document}